\documentclass[examplefnt,biber]{nowfnt} %

\usepackage[utf8]{inputenc}

\usepackage{tikz}

\usepackage{multirow}
\usepackage{bm}

\newcommand{\ie}{\emph{i.e.}}

\title{Optimization Frameworks \& Learning Algorithms for Discrete-Time Mean Field Games}

\subtitle{}

\maintitleauthorlist{
Xin Guo \\
UC Berkeley \\
xinguo@berkeley.edu
\and
Mathieu Lauri{\`e}re \\
NYU Shanghai \\
mathieu.lauriere@nyu.edu
}

\issuesetup
{ %
 copyrightowner={X.~Guo and M.~Lauri{\`e}re},
 copyrightyear = 2026,
 firstpage     = 1
 }

\usepackage{mwe}

\author[1]{Guo, Xin}
\author[2]{Lauri{\`e}re, Mathieu}

\affil[1]{UC Berkeley; xinguo@berkeley.edu}
\affil[2]{NYU Shanghai; mathieu.lauriere@nyu.edu}

\articledatabox{\strut}

\usepackage{amsmath,amssymb}
\usepackage{mathabx} %
\usepackage{algorithm}
\usepackage{algorithmicx}
\usepackage[noend]{algpseudocode}

\usepackage{imakeidx} %
\makeindex
\makeindex[name=notation,title={Index of notation}]
\newcommand{\boldindex}[1]{%
  \textbf{#1}\index{#1}%
}
\newcommand{\notationindex}[2]{\index[notation]{#1@#2}}

\newcommand{\states}{\mathcal{S}}
\newcommand{\actions}{\mathcal{A}}
\newcommand{\timesteps}{\mathcal{T}}

\newcommand{\strategies}{\Pi}
\newcommand{\state}{s}
\newcommand{\action}{a}
\newcommand{\supp}{\mathrm{supp}}

\DeclareMathOperator*{\argmax}{\arg\!\max}

\newcommand{\Poldist}{\operatorname{Pol}}
\newcommand{\OccupationMeasures}{\cD}
\newcommand{\softmax}{{\mathrm{softmax}}}
\newcommand{\indic}{{\mathbf{1}}}

\newcommand{\bstrategies}{\bsPi}

\newcommand{\bfa}{{\mathbf a}}

\newcommand{\bfs}{{\mathbf s}}

\newcommand{\cD}{{\mathcal D}}

\newcommand{\cT}{{\mathcal T}}

\newcommand{\cX}{{\mathcal X}}

\newcommand{\bbE}{{\mathbb E}}

\newcommand{\bbP}{{\mathbb P}}

\newcommand{\bbR}{{\mathbb R}}

\newcommand{\EE}{{\bbE}}

\newcommand{\PP}{{\bbP}}

\newcommand{\RR}{{\bbR}}

\newcommand{\bspi}{{\boldsymbol{\pi}}}

\newcommand{\bsPi}{{\boldsymbol{\Pi}}}

\newtheorem{Assumption}{Assumption}

\newcommand{\BR}{\mathrm{BR}}
\newcommand{\MF}{\mathrm{MF}}
\newcommand{\EMF}{\mathrm{EMF}}
\DeclareMathOperator{\Expl}{Expl}

\usepackage[normalem]{ulem}

\usepackage{subcaption}

\usepackage{booktabs}
\usepackage{tabularx}
\usepackage{array}
\input{arxiv-frontmatter.tex}

\begin{document}

\makeabstracttitle

\chapter*{Preface}
\addcontentsline{toc}{chapter}{Preface}
\markboth{\sffamily\slshape Preface}{\sffamily\slshape Preface}
The theory of mean field games (MFGs), pioneered by the seminal works of \citet{lasry2007mean} and \citet{huang2006large}, provides a tractable way of analyzing approximate Nash equilibria (NEs) for stochastic games that are otherwise difficult to solve. Since its inception, the literature on MFGs has grown in both theory and applications.
There are several excellent references on the subject, particularly for continuous time and space models, including the comprehensive two-volume textbook by \citet{carmona2018probabilistic}.

Meanwhile, building on works such as~\citep{cardaliaguet2017learning,guo2019learning}, learning algorithms for MFGs have developed into an active research direction, spanning the design and analysis of efficient algorithms as well as numerical experiments. This line of work has also reached the broader engineering and operations research communities.

As the literature continues to grow rapidly, it has become increasingly challenging for students and even experienced researchers to navigate both MFG theory and learning methods. Consequently, there is a growing demand for introductory notes that are easily accessible to a broader audience.

The present monograph is based on the tutorial given by the two authors at the 2024 AAAI Conference in Vancouver. It aims to provide a gentle introduction to MFGs in a discrete setting. This relatively simple framework allows us to present the most essential elements of MFG theory and its optimization and algorithmic learning aspects without delving into the mathematically demanding treatments required in the continuous-time setting.

\paragraph{Overview.}
This monograph covers discrete-time MFGs with illustrative examples and discusses the concept of mean field Nash equilibrium, adapted from the pioneering work of \citet{nash1950equilibrium}. The text addresses the fixed-point solution approach for MFGs, rooted in the classical theory of \citet{kakutani1941generalization}. The presentation is inspired by the classical optimization perspective for Markov decision processes (MDPs). The primal-dual formulation of MFGs is discussed, drawing an analogy to the linear programming formulation of MDPs, with early work on this viewpoint including~\citep{manne1960linear}. This perspective naturally leads to optimization-based algorithms for MFGs, complementing fixed-point iterations and learning dynamics such as fictitious play and online mirror descent.

The monograph is organized accordingly: it reviews MDPs and their linear-programming formulation, introduces the mean field framework with its equilibrium theory, presents fixed-point and optimization-based algorithms, and concludes with applications such as crowd motion, traffic routing, epidemics, and the computation of multiple equilibria. Each chapter ends with bibliographic notes for readers interested in a deeper theoretical understanding or exploring relevant research developments and related topics.

This monograph is intended for advanced undergraduate and graduate students, as well as researchers, in applied mathematics, engineering, operations research, economics, and related areas. A basic mathematical background, including probability, linear algebra, convex analysis, and linear programming, is assumed. The chapters are written to be read sequentially, but they can also be used selectively as a compact reference on discrete-time MFGs and learning algorithms.

\chapter{Markov Decision Processes and Markov Games}
\label{chap:mdp-markov-games}

\section{Markov decision processes}

This chapter introduces Markov decision processes and finite-population Markov games, which will serve as foundations for the theory of mean field games studied in the next chapters.

\subsection{Common notations}

For a finite set $S$, we denote by $|S|$\notationindex{cardinality}{$\lvert S\rvert$: cardinality of $S$} its cardinality, and by $\Delta(S)$\notationindex{DeltaS}{$\Delta(S)$: probability distributions on $S$} the set of probability distributions on $S$.
$\mathbb{P}_{x}$ denotes the probability distribution of a random variable $x$. We use $\PP$ and $\EE$ to denote respectively a probability and an expectation.
For a finite set $\mathcal{X}$, once an order of its elements has been fixed, we use $\mathbf{1}_{\mathcal{X}}$ to denote the all-ones vector in $\mathbb{R}^{|\mathcal{X}|}$. For $x \in \cX$, we use the notation $\indic_{x}$ to denote the function from $\cX$ to $\{0,1\}$, which takes the value $1$ if and only if the input is $x$. For $x,y \in \cX$, we will sometimes use the notation $\indic_{x=y} = \indic_x(y)$.

\subsection{Basic setup for Markov decision processes}
\label{sec:basic-MDP}

We present the setup of Markov decision processes. For simplicity, we shall focus on the framework of finite state and action spaces.

A \boldindex{Markov decision process} (or \boldindex{MDP} for short) is defined as a tuple $(\states,\actions,T,\mu_0,P,r,\Pi)$,
where
\begin{itemize}
	\item $\states$\notationindex{Sstate}{$\states$: state space} is the \boldindex{state space}, assumed non-empty and finite; a typical state is denoted by $\state$\notationindex{sstate}{$\state$: state};
	\item $\actions$\notationindex{Aaction}{$\actions$: action space} is the \boldindex{action space}, also assumed non-empty and finite; a typical action is denoted by $\action$\notationindex{aaction}{$\action$: action};
	\item $T$\notationindex{timehorizon}{$T$: time horizon} is a \boldindex{time horizon}, which could be finite or infinite (i.e., $T=+\infty$); the set of times is denoted by $\timesteps = \{0,1,\dots,T\}$\notationindex{Tcal}{$\timesteps$: set of times}, which should be understood as the set of all non-negative integers if $T=+\infty$; one generic time is denoted by $t$\notationindex{time}{$t$: time index};
	\item $\mu_0$\notationindex{mu0}{$\mu_0$: initial state distribution} is the initial distribution for the state;
	\item $P$\notationindex{Ptransition}{$P$: transition kernel} is the \boldindex{transition kernel} for the state dynamics:
	      $P: \timesteps \times \states \times \actions \to \Delta(\states),$
	      where $P_t(\cdot|\state, \action)$ denotes the probability distribution of the next state, when the current state is $s$, the action taken is $a$ and the current time is $t$;
	\item $r$\notationindex{rreward}{$r$: reward function} is the \boldindex{reward function}:
	      $r: \timesteps \times \states \times  \actions \to \RR,$
	      and $r_t(\state, \action)$ is the reward obtained when the state is $s$ and action $a$ is used at time $t$;
	\item $\strategies$\notationindex{Pi}{$\strategies$: policy space} is the set of policies\index{policy}, which are the sequences $ \pi = (\pi_t)_{t \in \timesteps}$\notationindex{pi}{$\pi$: policy} where for each $t \in \timesteps$, $\pi_t: \states \to \Delta(\actions)$. For $\pi \in \strategies$, $\pi_t(\cdot|s)$ is the distribution of the next action if the current time is $t$ and the current state is $s$. Note that this set contains the set of deterministic (pure) policies, i.e., policies $\pi$ for which there exists $a_t(s)$ such that $\pi_t(a|s) =  \indic_{a = a_t(s)}$.
\end{itemize}
For this MDP, the \boldindex{Markov decision problem} is to maximize, over all $\pi \in \strategies$, the total expected reward defined as
\begin{equation}
	\label{eq:def-MDP-J}
	J(\pi) = \EE\left[ \sum_{t \in \timesteps} r_t(s_t, a_t)\right],
\end{equation}
\notationindex{J}{$J$: total expected reward}
subject to the dynamics:
\[
	\begin{cases}
		s_0 \sim \mu_0
		\\
		a_t \sim \pi_t(\cdot|s_t), \quad t \in \timesteps,
		\\
		s_{t+1} \sim P_t(\cdot|s_t,a_t), \quad t \in \timesteps\backslash\{T\}.
	\end{cases}
\]
We say that $\pi \in \Pi$ is an \boldindex{$\epsilon$-suboptimal policy} for the MDP if:
\[
	J(\pi) \ge \sup_{\pi'} J(\pi') - \epsilon.
\]
A policy $\pi^*$ is an \boldindex{optimal policy} if and only if it is a $0$-suboptimal policy.

In the above definition, we assume that the expected sum of all rewards is well defined. If $T<+\infty$, it is sufficient to assume that $r$ is bounded. If $T=+\infty$, we usually take $r$ as a discounted reward; see Remark~\ref{rem:setting-statio}. In the finite horizon setting, the last action, $a_{T}$, does not have any consequence on the system so it is often assumed that the final reward $r_T$ depends only on the final state $s_{T}$.

In our setting (of finite state-action space), if $T<+\infty$, existence of an optimal policy is guaranteed and in fact there is an optimal policy that is \boldindex{pure}, i.e., deterministic. In other words, there exists $\pi^* \in \strategies$ such that for every $(t,s) \in \cT \times \states$, the support of the distribution $\pi^*_t(\cdot|s)$ is a singleton, and furthermore $J(\pi^*) = \sup_{\pi \in \strategies} J(\pi)$.
When the state-action space is infinite (and possibly continuous) or when the problem is with an infinite horizon, assumptions on the reward function and the transition kernel are needed. We refer to the bibliographic notes at the end of the chapter for more details.

\begin{remark}[Infinite horizon discounted and stationary setting]
	\label{rem:setting-statio}
	In the infinite horizon setting, i.e., $T = +\infty$, it is usually assumed that the reward is of the form $r_t(s,a) = \gamma^t r(s,a)$, where $\gamma \in (0,1)$ is a discount factor to ensure the total expected reward $J(\pi) = \EE\left[ \sum_{t=0}^{+\infty} \gamma^t r(s_t, a_t)\right]$ is well defined. If $P$ is independent of time, then one can show, under suitable conditions, that there exists an optimal deterministic (pure) policy which is also independent of time, i.e., it is a \boldindex{stationary policy}. This will considerably simplify the search for an optimal solution.
\end{remark}

\subsection{Value functions and Bellman equations}
\label{sec:mdp-value-bellman}
One of the fundamental principles in the theory of MDP is the \boldindex{dynamic programming principle} (\boldindex{DPP}), a.k.a. the time consistency property.
Here, we review the basic form of DPP for MDP. 
We consider separately the finite horizon and the infinite horizon discounted settings.

\paragraph{Finite horizon setting. }
The \boldindex{value function} associated to a given policy $\pi$ is defined as:
\begin{equation}
	\label{eq:V-pi-def}
	V^\pi_t(s) = \EE\left[ \sum_{t'=t}^{T} r_{t'}(s_{t'}, a_{t'})\right], \qquad t \in \timesteps, s \in \states,
\end{equation}
\notationindex{V}{$V^\pi$: value function under policy $\pi$}
where the trajectory is generated from $s_t=s$ using policy $\pi$: for $t'\in\{t,t+1,\dots,T\}$, $a_{t'} \sim \pi_{t'}(\cdot|s_{t'})$, and, for $t'<T$, $s_{t'+1} \sim P_{t'}(\cdot|s_{t'},a_{t'})$.

The \boldindex{optimal value function} is defined as:
\[
	V^*_t(s) = \sup_{\pi} V^\pi_t(s), \qquad t \in \timesteps, s \in \states.
\]
The total expected reward is the average value when starting with a state distributed according to $\mu_0$:
\[
	J(\pi) = \EE_{s \sim \mu_0}[V^\pi_0(s)],
\]
and hence:
\[
	\sup_{\pi} J(\pi) = \EE_{s \sim \mu_0}[V^*_0(s)].
\]

We define analogously the state-action value function by conditioning on the first state-action pair instead of the first state only. Namely, the \boldindex{state-action value function}, also called \boldindex{Q-function}, associated to a given policy $\pi$ is defined as:
\begin{equation}
	\label{eq:Q-pi-def}
	Q^\pi_t(s,a) = \EE\left[ \sum_{t'=t}^{T} r_{t'}(s_{t'}, a_{t'})\right], \qquad t \in \timesteps, s \in \states, a \in \actions,
\end{equation}
\notationindex{Q}{$Q^\pi$: state-action value function under policy $\pi$}
where $s_{t'}$ starts from $s_t = s$ at time $t$, the first action is $a_t = a$ and then the trajectory is generated using policy $\pi$. The \boldindex{optimal state-action value function} (or \boldindex{optimal Q-function}) is defined as:
\begin{equation}
	\label{eq:MDP-Qstar}
	Q^*_t(s,a) = \sup_{\pi} Q^\pi_t(s,a), \qquad t \in \timesteps, s \in \states, a \in \actions.
\end{equation}
The following relations hold between the $V$-functions and $Q$-functions:
\[
	V^\pi_t(s) = \EE_{a \sim \pi_t(\cdot|s)}[ Q^\pi_t(s,a)], \qquad V^*_t(s) = \max_a Q^*_t(s,a).
\]

The $V$-functions and the $Q$-functions satisfy dynamic programming principles, which leads to a \boldindex{Bellman equation} (also called \boldindex{dynamic programming equation}). For a given policy $\pi$, we have:
\[
	\begin{cases}
		V^\pi_{T}(s) = \EE_{a \sim \pi_T(\cdot|s)}[ r_T(s,a) ], \quad s \in \states
		\\
		V^\pi_{t}(s) = \EE_{\substack{a \sim \pi_t(\cdot|s), \\ s' \sim P_t(\cdot|s,a)}}[ r_t(s,a) + V^\pi_{t+1}(s')],
		\\
		\qquad\qquad\quad s \in \states, t = T-1, \dots, 0.
	\end{cases}
\]
and
\begin{equation}
	\label{eq:Q-pi-Bellman}
	\begin{cases}
		Q^\pi_{T}(s,a) =  r_T(s,a), \quad s \in \states, a \in \actions,
		\\
		Q^\pi_{t}(s,a) = \EE_{\substack{s' \sim P_t(\cdot|s,a), \\ a' \sim \pi_{t+1}(\cdot|s')}}[ r_t(s,a) + Q^\pi_{t+1}(s',a')],
		\\
		\qquad\qquad\quad s \in \states, a \in \actions, t = T-1, \dots, 0.
	\end{cases}
\end{equation}
These recursive equations  are useful  for policy evaluation in the policy iteration algorithm.

Correspondingly, the optimal value functions $V^*$ and $Q^*$ satisfy the following Bellman equations:
\[
	\begin{cases}
		\displaystyle
		V^*_{T}(s) = \max_{a \in \actions} r_T(s,a), \quad s \in \states
		\\
		\displaystyle
		V^*_{t}(s) = \max_{a \in \actions} \left\{ r_t(s,a) + \EE_{s' \sim P_t(\cdot|s,a)}[ V^*_{t+1}(s')] \right\},
		\\
		\qquad\qquad\quad s \in \states, t = T-1, \dots, 1, 0
	\end{cases}
\]
and
\begin{equation}
	\label{eq:Qstar-Bellman}
	\begin{cases}
		\displaystyle
		Q^*_{T}(s,a) =  r_T(s,a), \quad s \in \states, a \in \actions,
		\\
		\displaystyle
		Q^*_{t}(s,a) = \EE_{s' \sim P_t(\cdot|s,a)}[ r_t(s,a) + \max_{a' \in \actions} Q^*_{t+1}(s',a')], \\
		\qquad\qquad\quad s \in \states, a \in \actions, t = T-1, \dots, 1, 0.
	\end{cases}
\end{equation}
These equations are building blocks for algorithms such as the value iteration algorithm.

	Furthermore, given the optimal $Q$-function, an optimal policy can be constructed by selecting optimal actions at every time and state. More precisely, every policy $\pi$ satisfying
	\[
		\supp(\pi_t(\cdot|s)) \subseteq \argmax Q^*_t(s,\cdot),
		\qquad t \in \timesteps,\ s \in \states,
	\]
	is optimal. Such a policy is called a \boldindex{greedy policy} with respect to $Q^*$. Given a Q-function $Q: \timesteps \times \states \times \actions \to \RR$, we denote the set of greedy policies with respect to $Q$ by
	\begin{equation}
		\label{eq:greedy-policies}
		\mathrm{Greedy}(Q) := \{\pi \in \strategies \,:\, \forall t,s,a, \, \pi_t(a|s)>0 \Rightarrow a \in \argmax Q_t(s,\cdot)\}.
	\end{equation}

\paragraph{Infinite horizon discounted setting. }
In the infinite horizon discounted and stationary setting discussed in Remark~\ref{rem:setting-statio}, the definition of the optimal value function becomes
$$
	V^*_t(s) = \max_{\pi} \EE\left[ \sum_{t'=t}^{+\infty} \gamma^{t'}r(s_{t'}, a_{t'})\right],
$$
which depends on time due to the presence of the discount rate. Nevertheless, one can introduce the rescaled value function
$$
	\tilde{V}^*_t(s) = \gamma^{-t} V^*_t(s) = \max_{\pi} \EE\left[ \sum_{t'=t}^{+\infty} \gamma^{t'-t}r(s_{t'}, a_{t'})\right],
$$
and show that this function is independent of time. In fact, it satisfies the stationary Bellman equation:
\[
	\tilde{V}^*(s) = \max_{a \in \actions} \left\{ r(s,a) + \gamma \EE_{s' \sim P(\cdot|s,a)}[ \tilde{V}^*(s')] \right\}, \quad s \in \states.
\]
The corresponding optimal state-action value function is defined by
\[
	\tilde{Q}^*(s,a) = r(s,a) + \gamma \EE_{s' \sim P(\cdot|s,a)}[\tilde{V}^*(s')],
\]
and satisfies the stationary Bellman equation:
\[
	\tilde{Q}^*(s,a) = \EE_{s' \sim P(\cdot|s,a)}[ r(s,a) + \gamma \max_{a' \in \actions} \tilde{Q}^*(s',a')],     \quad s \in \states, a \in \actions.
\]
Policy-dependent value functions can also be introduced and, for stationary policies, the rescaled value functions also satisfy stationary Bellman equations, that we omit for the sake of brevity.

\subsection{Linear program formulation of MDPs}
\label{sec:LP-MDP}

In this section, we focus on the finite horizon setting, $T<+\infty$.
The classical \boldindex{linear program formulation} for MDP is built on the concept of an \boldindex{occupation measure}. 

To start,  denote  $\OccupationMeasures = \Delta(\states\times\actions)^\timesteps$\notationindex{Doccupation}{$\OccupationMeasures$: occupation-measure space} for sequences of state-action distributions. Then for a given policy $\pi \in \strategies$, the  occupation measure $d^\pi$\notationindex{dpi}{$d^\pi$: occupation measure induced by $\pi$} $\in \OccupationMeasures$ following this policy is  defined as
\begin{equation}\label{occupationmeasure}
	d^\pi_t(s,a)=\mathbb{P}(s_t=s,a_t=a)
\end{equation}
with $s_0\sim\mu_0$, $a_t\sim \pi_t(\cdot|s_t)$ for $t \in \timesteps$, and $s_{t+1}\sim P_t(\cdot|s_t,a_t)$ for $t \in \timesteps\backslash\{T\}$.

Conversely, given an occupation measure, one can retrieve a policy that gives rise to this occupation measure. Note that for states that have zero mass in the occupancy measure, any action could be used without changing the occupancy measure.
This yields a set-valued mapping, denoted by $\Poldist: \OccupationMeasures \to 2^{\strategies}$\notationindex{Poldist}{$\Poldist$: policies induced by an occupation measure}, defined by: for $(d_t)_{t\in\mathcal{T}} \in \OccupationMeasures$,
$\Poldist(d)$ is the set of policies $\pi$ satisfying:
\begin{equation}\label{eq:relation-d-to-pi}
	\pi_t(a|s)=\frac{d_t(s,a)}{\sum_{a'\in\mathcal{A}}d_t(s,a')}
\end{equation}
if $\sum_{a'\in\mathcal{A}}d_t(s,a')>0$, and otherwise $\pi_t(\cdot|s)$ is an arbitrary probability vector in $\Delta(\mathcal{A})$ if  $\sum_{a'\in\mathcal{A}}d_t(s,a')=0$.

Based on the notion of occupation measure, the MDP introduced above can be reformulated as a linear program as follows. For brevity, we use $\mathbb{P}^{\pi}(\cdot)$  for the state and/or action distribution generated by 
$s_0\sim\mu_0$, $a_t\sim\pi_t(\cdot|s_t)$ for $t=0,\dots,T$, and $s_{t+1}\sim P_t(\cdot|s_t,a_t)$ for $t=0,\dots,T-1$.

\begin{proposition}\label{prop:MDP_LP}
	Suppose that $\pi \in \strategies$ is an $\epsilon$-suboptimal policy for the MDP. Define $\mu_t(s):=\mathbb{P}^{\pi}(s_t=s)$ and $d_t(s,a):=\mu_t(s)\pi_t(a|s)$ for $s\in\mathcal{S}$ and $a\in\mathcal{A}$. Then $\pi\in\Poldist(d)$ and $d \in \OccupationMeasures$ is a feasible $\epsilon$-suboptimal solution to the following linear program:
	\begin{align}
		\max_{d \in \OccupationMeasures}\quad & \sum_{t\in\mathcal{T}}\sum_{s\in\mathcal{S}}\sum_{a\in\mathcal{A}} d_t(s,a)r_t(s,a) \notag                                                                               \\
\text{subj. to} \quad & \sum_{s\in\mathcal{S}}\sum_{a\in\mathcal{A}}d_t(s,a)P_t(s'|s,a)=\sum_{a\in\mathcal{A}}d_{t+1}(s',a),\notag\\
                     & \qquad\forall s'\in\mathcal{S},\ t\in\mathcal{T}\backslash\{T\},\notag\\
		                                      & \sum_{a\in\mathcal{A}}d_0(s,a)=\mu_0(s), \quad \forall s\in\mathcal{S},\notag                                                                                            \\
		                                      & d_t(s,a)\geq 0,\quad \forall s\in\mathcal{S},a\in\mathcal{A},t\in\mathcal{T}.\notag
	\end{align}
	Conversely, suppose that $d \in \OccupationMeasures$ is a feasible $\epsilon$-suboptimal solution to the above linear program. Then any $\pi\in\Poldist(d)$
	is an $\epsilon$-suboptimal policy for the MDP. In both directions, the linear-program objective evaluated at $d$ is equal to $J(\pi)$.
\end{proposition}

\section{The dual linear program and Bellman inequalities}
The occupation-measure formulation also has a useful dual interpretation
in terms of value functions. This interpretation provides the bridge
between the linear programming formulation of MDPs and the primal--dual
formulation of mean field games developed later.

For duality purposes, it is convenient to rewrite the flow constraints in
the following equivalent form:
\begin{equation}
  \label{eq:mdp-primal-flow-form}
  \begin{aligned}
    \underset{d}{\operatorname{maximize}}\quad
    & \sum_{t=0}^{T}\sum_{s\in\states}\sum_{a\in\actions}d_t(s,a)r_t(s,a)\\
    \operatorname{subject\ to}\quad
    & \sum_{a\in\actions}d_0(s,a)=\mu_0(s),\quad s\in\states,\\
    & \sum_{a\in\actions}d_t(s,a)=\sum_{s'\in\states}\sum_{a'\in\actions}
      d_{t-1}(s',a')P_{t-1}(s\mid s',a'),\\
    & \qquad s\in\states,\quad t=1,\ldots,T,\\
    & d_t(s,a)\geq 0,\\[-0.25ex]
    & \qquad s\in\states,\quad a\in\actions,\quad t=0,\ldots,T.
  \end{aligned}
\end{equation}
The variables in the dual program are unrestricted real numbers
$v_t(s)$, indexed by time and state. The dual linear program is
\begin{equation}
  \label{eq:mdp-dual-lp}
  \begin{aligned}
    \underset{v}{\operatorname{minimize}}\quad
    & \sum_{s\in\states}\mu_0(s)v_0(s)\\
    \operatorname{subject\ to}\quad
    & v_t(s)\geq r_t(s,a)\\[-0.25ex]
    & \qquad+\sum_{s'\in\states}P_t(s'\mid s,a)v_{t+1}(s'),\\
    & \qquad s\in\states,\quad a\in\actions,\quad t=0,\ldots,T-1,\\
    & v_T(s)\geq r_T(s,a),\quad s\in\states,\quad a\in\actions.
  \end{aligned}
\end{equation}
Thus, a dual-feasible vector $v$ is a supersolution of the Bellman
optimality equations. It gives an upper bound on the reward obtainable
from every time--state pair, and the dual objective minimizes the
expected initial upper bound under $\mu_0$.

For a dual-feasible vector $v$, define its Bellman slack by
\begin{equation}
  \label{eq:mdp-bellman-slack}
  z_t^v(s,a):=
  \begin{cases}
    \begin{aligned}[t]
      &v_t(s)-r_t(s,a)\\[-0.25ex]
      &\quad-\displaystyle\sum_{s'\in\states}P_t(s'\mid s,a)v_{t+1}(s')
    \end{aligned}
    & t=0,\ldots,T-1,\\[1.2ex]
    v_T(s)-r_T(s,a), & t=T.
  \end{cases}
\end{equation}
Dual feasibility is equivalent to $z_t^v(s,a)\geq 0$ for every
$(t,s,a)$.

\begin{proposition}[MDP duality and complementary slackness]
	\label{prop:mdp-dual-lp}
	Let $V^*=(V_t^*)_{t=0}^T$ be the optimal value function defined by the
	Bellman equations. Then $V^*$ is feasible and optimal for the dual
	program~\eqref{eq:mdp-dual-lp}. The primal and dual optimal values
	coincide and are equal to
	\[
		\sum_{s\in \states}\mu_0(s)V_0^*(s)
		=\sup_{\pi\in\strategies}J(\pi).
	\]
	Moreover, for every primal-feasible occupation measure $d$ and every
	dual-feasible vector $v$, one has the duality-gap identity
	\begin{equation}
		\label{eq:mdp-duality-gap}
		\sum_{s\in \states}\mu_0(s)v_0(s)
		-
		\sum_{t=0}^{T}\sum_{s\in \states}\sum_{a\in \actions}
		d_t(s,a)r_t(s,a)
		=
		\sum_{t=0}^{T}\sum_{s\in \states}\sum_{a\in \actions}
		d_t(s,a)z_t^v(s,a).
	\end{equation}
	Consequently, a primal-feasible $d$ and a dual-feasible $v$ are both
	optimal if and only if
	\begin{equation}
		\label{eq:mdp-complementary-slackness}
		d_t(s,a)z_t^v(s,a)=0,
		\qquad
		s\in \states,\quad a\in \actions,\quad t=0,\ldots,T.
	\end{equation}
	Equivalently, since both vectors are nonnegative, complementary
	slackness can be written as
	\[
		\sum_{t=0}^{T}\sum_{s\in \states}\sum_{a\in \actions}
		d_t(s,a)z_t^v(s,a)=0.
	\]
\end{proposition}

\begin{proof}
	Let $d$ be primal feasible. The initial and flow constraints imply the
	telescoping identity
	\begin{align*}
		\sum_{s\in \states}\mu_0(s)v_0(s)
		={} &
		\sum_{t=0}^{T-1}\sum_{s\in \states}\sum_{a\in \actions}
		d_t(s,a)
		\left(
		v_t(s)-\sum_{s'\in \states}P_t(s'\mid s,a)v_{t+1}(s')
		\right)      \\
		    & \quad+
		\sum_{s\in \states}\sum_{a\in \actions}d_T(s,a)v_T(s).
	\end{align*}
	Subtracting the primal objective from both sides gives
	\eqref{eq:mdp-duality-gap}. If $v$ is dual feasible, then every term on
	the right-hand side is nonnegative. This proves weak duality.

	The Bellman equations show that $V^*$ is dual feasible. Choose any
	policy that is greedy with respect to the optimal $Q$-function at every
	time--state pair, and let $d^*$ be its occupation measure. For every
	state--action pair to which $d^*$ assigns positive mass, the corresponding
	Bellman inequality is an equality. Hence the right-hand side of
	\eqref{eq:mdp-duality-gap} is zero for $(d^*,V^*)$. The two feasible
	solutions therefore attain the same objective value, proving strong
	duality and the optimality of $d^*$ and $V^*$.

	Finally, for arbitrary primal- and dual-feasible solutions, equality of
	their objective values is equivalent, by~\eqref{eq:mdp-duality-gap}, to
	the vanishing of a sum of nonnegative terms. This holds if and only if
	each product in~\eqref{eq:mdp-complementary-slackness} is zero.
\end{proof}

\paragraph{Interpretation of complementary slackness and matrix representation.}
For $v=V^*$, the slack has the particularly transparent form
\[
	z_t^{V^*}(s,a)=V_t^*(s)-Q_t^*(s,a),
	\qquad t=0,\ldots,T,
\]
where $Q_T^*(s,a)=r_T(s,a)$. Thus,
\[
	d_t^*(s,a)>0
	\quad\Longrightarrow\quad
	Q_t^*(s,a)=V_t^*(s)
	=\max_{a'\in \actions}Q_t^*(s,a').
\]
That is to say, an optimal occupation measure places positive mass only on
Bellman-optimal actions. Conversely, if a feasible occupation measure is
supported only on Bellman-optimal actions, then it is optimal.

This statement also clarifies the role of states with zero occupation
mass. If $\sum_{a\in \actions}d_t^*(s,a)=0$, complementary slackness imposes no
restriction on the action selected at $(t,s)$. Hence a policy that is
optimal for the objective $J$ under the fixed initial distribution
$\mu_0$ may be specified arbitrarily at states that are never reached.
A policy that is greedy with respect to $Q^*$ at every time--state pair
is stronger: it is Bellman optimal from every time and state, and in
particular is optimal under $\mu_0$.

The same relations can be written compactly in matrix form.
After
changing signs from reward maximization to cost minimization, let
$c$ be the vector obtained by stacking $-r_t(s,a)$, and let $A$ denote
the linear operator encoding the initial and flow constraints. The
primal program becomes
\[
	\underset{d}{\operatorname{minimize}}\ c^\top d
	\qquad
	\operatorname{subject\ to}\qquad
	Ad=b,\quad d\geq 0.
\]
Its primal--dual optimality conditions are
\[
	Ad=b,
	\qquad
	A^\top y+z=c,
	\qquad
	d\geq 0,
	\qquad
	z\geq 0,
	\qquad
	z^\top d=0.
\]
If, for this matrix representation, each dynamic flow row is written as incoming flow minus the current state mass, and the dynamic rows precede the initial-distribution rows, then with the block ordering
\[
	y=(v_1,\ldots,v_T,-v_0),
\]
the vector $z$ is exactly the stacked Bellman slack in
\eqref{eq:mdp-bellman-slack}.

As will be clear later, this matrix representation is convenient for a mean field game:  fixing a mean-field
flow produces an MDP of this form, but the operator $A$ and the cost
vector $c$ depend on that flow. Adding the consistency condition that
the representative player's occupation measure equals the mean field
then leads directly to the primal--dual MFG formulation.

\section{Finite-population games}

We now proceed from MDPs, where there is a single agent making decisions, to finite-population games. Building on the framework introduced earlier for MDPs, we present a basic setup of stochastic, dynamic games with finitely many players.

\subsection{Finite-population game setup}
\label{sec:finite-pop-game-setup}

We consider $N$ players, indexed by $i=1,\dots,N$. We denote $[N] = \{1,\dots,N\}$ and we will use bold symbols for vectors involving multiple players. Motivated by mean field games, we focus on homogeneous games (in which the players have similar dynamics and reward functions) in which the interactions are symmetric (i.e., through the empirical distribution of states and actions).

The \boldindex{finite-population game} is defined by the tuple
\[
  (N,\states,\actions,T,\mu_0,P,r,\Pi).
\]
Its components are as follows:
\begin{itemize}
	\item $N$ is the number of players;
	\item $\states$, $\actions$, and $T$ are as in the MDP setting discussed above;
	\item $\mu_0 \in \Delta(\states)$ is the initial distribution of each player's state;
	\item $P: \timesteps \times \states \times \actions \times \Delta(\states\times\actions) \to \Delta(\states)$ is the transition kernel for the state dynamics of each player;
	      $P_t(\cdot|s, a, \nu)$ denotes the probability distribution of the next state, when, at time $t$, the player's state and action are respectively $s$ and $a$, and the state-action distribution is $\nu$;
	\item $r: \timesteps \times \states \times  \actions \times \Delta(\states\times\actions) \to \RR^N,$ is the reward function;
	      $r^i_t(\state, \action, \nu)$ is the reward obtained at time $t$ by player $i$ in state $s$ using action $a$ if the state-action population distribution is $\nu$;
	\item $\strategies$ is the set of policies for an individual player, which are the sequences $\pi = (\pi_t)_{t \in \timesteps}$ where for each $t \in \timesteps$, $\pi_t: \states^N \to \Delta(\actions)$; we denote by $\bstrategies=\strategies^N$\notationindex{boldPi}{$\bstrategies$: policy-profile space} the set of policy profiles\index{policy profile}.
\end{itemize}
In the homogeneous examples below, all players share the same scalar reward function; equivalently, $r^i_t=r_t$ for every $i$, and we simply write $r_t(s,a,\nu)$.

Notice that we consider only \boldindex{Markovian policies}. This is motivated by the fact that MDPs (as described in the previous section) admit optimal Markovian policies so, given the other players' policies, each player can find an optimal policy that is Markovian.

If the players use the policy profile $\bspi = (\pi^1,\dots,\pi^N) \in \bstrategies$,
the total expected reward of player $i$ is:
\begin{equation}
	J^i(\bspi) = \EE\left[ \sum_{t=0}^T r^i_t(s^i_t, a^i_t, \nu^N_t)\right],
\end{equation}
where for each $t \in \timesteps$, $\nu^N_t(s,a)$\notationindex{nuN}{$\nu^N_t$: empirical population distribution} $= \frac{1}{N} \sum_{j=1}^N \indic_{(s^j_t,a^j_t)}(s,a)$ is the \boldindex{empirical population distribution}, subject to:
\[
	\begin{cases}
		s_0^i \sim \mu_0,
		\\
		a^i_t \sim \pi^i_t(\cdot|\bfs_t), \quad t \in \timesteps,
		\\
		s^i_{t+1} \sim P_t(\cdot|s^i_t,a^i_t,\nu^N_t), \quad t \in \timesteps\backslash\{T\},
	\end{cases}
\]
where the initial states are independent with common law $\mu_0$; at each time, the players' actions are conditionally independent given the current state profile, and the next states are conditionally independent given the current state and action profiles.

Given a policy profile $\bspi \in \bstrategies$ and a policy $\tilde\pi^i \in \strategies$ for player $i$, we will denote by $J^i(\tilde\pi^i, \bspi^{-i}) = J^i(\pi^1,\dots,\pi^{i-1},\tilde\pi^i,\pi^{i+1}, \dots,\pi^{N})$ the total reward obtained by player $i$ when she uses policy $\tilde\pi^i$ and the other players use $\bspi^{-i} = (\pi^1,\dots,\pi^{i-1},\pi^{i+1}, \dots,\pi^{N})$.

\subsection{Nash equilibrium}

There are multiple optimality criteria for finite-population games. We will mainly focus on a non-cooperative setting and adopt the  notion of Nash equilibrium, which  corresponds to a situation in which each player's policy is optimal for her MDP, given the behavior of the other players.
Note that when the policies of the other players are fixed, each player faces a standard MDP.
\begin{definition}[Nash equilibrium in finite-agent game]
	For $\epsilon\ge 0$, a policy profile $\hat\bspi = (\hat\pi^1,\dots,\hat\pi^N) \in \bstrategies$ is an \boldindex{$\epsilon$-Nash equilibrium} if for all $i \in [N]$,
	\[
		J^i(\hat\bspi) \ge \max_{\pi^i \in \strategies} J^i(\pi^i,\hat\bspi^{-i})-\epsilon.
	\]
	A $0$-Nash equilibrium is simply called a \boldindex{Nash equilibrium} (\boldindex{NE}).
\end{definition}
In other words, a policy profile is a Nash equilibrium if no player can improve her payoff only by changing her own policy (i.e., by a \boldindex{unilateral deviation}).

Throughout the rest of this monograph, unless explicitly stated otherwise, we work with finite state and action spaces and a finite time horizon. In this setting, the stochastic games considered in this chapter have at least one Nash equilibrium; later chapters use the same finite-dimensional convention when formulating discrete MFGs and algorithms.
\begin{theorem}[Existence of NE]
	\label{thm:existence-NE}
	There exists at least one Nash equilibrium for the game described above.
\end{theorem}
The proof is provided in Appendix~\ref{app:proof-existence-NE}.
When the spaces are continuous, assumptions on the transitions and the reward functions are needed in order to ensure existence of a Nash equilibrium.

In line with Remark~\ref{rem:setting-statio}, one can also consider an infinite horizon, discounted and stationary setting in the context of $N$-player games.

\paragraph{Example 1.}
	Finite-agent games can have several Nash equilibria, even in simple one-stage examples. As a first illustration, consider a game with two players, a single state $\states=\{s_0\}$, horizon $T=0$, and action space $\actions=\{\text{Apple},\text{Banana}\}$. If both players choose Apple, or if both players choose Banana, each receives payoff $1$; otherwise each receives payoff $0$. This game has two pure Nash equilibria, corresponding to the two coordinated action profiles.

	A slightly richer variant keeps the same one-stage structure but allows several equilibria. Consider a game with two players, Alex and Bob. Fix an integer $M\ge 1$. There is a single state, $\states=\{s_0\}$, a horizon $T=0$, and action space $\actions=\{0,1,\dots,M\}$. If a player chooses action $a$, she requests $a$ apples and $M-a$ bananas. There are $M$ apples and $M$ bananas available in total. The initial distribution is $\mu_0=\delta_{s_0}$; since $T=0$, the transition kernel is irrelevant.

	Alex prefers apples to bananas, while Bob has the opposite preference. For an action profile $(a^1,a^2)$, let $\nu^2$ denote the empirical distribution of this two-player game. It satisfies
	\[
		2\sum_{\bar a=0}^M \bar a\,\nu^2(s_0,\bar a)=a^1+a^2.
	\]
	If the total request is feasible, namely if $a^1+a^2=M$, the rewards are
	\[
		r^1_0(s_0,a^1,\nu^2)=2a^1-M,
		\qquad
		r^2_0(s_0,a^2,\nu^2)=M-2a^2.
	\]
	If $a^1+a^2\neq M$, then at least one type of fruit is over-requested, and both rewards are set equal to the penalty $-(M+1)$.
	The pure Nash equilibria are precisely the action profiles $(a^1,a^2)=(m,M-m)$, for $m=0,\dots,M$. At such a profile, the requests exactly exhaust the available fruits. If one player deviates unilaterally, the feasibility condition fails and the deviating player receives the penalty $-(M+1)$, which is smaller than any feasible payoff. Conversely, if $a^1+a^2\neq M$, either player can switch to the unique action that restores feasibility and obtain a payoff at least $-M$, so the profile is not a Nash equilibrium.

	This example illustrates multiplicity of Nash equilibria; the MFG models considered later use more symmetric interaction structures.

\paragraph{Example 2.}
A well-known game that demonstrates the difference between NE and other types of equilibrium criteria is the Prisoner's dilemma. The game involves two players, say
Ted and Jack, who are arrested because of a crime they committed together. They will be interrogated separately, and each chooses whether to confess or remain silent.  If Ted confesses and Jack remains silent, then Jack will be sentenced to 5 years in prison while Ted can be freed, and vice versa. Meanwhile, if both of them confess, then both get a sentence of 3 years in prison; if both of them keep silent, then both get a sentence of 1 year in prison. The payoff matrix is given in Table~\ref{table:prisoner-dilemma}: each cell contains a pair $(a,b)$, where $a$ represents the payoff of Ted and $b$ represents the payoff of Jack.
\begin{table}
	\centering
	\begin{tabular}{|c|c|c|}
		\hline
		              & Jack confesses & Jack silent \\
		\hline
		Ted confesses & (-3,-3)        & (0, -5)     \\
		\hline
		Ted silent    & (-5,0)         & (-1,-1)     \\
		\hline
	\end{tabular}
	\caption{Payoff matrix for the prisoner's dilemma \label{table:prisoner-dilemma}}
\end{table}
This game is static. It can be recast in the dynamic game framework discussed above with just one state and one time step: $N=2$, $\states = \{0\}$, $\actions = \{0,1\}$ (where $a=0$ represents confess and $a=1$ represents staying silent), $T=0$, $\mu_0 = \delta_0$, $P(s'|s,a,\nu) = \mathbf{1}_{\{s'=0\}}$, and the reward, where $m=\nu(0,0)$ denotes the fraction of players who confess:
\[
	r(s,a,\nu) =
	\begin{cases}
		-3, & \hbox{ if $a=0$, $m = 1$,}
		\\
		0,  & \hbox{ if $a=0$, $m = 1/2$,}
		\\
		-5, & \hbox{ if $a=1$, $m = 1/2$,}
		\\
		-1, & \hbox{ if $a=1$, $m = 0$.}
	\end{cases}
\]
In this example, confessing is strictly better for each player regardless of the other player's choice: it gives $0$ rather than $-1$ against silence, and $-3$ rather than $-5$ against confession. Hence both players confessing is the unique Nash equilibrium.
However, this strategy is suboptimal in comparison to the case where they collaborate and both remain silent: they get $-3$ if they both confess whereas they get $-1$ if they both stay silent. This latter strategy is called a social optimum strategy, as collaboration between them optimizes the overall cost of two years in prison instead of a total of six years in prison for the Nash equilibrium. As a result, in a Nash equilibrium, both players find themselves in a worse state than if they had cooperated with each other in the decision-making process.

\paragraph{Example 3.}
Coordination games are typical examples exhibiting multiple Nash equilibria. For example, consider the same two actions as in the prisoner's dilemma, but with the payoff matrix in Table~\ref{table:coordination-game}: coordination is rewarded relative to mis-coordination, and the two coordinated outcomes give the same payoff.
\begin{table}
	\centering
	\begin{tabular}{|c|c|c|}
		\hline
		              & Jack confesses & Jack silent \\
		\hline
		Ted confesses & (-3,-3)        & (-5,-5)     \\
		\hline
		Ted silent    & (-5,-5)        & (-3,-3)     \\
		\hline
	\end{tabular}
	\caption{Payoff matrix for the coordination game \label{table:coordination-game}}
\end{table}
In this case, there are two pure Nash equilibria: both confess, and both stay silent. In each case, each player can only be worse off by deviating to the other pure action.

\section{Bibliographic notes}

For more details on the theory of MDPs and in particular the existence of optimal policies, we refer to the monographs~\citep{puterman2014markov,bertsekas1996stochastic}. See, for instance, \citep[Proposition 4.4.3]{puterman2014markov} for the existence of pure policies in a finite-horizon setting as discussed in Section~\ref{sec:basic-MDP}.

The DPP was  first established for the value function of MDPs in \citep{bellman1957markovian,bellman1966dynamic}, and later
derived for the Q-function in a learning framework of MDP in \citep{watkins1989learning,watkins1992q}. We refer for instance to the books of~\citep{sutton2018reinforcement} and~\citep{bertsekas1996stochastic} for proofs that the value functions $V^*$ and $Q^*$ are indeed the (unique) solutions of the Bellman equations discussed in Section~\ref{sec:mdp-value-bellman}.
We  also refer to~\citep{bertsekas1996neuro} and~\citep{sutton2018reinforcement} for more background and details in the context of learning, and to~\cite{bertsekas1996stochastic} for the mathematical theory. There are further extensions such as partially observable MDPs; see~\cite{aastrom1965optimal}.

The original form of linear program for Markov decision process first appears in \citep{manne1960linear} for discrete-time MDPs with finite state-action spaces. It is then extended  to semi-MDPs (\cite{osaki1968linear}) and later to approximate dynamic programming (\cite{schweitzer1985generalized,de2003linear}), as well as to continuous-time stochastic control  (\cite{stockbridge1990time,bhatt1996occupation,kurtz1998existence}) and singular control (\cite{taksar1997infinite,kurtz2001stationary}) with continuous state-action spaces, among
others. See also \cite{wolfe1962linear,derman1962sequential,denardo1970linear,hordijk1979linear}.
The proof of Proposition~\ref{prop:MDP_LP} and the linear program formulation for MDP in the discounted and average reward settings can be found for example in~\citep{de2003linear,huang1997lp,wang2016online}.

The classic book~\citep{vonneumann1947} laid the foundations of modern game theory. The concept of Nash equilibrium was introduced by~\citep{nash1950equilibrium} and has found a wide range of applications, especially in economics. We refer, for example, to the monograph~\citep{fudenberg1991game}. Other notions of equilibrium, such as correlated equilibria, are also relevant to algorithmic methods; see~\citep{roughgarden2010algorithmic} for more details. Stochastic games were introduced by~\citep{shapley1953stochastic}. For the existence of Nash equilibria in continuous spaces, we refer the interested reader to~\citep[Section 12.2]{fudenberg1991game}. For the existence of Nash equilibria in the infinite-horizon discounted stationary setting, we refer, for example, to~\citep{fink1964equilibrium,takahashi1964equilibrium,sobel1971noncooperative}. In the machine-learning literature, Markov games are a standard framework for multi-agent reinforcement learning; see, for example,~\citep{littman1994markov,hu2003nash}.

For zero-sum games with continuous action spaces, convex-concave assumptions are often used to obtain saddle-point formulations and minimax identities. From an optimization perspective, these assumptions also ensure the absence of a duality gap, as in Sion's generalized minimax theorem~\citep{sion1958general} and the classical minimax theorem of~\citep{von1928theorie}.

\chapter{Discrete Mean Field Games}
\label{chap:discrete-mfg}

\section{Motivating examples}
\label{sec:motivating-ex}

In many applications, multi-player systems involve a  large number of agents who interact and make decisions while trying to anticipate other players' behavior. For each player, it is rather complicated to keep track of the interactions between all players and, in many cases, it is not necessary to do that in order to find an approximately optimal behavior: the player only needs to predict the evolution of some aggregate quantities.

For instance, in macro-economic models such as the Aiyagari model with infinitely many heterogeneous agents, the players are households who interact through aggregate quantities such as the wage and the interest rate. These quantities are determined through equilibrium conditions which involve the whole population. When the population is  large,  players behave as price takers, in the sense that they decide their consumption level believing that the aggregate quantities are fixed and not affected by their individual behavior. Strictly speaking, when there is a finite number of households in the market, prices \emph{are} affected by the decisions of each player. Nevertheless, the impact of a single household may be negligible when the population is  large.

Another example is crowd motion. When a pedestrian tries to find her trajectory in order to reach a target location, she observes the crowd and tries to anticipate how the crowd will evolve, in order to avoid congested areas. To this end, the pedestrian usually takes the crowd's behavior as given and assumes that her own behavior will not change that of the crowd. This assumption is an approximation; nevertheless, it simplifies the computation of the optimal trajectory and is usually sufficient for practical purposes. Indeed, when the crowd is  huge, its evolution is barely affected if a single player decides to change her trajectory.

In the mean field game literature, one of the first examples is about finding the time when a meeting should start. This is introduced as a static game, in which the players are individuals who would like to attend a meeting and decide their arrival time, so that their actual arrival time is their chosen time up to a small perturbation. The scheduled meeting time is announced to everyone, however,  the meeting only starts when  a certain proportion of individuals has arrived. The players want to balance costs induced by the social pressure of being late and the inconvenience of arriving  earlier than necessary.

In this example, the notion of mean-field interaction is introduced in that each individual player interacts only through the quantile (i.e., actual starting time of the meeting)  or the distribution of the arrival time of the entire  population.
In the case of an infinite population, each single agent has no influence on averaged quantities, as
the population distribution does not change if a single agent changes her state and/or action. Hence the interaction is  {\it weak}. This weak-interaction viewpoint is one of the main ideas behind MFGs and simplifies the analysis of finite-agent games.

More examples can be found in finance where traders interact through prices, carbon emissions where electricity producers interact through the pollution level, and epidemic management, where citizens interact through infection probabilities. We discuss several applications in Chapter~\ref{chap:applications}.

\section{Intuition from $N$-player games}
\label{sec:sym-games-finite-infi}

Before introducing the  mean field game setup, let us consider first an $N$-player game. Notice that we can rewrite the total reward of player $i\in [N]$ as follows.
If  players use the policy profile $\bspi = (\pi^1,\dots,\pi^N) \in \bstrategies$,
the total expected reward of player $i$ is:
\begin{equation}
	J^i(\pi^i,\bspi^{-i})
	= J^i(\bspi)
	= \EE\left[ \sum_{t=0}^T r_t(s^i_t, a^i_t, \nu^N_t)\right]
\end{equation}
subject to:
\[
	\begin{cases}
		s^j_0 \sim \mu_0, \quad j \in [N]
		\\
		a^j_t \sim \pi^j_t(\cdot|\bfs_t), \quad j \in [N], t=0,\dots,T \\
		s^j_{t+1} \sim P_t(\cdot|s^j_t,a^j_t,\nu^N_t), \quad j \in [N], t=0,\dots,T-1
		\\
		\nu^N_t(s,a) = \frac{1}{N} \sum_{j=1}^N \indic_{(s^j_t,a^j_t)}(s,a), \quad (s,a)\in\states\times\actions,\ t=0,\dots,T.
	\end{cases}
\]

Since player $i$ interacts at time $t$ with the other players through empirical distributions, it is natural to look for symmetric Markovian policies. When the aggregate flow has a deterministic mean-field limit, that limit can be treated as a prescribed function of time. This motivates the decentralized policy $\pi^j_t(\cdot|\bfs_t)=\pi_t(\cdot|s^j_t)$ used below.

If the players' individual evolution is stochastic, then $\nu^N$ evolves stochastically. Under homogeneity, as $N$ tends to infinity, \boldindex{propagation of chaos} results often imply that $\nu^N$ converges to a deterministic limit. This motivates the search for \boldindex{decentralized policies} (i.e., taking as input only the individual state and not the mean field). Conversely, under suitable assumptions, policies of the form $\pi^j_t(\cdot|\bfs_t) = \pi_t(\cdot|s^j_t)$ can yield an $\epsilon$-Nash equilibrium with $\epsilon \to 0$ as $N \to +\infty$. MFG theory provides a framework to formalize this idea in a mathematically precise way.

\section{Mean field setup}
\label{sec:mf-setup}

We now consider an infinite population of homogeneous and rational players, corresponding at least intuitively to the limit of the game described in the previous section. We model the interactions between a representative player and the population, represented by a probability distribution. Many variations have been considered in the literature and will be mentioned in the bibliography notes at the end of this chapter.
Formally, a \boldindex{mean field game} (\boldindex{MFG}) is defined by a tuple $(\states,\actions,T,\mu_0,P,r,\Pi)$, where
\begin{itemize}
	\item $\states,\actions,T,\mu_0,P,$ and $r$ are as in the finite-population game of Section~\ref{sec:finite-pop-game-setup} but the reward function takes scalar values;
	\item $\strategies$ is the set of policies, which are the sequences $ \pi = (\pi_t)_{t \in \timesteps}$ where for each $t \in \timesteps$, $\pi_t: \states \to \Delta(\actions)$.
\end{itemize}
Here, we focus on policies depending only on the individual state. We refer to such policies as \boldindex{decentralized policies}.
The set $\Delta(\states \times \actions)$ of probability distributions on $\states \times \actions$ is interpreted as the set of values for the state-action distribution, also called \boldindex{mean field}; a generic element is denoted by $\nu$\notationindex{nu}{$\nu$: mean field}.
Of particular interest are examples where the interaction depends on the mean field only through its state marginal. In those cases, for $\nu_t\in\Delta(\states\times\actions)$, we write $\mu_t(s)=\sum_{a\in\actions}\nu_t(s,a)$ and may use $\mu_t$ in formulas for $P_t$ or $r_t$ as shorthand for this dependence on $\nu_t$.

Consider a \boldindex{representative player} facing a given mean field $\nu$. This gives rise to an MDP in which the player tries to find her best policy while the mean field is fixed. More precisely, if the mean field sequence is $\nu = (\nu_t)_{t \in \timesteps}$ and the  player adopts the policy sequence $\pi = (\pi_t)_{t \in \timesteps} \in \Pi$,
the total expected reward for this player is:
$$
	J(\pi; \nu) = \EE\left[ \sum_{t \in \timesteps} r_t(s_t, a_t, \nu_t)\right]
$$
subject to:
\[
	\begin{cases}
		s_0 \sim \mu_0,
		\\
		a_t \sim \pi_t(\cdot|s_t), \quad t \in \timesteps,
		\\
		s_{t+1} \sim P_t(\cdot|s_t,a_t,\nu_t), \quad t \in \timesteps\backslash\{T\}.
	\end{cases}
\]
For a given policy sequence $\pi$, we will denote by $\mu^\pi$\notationindex{mupi}{$\mu^\pi$: state mean field induced by $\pi$} the \boldindex{state mean field} induced by this policy, i.e., the state distribution sequence obtained when every player uses $\pi$. Its evolution is given by:
\begin{equation}
	\label{eq:def-mu-pi}
	\begin{cases}
		\mu_0^\pi(s) = \mu_0(s), \quad s \in \states
		\\
		\displaystyle \mu_{t+1}^\pi(s) = \sum_{s' \in \states} \sum_{a \in\actions} \mu_t^\pi(s') \pi_t(a|s') P_t(s|s',a,\nu_t^\pi), \quad s \in \states, t \in \timesteps\backslash\{T\},
	\end{cases}
\end{equation}
where we denote by $\nu^\pi$\notationindex{nupi}{$\nu^\pi$: state-action mean field induced by $\pi$} the \boldindex{state-action mean field}, defined, for every $(s,a) \in \states \times \actions$, as:
\begin{equation}
	\label{eq:def-nu-pi}
	\nu^\pi_t(s,a) = \mu^\pi_t(s) \pi_t(a|s).
\end{equation}

Here again, in line with Remark~\ref{rem:setting-statio}, one can also consider an infinite horizon, discounted and stationary setting in the context of MFGs.

\section{Solution concept: mean field Nash equilibrium}

We now turn to extensions of the notions of Nash equilibrium to the mean field setup.

Recall that in a non-cooperative game, when all the players are indistinguishable as in the mean field setting, a Nash equilibrium corresponds to a situation in which no player can be better off by deviating from the behavior of the rest of the population. This is formalized as follows.

\begin{definition}[Nash equilibrium in mean field setup]
	\label{def:mfNE}
	For $\epsilon\ge 0$, a policy $\hat\pi \in \strategies$ is an \boldindex{$\epsilon$-mean field Nash equilibrium} if:
	\[
		J(\hat\pi; \nu^{\hat\pi}) \ge \max_{\pi \in \strategies} J(\pi; \nu^{\hat\pi})-\epsilon,
	\]
	where $\nu^{\hat\pi}$ is defined by~\eqref{eq:def-nu-pi} with policy $\hat\pi$.
	Equivalently, the pair $(\hat\pi,\hat\nu)$ is an $\epsilon$-mean field Nash equilibrium if the following two conditions are satisfied:
	\begin{enumerate}
		\item Approximate optimality: $\hat\pi$ is $\epsilon$-optimal against mean field $\hat\nu$.
		\item Consistency: $\hat\nu = \nu^{\hat\pi}$, which is the mean field generated when everyone uses policy $\hat\pi$.
	\end{enumerate}

	A $0$-mean field Nash equilibrium is simply called a \boldindex{mean field Nash equilibrium} (\boldindex{MFNE}).
\end{definition}
This corresponds to a fixed point problem for the map which is the composition of the best response operator and the mean field operator. We will elaborate more in Chapter~\ref{chap:analytical-MFG}.

\section{Examples}

We now describe three examples of mean field models. For simplicity and to be consistent with the framework described earlier, we consider discrete-time models with finite state and action spaces. Extensions to continuous spaces and/or continuous time will be discussed in the bibliography notes.

\subsection{Crowd motion}
\label{sec:ex-mf-crowd-motion}
In this example describing the evolution of agents in a grid world,  the state of a player is her position. During one time step, each player can move to one of the neighboring positions. The reward of a typical player consists of several terms: one term penalizing movement suggesting that the player has a preference for not moving, one term encouraging staying in some state(s) indicating that the player has some preferred location, and one term penalizing a high population density at her state assuming that  the player tends to avoid crowded states. Since the objective accumulates these terms over several time steps, the Nash equilibrium depends on the balance between movement costs, target preferences, congestion, and the dynamics.

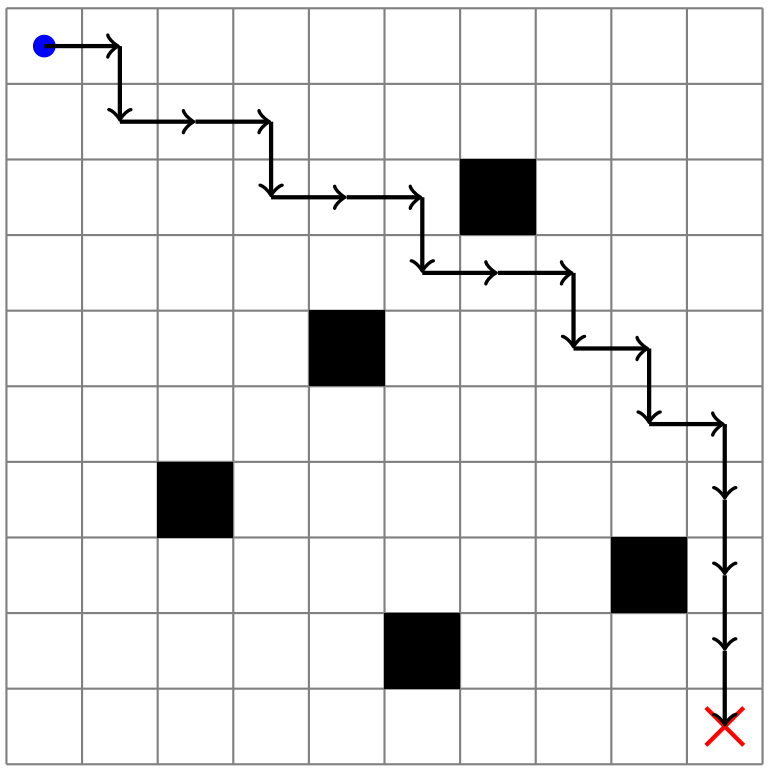
\begin{figure}
	\begin{center}
		\begin{tikzpicture}
			\def\cellsize{0.5}
			\def\gridsize{10}

			\foreach \x in {0,1,...,10} {
					\draw[gray] (\x*\cellsize,0) -- (\x*\cellsize,\gridsize*\cellsize);
					\draw[gray] (0,\x*\cellsize) -- (\gridsize*\cellsize,\x*\cellsize);
				}

			\foreach \x/\y in {2/3, 4/5, 6/7, 8/2, 5/1} {
					\fill[black] (\x*\cellsize,\y*\cellsize) rectangle (\x*\cellsize+\cellsize, \y*\cellsize+\cellsize);
				}

			\draw[red, thick] (9.25*\cellsize, 0.25*\cellsize) -- (9.75*\cellsize, 0.75*\cellsize);
			\draw[red, thick] (9.75*\cellsize, 0.25*\cellsize) -- (9.25*\cellsize, 0.75*\cellsize);

			\fill[blue] (0.5*\cellsize, 9.5*\cellsize) circle (0.15*\cellsize);

			\draw[thick,->] (0.5*\cellsize, 9.5*\cellsize) -- (1.5*\cellsize, 9.5*\cellsize);
			\draw[thick,->] (1.5*\cellsize, 9.5*\cellsize) -- (1.5*\cellsize, 8.5*\cellsize);
			\draw[thick,->] (1.5*\cellsize, 8.5*\cellsize) -- (2.5*\cellsize, 8.5*\cellsize);
			\draw[thick,->] (2.5*\cellsize, 8.5*\cellsize) -- (3.5*\cellsize, 8.5*\cellsize);
			\draw[thick,->] (3.5*\cellsize, 8.5*\cellsize) -- (3.5*\cellsize, 7.5*\cellsize);
			\draw[thick,->] (3.5*\cellsize, 7.5*\cellsize) -- (4.5*\cellsize, 7.5*\cellsize);
			\draw[thick,->] (4.5*\cellsize, 7.5*\cellsize) -- (5.5*\cellsize, 7.5*\cellsize);
			\draw[thick,->] (5.5*\cellsize, 7.5*\cellsize) -- (5.5*\cellsize, 6.5*\cellsize);
			\draw[thick,->] (5.5*\cellsize, 6.5*\cellsize) -- (6.5*\cellsize, 6.5*\cellsize);
			\draw[thick,->] (6.5*\cellsize, 6.5*\cellsize) -- (7.5*\cellsize, 6.5*\cellsize);
			\draw[thick,->] (7.5*\cellsize, 6.5*\cellsize) -- (7.5*\cellsize, 5.5*\cellsize);
			\draw[thick,->] (7.5*\cellsize, 5.5*\cellsize) -- (8.5*\cellsize, 5.5*\cellsize);
			\draw[thick,->] (8.5*\cellsize, 5.5*\cellsize) -- (8.5*\cellsize, 4.5*\cellsize);
			\draw[thick,->] (8.5*\cellsize, 4.5*\cellsize) -- (9.5*\cellsize, 4.5*\cellsize);
			\draw[thick,->] (9.5*\cellsize, 4.5*\cellsize) -- (9.5*\cellsize, 3.5*\cellsize);
			\draw[thick,->] (9.5*\cellsize, 3.5*\cellsize) -- (9.5*\cellsize, 2.5*\cellsize);
			\draw[thick,->] (9.5*\cellsize, 2.5*\cellsize) -- (9.5*\cellsize, 1.5*\cellsize);
			\draw[thick,->] (9.5*\cellsize, 1.5*\cellsize) -- (9.5*\cellsize, 0.5*\cellsize);

		\end{tikzpicture}
	\end{center}
	\caption{Grid world example: $10 \times 10$ grid with obstacles represented by black cells and a target state represented by a red cross. The arrows represent a possible path from an initial state (blue dot) to the target state. \label{fig:grid-world-diagram}}
\end{figure}

To be specific, using the notations introduced above, the game model can be described as follows:
\begin{itemize}
	\item State space $\states = \{1,\dots,d_1\} \times \{1,\dots,d_2\}$, and action space $\actions = \{-1,0,1\} \times \{-1,0,1\}$. Some states can be tagged as ``forbidden'' states (e.g., walls) and some others as ``target'' states (e.g., locations of interest). Let us denote by $\states_{\mathrm{target}}$ the set of target states.
	\item The transition kernel $P$ is defined through a transition function and an additive noise.
	      \begin{itemize}
		      \item The transition function is defined as: $F: \timesteps \times \states \times \Delta(\states\times\actions) \times \actions  \to \states,$ where $F_t(\state, \nu, \action)$ denotes the next state that the player can expect to reach in the absence of noise. In this example, the transition does not depend on the mean field; the argument $\nu$ is kept only to match the general framework. If there are no obstacles, the next state is $\state + \action = (\state_1+\action_1, \state_2+\action_2)$, whereas if this new state is forbidden, then the player stays in $\state$, so $F_t(\state, \nu, \action) = \state$.
		      \item The noise term at time $t$ is a random variable distributed over the action space; it is added to the state computed by $F_t$, provided the resulting state is not forbidden.
	      \end{itemize}
	\item The individual reward is of the form $r(\state,\action,\mu) = -\|\action\| + \indic_{\state \in \states_{\mathrm{target}}} - \mu(\state)$. Each can be multiplied by a coefficient to weigh its importance compared with other ones. In particular, the first term can be multiplied by a coefficient which increases with $\mu(\state)$ in order to model \boldindex{congestion} effects (i.e., the fact that moving through a crowded state is more expensive than moving through an empty state).
\end{itemize}
Including forbidden states (obstacles) allows to construct examples with complex geometries. Furthermore, other dynamics and reward functions can be considered. An example of a grid world state space is depicted in Figure~\ref{fig:grid-world-diagram}: it is a $10 \times 10$ example, with 5 forbidden states and one target state. An example of a path from an initial state to the target set is represented by arrows.

We will present a numerical illustration in Section~\ref{sec:crowd-motion-numerics}.

\subsection{Epidemics}
\label{sec:ex-mf-epidemics}

Another natural application of games with a large number of players is epidemic management. Here, following the well-known \boldindex{susceptible-infected-susceptible model} (\boldindex{SIS model}), we focus on a simple game with two states: susceptible ($S$) and infected ($I$). Each agent can choose her level of socialization: $a=D$ (perfect social distancing) or $a=H$ (high level of socialization). A susceptible agent can become infected, with a probability that increases with the proportion of infected individuals in the population and with her level of socialization. After being infected, an agent recovers and becomes susceptible again with a fixed probability. Such transitions are a discrete-time approximation of the continuous-time SIS model, taking into account the fact that agents can influence their probability of getting infected. The reward function is defined to reflect the idea that each individual wants to socialize but does not want to be infected.

\begin{figure}
	\begin{center}
		\begin{tikzpicture}
			\node[draw, circle, minimum size=1cm] (S) at (0,0) {\large $S$};
			\node[draw, circle, minimum size=1cm] (I) at (4,0) {\large $I$};

			\draw[->, thick, bend left=30] (S) to node[above] {\large $w_{I} \mu(I) \indic_{a=H}$} (I);

			\draw[->, thick, bend left=30] (I) to node[below] {\large $p_{I \to S}$} (S);
		\end{tikzpicture}
	\end{center}
	\caption{Epidemics example. Each agent has two possible states: susceptible (S) and infected (I). The arrows represent the transition probabilities from $S$ to $I$ (infection) and from $I$ to $S$ (recovery). \label{fig:epidemics-diagram}}
\end{figure}
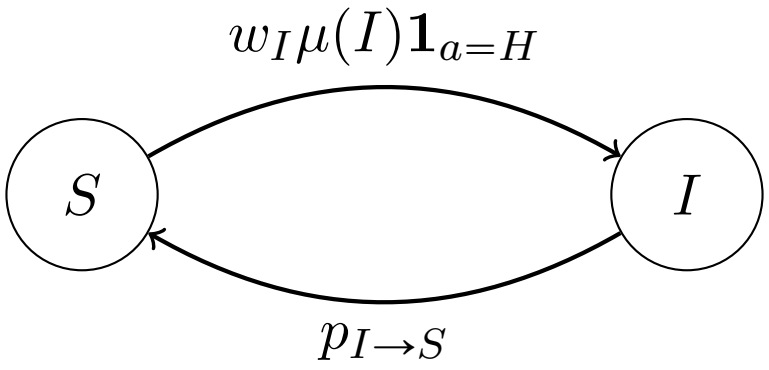

The model can be described as follows:
\begin{itemize}
	\item State space $\states = \{S, I\}$, and action space $\actions = \{D,H\}$;
	\item The transition kernel $P$ is such that:
	      \begin{align*}
		       & P(I|S, a, \nu) = w_{I} \mu(I)\indic_{a=H},
		      \qquad P(S|S,a,\nu) = 1 - P(I|S,a,\nu),
		      \\
		       & P(S|I,a,\nu) = p_{I \to S},
		      \qquad P(I|I,a,\nu) = 1 - p_{I \to S}.
	      \end{align*}
	\item The individual reward is of the form $r(\state,\action,\nu) = -w_{D}\indic_{\action=D} -\indic_{\state = I}$.
\end{itemize}
Here $w_D\ge 0$, $w_I\in[0,1]$, and $p_{I \to S} \in [0,1]$.
The transitions depend on the mean field $\nu$ only through the state distribution, $\mu(s) = \sum_a \nu(s,a)$. Figure~\ref{fig:epidemics-diagram} displays a representation of this model.
We will present the numerical analysis in Section~\ref{sec:appli-sis} for this model.

	\subsection{Traffic routing}
	\label{sec:ex-mf-traffic-routing}

	A third important class of discrete MFG models comes from traffic routing. The road network is represented by a directed graph, and each agent is a driver trying to reach a destination while avoiding congested links. In a discrete-time model, the state of an agent at time $t$ can be described by $s_t=(\ell_t,w_t)$, where $\ell_t$ is the current link and $w_t$ is the remaining number of time steps before the driver can move again. The action is the next link to take, restricted to the successors of $\ell_t$.

	The mean field enters through the congestion on each link. If $\mu_t(\ell)$ denotes the mass of drivers on link $\ell$ at time $t$, the waiting time after entering link $\ell$ can be modeled by a function $c(\ell,\mu_t(\ell))$. The reward typically penalizes time spent before reaching the destination, for instance through a term of the form $-\mathbf{1}_{\ell\neq D}$. This gives a finite-state, discrete-time model in which congestion affects the transition dynamics. A numerical illustration based on a Braess-type network is presented in Section~\ref{sec:appli-traffic}.

\section{Links between mean field and finite-agent models}

In Section~\ref{sec:sym-games-finite-infi}, we have justified the mean field limit in an intuitive way. In fact, the connections between the finite-agent models and their mean field counterparts can be justified in a rigorous way under suitable assumptions. Indeed, there are two ways to connect these models. First, one can start from a sequence of finite-agent models with a growing number of agents and show that the sequence of equilibria  converges to the equilibrium  of the mean-field model. Secondly, one can show that the equilibrium  in the mean field model provides an approximate equilibrium  in a corresponding finite-agent model, and the approximation quality improves with the number of agents.

The latter type of connection has several advantages: first, it does not require the existence of an (exact) equilibrium  in the finite-agent model; second, it is of practical importance since it justifies the interest for solving mean field problems. In this section, we will thus focus on this direction.

We start with the Nash equilibrium concept. We will discuss in Section~\ref{sec:existence-mfne} existence and uniqueness results for MFGs but for now we provide the following result which shows that, under a simplifying assumption on the dynamics, any MFNE policy, assuming there exists one, yields an approximate NE for the $N$-player game against decentralized deviations when all the agents use this policy.

\begin{theorem}[Approximate Nash equilibrium for decentralized deviations]
	\label{thm:approx-Nash-mf}
	Assume that the transition kernel does not depend on the mean field, \ie, $P_t(\cdot|s,a,\nu)=P_t(\cdot|s,a)$ for every $t\in\timesteps\backslash\{T\}$, $s\in\states$, $a\in\actions$, and $\nu\in\Delta(\states\times\actions)$. Assume also that, for every $t\in\timesteps$, $s\in\states$, and $a\in\actions$, the map $\nu\mapsto r_t(s,a,\nu)$ is continuous on $\Delta(\states\times\actions)$. Let $\hat\pi$ be a MFNE policy. Then, for any $\epsilon>0$, there exists an integer $N'$ such that for all $N \ge N'$, the policy profile $(\hat\pi,\hat\pi,\dots,\hat\pi)$ is an $\epsilon$-Nash equilibrium, within the class of decentralized policies, for the corresponding $N$-player game. In other words, no player can gain more than $\epsilon$ by switching to another policy depending only on time and on her own state.
\end{theorem}

	\begin{example}[A congestion reward covered by Theorem~\ref{thm:approx-Nash-mf}]
		Consider a finite-horizon model with $\states=\{s_0\}$ and $\actions=\{0,1\}$. Let the transition be independent of the mean field, for instance $P_t(s_0|s_0,a,\nu)=1$ for every $t$, $a$, and $\nu$. Let
		\[
			r_t(s_0,a,\nu)=u_t(a)-\lambda \nu(s_0,a),
			\qquad \lambda\geq 0,
		\]
		where $u_t(a)$ is a baseline reward. The map $\nu\mapsto r_t(s_0,a,\nu)$ is continuous, so the assumptions of Theorem~\ref{thm:approx-Nash-mf} are satisfied. Thus any MFNE policy of this model yields an approximate Nash equilibrium for the corresponding $N$-player game, within decentralized deviations, when $N$ is large enough.
	\end{example}

	\section{Potential MFGs}
	\label{sec:potential-mfg}

	A useful special class of MFGs is given by \boldindex{potential mean field games}. To state the idea in a simple discrete setting, assume that the transition kernel does not depend on the mean field and that the reward can be decomposed as
	\[
		r_t(s,a,\nu_t)=r^0_t(s,a)+F_t(s,a,\nu_t).
	\]
	We say that the coupling is \boldindex{potential}\index{potential coupling} if, after identifying $\nu_t$ with a vector in $\mathbb{R}^{|\states||\actions|}$, there exists a differentiable \boldindex{potential function} $\mathcal F_t$\notationindex{Ftpotential}{$\mathcal F_t$: potential function} such that
	\[
		F_t(s,a,m)=\frac{\partial \mathcal F_t}{\partial m(s,a)}(m),
		\qquad m\in\Delta(\states\times\actions).
	\]
	In this case, the effect of a unilateral infinitesimal change in the representative player's occupation measure is aligned with the variation of the scalar potential $\sum_t \mathcal F_t(\nu_t)$.

	This structure is useful because equilibrium computation can sometimes be related to an optimization problem over occupation measures. Potential MFGs are the mean field analogue of potential games~\citep{monderer1996potential}, and they help explain why optimization algorithms such as fictitious play, Frank-Wolfe-type methods, and smoothed policy iteration can have convergence guarantees in special cases; see, for example,~\citep{hadikhanloo2017learning,geist2022concave,tang2024learning}.
	Analytical frameworks for dynamic potential games, the finite-agent counterparts of potential MFGs, are developed in~\citep{guozhang2025}.

	\begin{example}[A potential MFG with local congestion]
		Let $\lambda_t>0$ and consider the coupling
		\[
			F_t(s,a,m)=-\lambda_t m(s,a),
			\qquad m\in\Delta(\states\times\actions).
		\]
		Then the coupling is potential with
		\[
			\mathcal F_t(m)
			=
			-\frac{\lambda_t}{2}
			\sum_{(s,a)\in\states\times\actions} m(s,a)^2.
		\]
		Indeed,
		\[
			\frac{\partial \mathcal F_t}{\partial m(s,a)}(m)
			=
			-\lambda_t m(s,a)
			=
			F_t(s,a,m).
		\]
		This example corresponds to a local congestion term: using a state-action pair with larger population mass lowers the reward at that same state-action pair.
	\end{example}

\section{Bibliographic notes}

The framework of MFGs was introduced by~\citep{lasry2006jeux,lasry2006jeux2,lasry2007mean}; the idea of finding approximate Nash equilibria using a limiting system was introduced around the same time and called Nash certainty equivalence by~\citep{huang2005nash,huang2006nash,huang2006large}. Lions extensively developed the mathematical theory in his lectures at Coll{\`e}ge de France~\citep{Lions2012CDF}. We refer to the monographs~\citep{bensoussan2013mean} and~\citep{carmona2018probabilistic} for more background in continuous time. The connection with finite-agent games relies on the propagation of chaos; see~\citep{sznitman1991topics}.

MFGs in discrete time and space were first studied in~\citep{gomes2010discrete}. Counterexamples to the convergence of finite-population Nash equilibria to mean field equilibria have been studied in~\citep{doncel2019discrete}; see also~\cite{cardaliaguet2020example} for examples in continuous time and space MFGs.

The terminology of decentralized policies was used e.g. in~\citep{huang2006large,huang2010social} to describe policies that are functions of time and the individual state only (and not functions of the mean field or any other player's state). This is based on the idea that, in a very large system, each player can only observe their own state. If the policy depends on the mean field, it requires information about the population distribution. A different line of work studies population-dependent policies, including master policies, which take the population distribution as an input and aim to generalize across initial distributions or population flows~\citep{perrin2022generalization,wu2025population,lambrecht2026populationaware}. This direction is outside the basic setting of these notes, but is useful for scalable learning algorithms when the initial distribution varies or when common noise makes the population flow random.

The problem of convergence of finite-population Nash equilibria to mean field Nash equilibria has been well studied in continuous time, see e.g.~\citep{lacker2016general,cardaliaguet2019master,lacker2020convergence}. In discrete time, the question has been less explored and we refer the interested reader to e.g.~\citep{doncel2019discrete} for results in this direction.

Models of crowd motion have been studied through the lens of MFGs in~\citep{lachapelle2011mean,achdou2019meancrowd} in the continuous time and space setting. In the discrete setting, models similar to the one presented above have been studied for instance in~\citep{lauriere22scalable-pmlr,perrin2022generalization}.

MFG models for epidemics were studied in the continuous time and finite state setting e.g. in~\citep{laguzet2015individual,hubert2018nash}, with a special focus on vaccination policies.
The SIS  model for epidemics presented in this chapter has been proposed in ~\citep{cui2021approximately}, and  further studied for instance in~\citep{lauriere22scalable-pmlr,guo2024mf}.

The Aiyagari model mentioned in Section~\ref{sec:motivating-ex} was introduced in~\citep{aiyagari1994uninsured}. Other early works in economics considered dynamic games with a continuum of anonymous players, such as~\citep{jovanovic1988anonymous}. Since then, the Aiyagari model and several other macro-economics models have been revisited in the framework of MFGs e.g. in~\citep{achdou2014partial,achdou2022income}. 

The example of finding the starting time of a meeting has been introduced and analyzed in~\citep{gueant2011mean}, which also contains other motivating examples.

The range of applications of MFGs is much broader than the examples discussed in this chapter. Discrete and learning-oriented models have been used for traffic routing~\citep{cabannes2022solving} and flocking~\citep{Perrin-2020-MFGflock}. Other application domains include finance, for instance inter-bank borrowing and lending~\citep{angiuli2022reinforcement}, and engineering applications surveyed in~\citep{barreiro2021mean}.

The fact that mean field equilibrium policies provide approximate equilibrium policies in finite-population games is fundamental to motivate the development of methods to solve MFGs. In line with Theorem~\ref{thm:approx-Nash-mf},
we refer for instance to~\citep[Theorem 4.1]{saldi2018markov}  and~\citep[Theorem 1]{cui2021approximately} for similar results in the case where the interactions are through the state distribution only. See also~\citep{cont2025homogenization} for results in heterogeneous games, with rates of convergence. The present notes focus on homogeneous anonymous MFGs, where players share the same primitives and interact through empirical distributions. Several extensions relax one or both restrictions, including multi-population MFGs, graphon MFGs, major-minor MFGs, Stackelberg MFGs, and models with cooperative subgroups; see~\citep{lauriere2025overview} for a pedagogical overview. For learning algorithms in discrete-time major-minor MFGs, see~\citep{cui2024majorminor}.

For more discussions on the various existing settings of discrete time MFGs, see the survey~\citep{lauriere2022learningsurvey}.

\chapter{Theoretical Foundations of Discrete Mean Field Games}
\label{chap:analytical-MFG}

\section{Fixed point formulation of MFGs} \label{sec:MFG}
In this section, we will present the most established  theory  for
MFGs in a discrete-time framework. For ease of exposition, we will focus on a finite-time horizon $T$.

\paragraph{Iterative solution approach for Nash equilibrium.}
To analyze MFGs, we will focus on the solution concept of the Nash equilibrium, introduced in Definition~\ref{def:mfNE}.
Note that this definition imposes a \boldindex{consistency condition}: the mean-field flow used in the representative agent's control problem must coincide with the flow generated when the population uses the resulting policy.

Based on this definition of Nash equilibrium, the fixed-point formulation first fixes a mean-field flow $(\nu_t)_{t\in\mathcal{T}}$. The representative player's objective is then to maximize her accumulated rewards, \textit{i.e.,} to solve the following MDP, denoted MDP$(\nu)$:
\begin{equation}\label{setup:MFG-MDP-nu}
	\begin{split}
		\text{maximize}_{(\pi_t)_{t\in\mathcal{T}}}  \quad &\mathbb{E}\left[\sum_{t=0}^T r_t(s_t, a_t,\nu_t)\Big|s_0\sim \mu_0\right] \\
\text{subject to} \quad & s_{t+1}\sim P_t(\cdot|s_t,a_t,\nu_t),\, t\in\mathcal{T}\backslash\{T\},\\
                       & a_t\sim \pi_t(\cdot|s_t),\, t\in\mathcal{T}.
	\end{split}
\end{equation}
Once an optimal solution for~\eqref{setup:MFG-MDP-nu} is found, the corresponding state dynamics yields a new mean-field flow $\nu'$. One can then repeat the procedure with MDP$(\nu')$. If this iteration converges, its limit satisfies the consistency condition.
This fixed-point idea is central for both the theoretical analysis and algorithm design of Nash equilibria for MFGs.

To formalize this, let us introduce the following concepts:
\begin{itemize}
	\item Let us denote by $\BR: \Delta(\states\times\actions)^\timesteps \to 2^\strategies$\notationindex{BR}{$\BR$: best response map} the \boldindex{best response map}, such that $\BR(\nu)$ is the set of best responses against $\nu$, i.e., the set of optimal policies for~\eqref{setup:MFG-MDP-nu}.
	\item We denote by $\MF: \strategies \to \Delta(\states)^\timesteps$\notationindex{MF}{$\MF$: state mean field map} the \boldindex{mean field map} defined by $\MF(\pi)=\mu^\pi$, where $\mu^\pi_0=\mu_0$ and, for $t \in \timesteps \backslash \{T\}$,
	      \[
		      \mu_{t+1}^\pi(s') = \sum_{s \in \states} \sum_{a \in \actions} P_t(s'|s, a, \nu_t^\pi) \pi_t(a|s) \mu_t^\pi(s).
	      \]
	      The state-action mean field can be recovered as follows. We denote the \boldindex{extended mean field function} by $\EMF: \strategies \to \Delta(\states\times\actions)^\timesteps$\notationindex{EMF}{$\EMF$: extended mean field function}, which is defined as: $\EMF(\pi) = \nu^\pi$ where $\nu^\pi_t(s,a) = \mu^\pi_t(s) \pi_t(a|s)$, with $\mu^\pi = \MF(\pi)$. In other words,
	      \begin{equation}\label{consistency_pop_flow}
		      \nu_t(s,a)=\mathbb{P}(s_t=s,a_t=a), \qquad (s,a)\in\states\times\actions,\ t\in\mathcal{T},
	      \end{equation}
	      where $s_{t+1}\sim P_t(\cdot|s_t,a_t,\nu_t),$ $t\in\mathcal{T}\backslash \{T\}$, $a_t\sim\pi_t(\cdot|s_t),\, s_0\sim \mu_0,\, t\in\mathcal{T}.$

	\item We denote by $\Gamma: \Pi \to 2^\Pi$\notationindex{Gamma}{$\Gamma$: fixed-point correspondence} the \boldindex{fixed-point correspondence}, i.e., $\Gamma(\pi) = \BR(\EMF(\pi)) = \BR(\nu^\pi)$.
\end{itemize}

This leads to the following fixed-point reformulation of Definition~\ref{def:mfNE}.
\begin{definition}[Mean field Nash equilibrium (rephrased)]\label{def:mfne-fixedpoint}
	A pair $(\pi,\nu) = ((\pi_t)_{t\in\mathcal{T}}, (\nu_t)_{t\in\mathcal{T}})$ of policy and mean field is a mean field Nash equilibrium if the following two conditions are satisfied:
	\begin{itemize}
		\item[1)] (Optimality) $\pi \in \BR(\nu)$,
		\item[2)] (Consistency) $\nu = \EMF(\pi)$.
	\end{itemize}
	$\pi$ is an MFNE policy if $\pi \in \Gamma(\pi)$.
\end{definition}

\section{Existence  of MFNE}
\label{sec:existence-mfne}

In this section, we first present an existence and uniqueness result for MFNE based on the classical Banach fixed-point theorem, under a uniqueness assumption for the best response. Although for finite-state MFGs with policies, best responses are generally not unique, this result highlights the fixed-point approach for MFGs. When the best response is not unique, we then present an existence result based on Kakutani's fixed-point theorem for set-valued maps.

We use the sup norms
\[
	\|\nu\|_\infty=\max_{t\in\timesteps,s\in\states,a\in\actions}|\nu_t(s,a)|,
	\qquad
	\|\pi\|_\infty=\max_{t\in\timesteps,s\in\states,a\in\actions}|\pi_t(a|s)|.
\]
These norms induce the metric spaces used below; since the spaces are finite-dimensional, any equivalent norm would give the same continuity statements.
\begin{theorem}[Existence by Banach fixed point]
	\label{thm:existenceMFNE-Banach}
	Assume that $P$ is Lipschitz continuous in the mean field argument with constant $L_P$, for the sup norm defined above. Assume also that for each mean field there is a unique best response, so that $\BR$ can be viewed as a map from $\Delta(\states\times\actions)^\timesteps$ to $\strategies$. Further assume that $\BR$ is Lipschitz continuous with constant $L_{\BR}$. Let
	\[
		L_{\EMF}:=\sum_{k=0}^{T}(L_P+|\states||\actions|)^k.
	\]
	If $L_{\EMF}L_{\BR}<1$, then there exists a unique MFNE.
\end{theorem}
The proof is given in Appendix~\ref{Appendix: MFNEexistence-Banach}. This result provides existence and uniqueness at once. However, the assumption is  restrictive as in general the best response function may not  be Lipschitz continuous, or the smallness condition is not satisfied.

The following existence result for MFNE holds as long as the transitions and rewards are continuous in $\nu$. It is established by the Kakutani fixed-point theorem. 
\begin{theorem}
	[Existence of MFNE]
	\label{thm:existenceMFNE-Kakutani}
	Suppose that, for every $t\in\mathcal{T}\backslash\{T\}$, $s,s'\in\mathcal{S}$, and $a\in\mathcal{A}$, the map $\nu \mapsto P_t(s'|s,a,\nu)$ is continuous, and that, for every $t\in\mathcal{T}$, $s\in\mathcal{S}$, and $a\in\mathcal{A}$, the map $\nu \mapsto r_t(s,a,\nu)$ is continuous. Then there exists at least one mean field Nash equilibrium.
\end{theorem}

The proof is provided in Appendix~\ref{Appendix: MFNEexistence-Kakutani}.

	\begin{example}[Simple cases covered by the existence assumptions]
		The contraction result of Theorem~\ref{thm:existenceMFNE-Banach} applies in particular to dominant-action models. For instance, take $\states=\{s_0\}$, $\actions=\{0,1\}$, $P_t(s_0|s_0,a,\nu)=1$, and
		\[
			r_t(s_0,1,\nu)=1,
			\qquad
			r_t(s_0,0,\nu)=0.
		\]
		The unique best response is to choose action $1$ at every time, independently of $\nu$. Hence $\BR$ is a constant map, so it is Lipschitz with constant $L_{\BR}=0$, and the smallness condition in Theorem~\ref{thm:existenceMFNE-Banach} is satisfied.

		The Kakutani result of Theorem~\ref{thm:existenceMFNE-Kakutani} covers a broader class. For example, let
		\[
			P_t(\cdot|s,a,\nu)
			=
			(1-\eta)P_t^0(\cdot|s,a)
			+
			\eta\sum_{(\bar s,\bar a)\in\states\times\actions}
			\nu(\bar s,\bar a)P_t^1(\cdot|s,a,\bar s,\bar a),
		\]
		where $\eta\in[0,1]$ and $P_t^0(\cdot|s,a)$ and $P_t^1(\cdot|s,a,\bar s,\bar a)$ are probability distributions on $\states$. Let also
		\[
			r_t(s,a,\nu)
			=
			r_t^0(s,a)
			+
			\sum_{(\bar s,\bar a)\in\states\times\actions}
			c_t(s,a,\bar s,\bar a)\nu(\bar s,\bar a).
		\]
		Then $P_t$ and $r_t$ are continuous in $\nu$, and the assumptions of Theorem~\ref{thm:existenceMFNE-Kakutani} are satisfied.
	\end{example}

\section{Uniqueness of MFNE}
\label{sec:uniqueness-mfne}
One of the most established conditions for uniqueness of an MFNE is the following. For a reward function $r: \timesteps \times \states \times \actions \times \Delta(\states\times\actions) \to \RR$, $t \in \timesteps$, and two mean fields $\nu, \nu' \in \Delta(\states\times\actions)$, we use the notation
\[
	\langle r_t(\nu), \nu' \rangle := \sum_{(s,a) \in \states \times \actions} r_t(s,a,\nu) \nu'(s,a),
\]
which is the dot product of $r_t(\nu)$ and $\nu'$ viewed as vectors indexed by $(s,a)$.

\begin{definition}
	\label{def:LL-mono}
	The MFG satisfies the \boldindex{Lasry-Lions monotonicity assumption} if the following two conditions are satisfied:
	\begin{itemize}
		\item  the transition kernel $P$ is independent of the mean field, and
		\item the reward is such that for all distinct $\nu,\nu' \in \Delta(\states \times \actions)$, for all $t \in \timesteps$,
		      \(\langle r_t(\nu) - r_t(\nu'), \nu-\nu' \rangle\)
		      is strictly negative.
	\end{itemize}
\end{definition}

	\begin{example}[Local congestion and Lasry-Lions monotonicity]
		Assume that the transition kernel is independent of the mean field and that, for each $t\in\timesteps$, $\lambda_t>0$ and
		\[
			r_t(s,a,\nu)=r_t^0(s,a)-\lambda_t \nu(s,a).
		\]
		Then, for any $\nu,\nu'\in\Delta(\states\times\actions)$,
		\[
			\begin{aligned}
				\langle r_t(\nu)-r_t(\nu'),\nu-\nu'\rangle
				 & =
				-\lambda_t
				\sum_{(s,a)\in\states\times\actions}
				\big(\nu(s,a)-\nu'(s,a)\big)^2 \\
				 & =
				-\lambda_t\|\nu-\nu'\|_2^2.
			\end{aligned}
		\]
		Hence the Lasry-Lions monotonicity condition holds whenever $\nu\neq\nu'$. The sign corresponds to congestion: increasing the mass at a state-action pair lowers the reward for using that state-action pair.
	\end{example}

\begin{theorem}[Uniqueness of MFNE]
	\label{thm:unique-mfne-monotonicity}
	Assume that the MFG satisfies the Lasry-Lions monotonicity condition. Then any two MFNEs have the same mean-field flow. In particular, the MFNE is unique at the level of $\nu$; if, in addition, the best response to this flow is unique, then the MFNE policy is unique.
\end{theorem}
The proof is provided in Appendix~\ref{app:auxi-uniqueness-mfne}.

\section{Example of multiple MFNE} \label{Example: multipleNE}
In this section, we present a simple example of an MFG with multiple Nash equilibria. In this game,  the Lasry-Lions monotonicity condition breaks down.

Consider a model of behaviors for a group of animals that try to avoid their predators, via a discrete-time MFG with  finite states and actions.
Let
$\mathcal{S}=\mathcal{A}=\{1,\dots,n\}$ with $n\ge 2$, let $T\ge 1$, assume an arbitrary initial distribution $\mu_0$, and fix reward parameters $r^1,\dots,r^n\ge 0$, at least one of which is strictly positive.
For brevity, write $\delta_{(i,i)}$ for the Dirac mass on the state-action pair $(i,i)$.
At time $0$, for any $\nu_0\in\Delta(\mathcal{S}\times\mathcal{A})$  and any $i,j=1,\dots,n$, the rewards $r_0(i,j,\nu_0)=0$, and the transitions are deterministic
such that
$P_0(i'|i,j,\nu_0)=\mathbf{1}_{\{i'=j\}}$.
For $t=1,\dots,T$, the rewards depend on $\nu_t$ as follows; for $t=1,\dots,T-1$, the transition probabilities are defined similarly. If the population is concentrated on $(i,i)$, \ie, if $\nu_t=\delta_{(i,i)}$, then $r_t(i,j,\nu_t)=r^i\mathbf{1}_{\{i=j\}}$ and, for $t<T$, $P_t(i'|i,j,\nu_t)=\mathbf{1}_{\{i'=j\}}$. Intuitively, staying with the group is safest for an individual agent. Away from this concentrated flow, we use the following formulas: 

\begin{itemize}
	\item $r_t(i,j,\nu_t)=\mathbf{1}_{\{i=j\}}(r^i-r^i\|\nu_t-\delta_{(i,i)}\|_2^2/2)$, for any $i,j=1,\dots,n$ and  $t=1,\dots,T$.
	\item $P_t(i'|i,j,\nu_t)=\dfrac{\mathbf{1}_{\{i'=j\}}+C_{t}\|\nu_t-\delta_{(i,i)}\|_2^2}{1+nC_{t}\|\nu_t-\delta_{(i,i)}\|_2^2}$ for any $i',i,j=1,\dots,n$ and $t=1,\dots,T-1$.
\end{itemize}
Here $C_t$ $(t=1,\dots,T-1)$ are some non-negative constants.

One can verify that for any
$j^\star\in\text{argmax}_{i=1,\dots,n}r^i$, 
if we define \begin{itemize}
	\item $\pi_t(a|s)=\mathbf{1}_{\{a=j^\star\}}$ for any $s,a=1,\dots,n$ and $t=0,\dots,T$;
	\item $\nu_0(s,a)=\mu_0(s)\pi_0(a|s)$ and $\nu_t(s,a)=\mathbf{1}_{\{s=j^\star,a=j^\star\}}$, for any $s,a=1,\dots,n$ and $t=1,\dots,T$,
\end{itemize}
then $\pi=(\pi_t)_{t\in\mathcal{T}}$ and $\nu=(\nu_t)_{t\in\mathcal{T}}$ constitute a Nash equilibrium.

To see the optimality, fix this mean-field flow. At terminal time $T$, action $j^\star$ is optimal: at state $j^\star$ it yields reward $r^{j^\star}$, while at any other state all actions yield zero reward. For $1\leq t<T$, the same conclusion follows by backward induction. If the representative agent is at $j^\star$, action $j^\star$ gives the only positive immediate reward and keeps the agent at $j^\star$. If the representative agent is at another state, all immediate rewards are zero under the flow $\delta_{(j^\star,j^\star)}$, and action $j^\star$ maximizes the transition probability of moving to $j^\star$, the state with maximal continuation value. At time $0$, all immediate rewards are zero, and action $j^\star$ moves the agent deterministically to $j^\star$, which has maximal continuation value. Thus the policy defined above is a best response. Consistency follows from the transition rule at time $0$ and from the deterministic transition at $\nu_t=\delta_{(j^\star,j^\star)}$ for $1\le t\le T-1$.

	Note that the Lasry-Lions monotonicity condition fails. Indeed, for any $t=1,\dots,T$, if $i\neq j$ are such that $r^i+r^j>0$ and $\nu=\delta_{(i,i)}$, $\nu'=\delta_{(j,j)}$, then
	\[
		\langle r_t(\nu)-r_t(\nu'),\nu-\nu'\rangle = r^i+r^j>0,
	\]
	whereas Definition~\ref{def:LL-mono} would require this quantity to be strictly negative. In particular, in the one-step case $T=1$, if the maximizer of $i\mapsto r^i$ is not unique, then the construction above yields multiple distinct Nash equilibria, indexed by $j^\star\in \text{argmax}_{i=1,\dots,n}r^i$ and denoted by $\pi^{(j^\star)}$ and $\nu^{(j^\star)}$.

\section{Characterization of all MFNE}
\label{sec:mfne-characterization}

When an MFG admits multiple Nash equilibria, it is natural to ask whether one can characterize the whole equilibrium set. Such a characterization is useful for studying the robustness of equilibrium strategies, for instance in learning environments or under model uncertainty. The approach presented below relies on an optimization framework for MFGs.
The key idea is to translate the problem of finding Nash equilibria in discrete-time MFGs into an equivalent optimization problem. 
This characterization of the full set of Nash equilibria is established without assuming uniqueness of the Nash equilibrium, contractivity, or the monotonicity conditions discussed in the previous section.

The first step addresses optimality in the representative-agent problem MDP$(\nu)$. It follows the classical linear programming reformulation of MDPs, recalled in Section~\ref{sec:LP-MDP}, by introducing the occupation measure of the representative player.
The second step imposes consistency by setting this occupation measure equal to the mean field. Dual variables and strong duality then express the optimality condition through primal-dual constraints. Together, these conditions characterize Nash equilibria as a constrained optimization problem.
This constrained formulation can then be refined into an optimization problem with bounded variables and convex constraints, which is more convenient computationally.

\section{Primal-dual formulation for Nash equilibria}\index{primal-dual formulation}
\label{sec:maths-MFG}

We work with the finite-horizon MFG setup of Section~\ref{sec:MFG}, with state space $\mathcal{S}$ and action space $\mathcal{A}$ satisfying $|\mathcal{S}|=S<\infty$ and $|\mathcal{A}|=A<\infty$.
For notational simplicity, for a given mean-field flow $\nu=(\nu_t)_{t\in\mathcal{T}}$, we denote by $\mathcal{M}(\nu)$ the mean-field-induced MDP~\eqref{setup:MFG-MDP-nu}. We write $V_t^\star(\nu)\in\mathbb{R}^S$ for its optimal value function from time $t$, with $[V_t^\star(\nu)]_s$ the optimal expected reward starting from state $s$ at time $t$, and $V_{\mu_0}^\star(\nu)=\sum_{s\in\mathcal{S}}\mu_0(s)[V_0^\star(\nu)]_s$ for the corresponding value under the initial distribution $\mu_0$.
For a policy $\pi=(\pi_t)_{t\in\mathcal{T}}$, we similarly write $V_t^{\pi}(\nu)\in\mathbb{R}^S$, $[V_t^{\pi}(\nu)]_s$, and $V_{\mu_0}^{\pi}(\nu)=\sum_{s\in\mathcal{S}}\mu_0(s)[V_0^{\pi}(\nu)]_s$ for the values in $\mathcal{M}(\nu)$ when the representative agent uses $\pi$.

By the occupation measure formulation of Proposition~\ref{prop:MDP_LP}, a pair $(\pi,\nu)$ is a Nash equilibrium if and only if there exists an occupation measure $d^\star=(d_t^\star(s,a))_{t\in\mathcal{T}}$ such that $\pi\in\Poldist(d^\star)$ and the following two conditions hold.

\begin{itemize}
	\item[\bf (A) ] Optimality of the representative agent: $d^\star=(d_t^\star(s,a))_{t\in\mathcal{T}}$ solves the following linear program
		\begin{align}
			\text{maximize}_{d}\quad & \sum_{t\in\mathcal{T}}\sum_{s\in\mathcal{S}}\sum_{a\in\mathcal{A}} d_t(s,a)r_t(s,a,\nu_t)                  \\
			\text{subject to} \quad  & \sum_{s\in\mathcal{S}}\sum_{a\in\mathcal{A}}d_t(s,a)P_t(s'|s,a,\nu_t)=\sum_{a\in\mathcal{A}}d_{t+1}(s',a),
			\label{cons_occ}                                                                                                                      \\ &\hspace{4.2cm}\forall s'\in\mathcal{S},t\in\mathcal{T}\backslash\{T\},\notag\\
			                         & \sum_{a\in\mathcal{A}}d_0(s,a)=\mu_0(s), \quad\forall s\in\mathcal{S},\notag                               \\
			                         & d_t(s,a)\geq 0,\quad \forall s\in\mathcal{S},a\in\mathcal{A},t\in\mathcal{T}.\notag
		\end{align}

	\item[\bf (B)] Consistency between the agent and the mean-field:
		$d_t^\star(s,a)=\nu_t(s,a)$ for any $s\in\mathcal{S},a\in\mathcal{A},t\in\mathcal{T}$.
\end{itemize}
Condition (A) encodes optimality of the representative agent against the fixed mean-field flow $\nu$. Condition (B) encodes consistency: the occupation measure generated by the representative agent matches the population state-action flow.

Equivalently, the linear program in condition (A) can be written as
\begin{equation}\label{cond_A_rewrite}
	\begin{array}{ll}
		\text{minimize}_d & c_\nu^\top d            \\
		\text{subject to} & A_\nu d=b,\quad d\geq0,
	\end{array}
\end{equation}
where $c_\nu\in\mathbb{R}^{SA(T+1)}$ is
\begin{equation}\label{kkt_param1}
	c_\nu=\left[
		\begin{array}{c}
			-r_0(\cdot,\cdot,\nu_0) \\
			\vdots                  \\
			-r_T(\cdot,\cdot,\nu_T)
		\end{array}
		\right],\quad
	b=\left[
		\begin{array}{c}
			0      \\
			\vdots \\
			0      \\
			\mu_0
		\end{array}
		\right],
\end{equation}
$r_t(\cdot,\cdot,\nu_t)\in\mathbb{R}^{SA}$ is a flattened vector (with column-major order),  
and the matrix $A_\nu\in\mathbb{R}^{S(T+1)\times SA(T+1)}$ is
\begin{equation}\label{kkt_param2}
	A_\nu=\left[
		\begin{array}{ccccccc}
			W_0(\nu_0) & -Z         & 0          & 0      & \cdots & 0                  & 0      \\
			0          & W_1(\nu_1) & -Z         & 0      & \cdots & 0                  & 0      \\
			0          & 0          & W_2(\nu_2) & -Z     & \cdots & 0                  & 0      \\
			\vdots     & \vdots     & \vdots     & \vdots & \ddots & \vdots             & \vdots \\
			0          & 0          & 0          & 0      & \cdots & W_{T-1}(\nu_{T-1}) & -Z     \\
			Z          & 0          & 0          & 0      & \cdots & 0                  & 0
		\end{array}\right].
\end{equation}
Here, 
$W_t(\nu_t)\in \mathbb{R}^{S\times SA}$ is the matrix with the $l$-th row ($l=1,\dots,S$) being the  flattened vector  $[P_t(l|\cdot,\cdot,\nu_t)]\in\mathbb{R}^{SA}$ (with column-major order),
and the matrix $Z$ is defined as
\begin{equation}\label{eq:defn-Z}
	Z:=[\overbrace{I_{S},\dots,I_{S}]}^{A}\in \mathbb{R}^{S\times SA},
\end{equation}
where $I_S$ is the $S\times S$ identity matrix. The variable $d$ is flattened in the same order as $c_\nu$. In what follows, $d^\star=(d^\star_t)_{t\in\mathcal{T}}$ and $\nu=(\nu_t)_{t\in\mathcal{T}}$ are viewed either as vectors in $\mathbb{R}^{SA(T+1)}$ or as sequences of $(T+1)$ vectors in $\mathbb{R}^{SA}$, always using the same column-major ordering. We use $\nu_t(s,a)$ and $\nu_{s,a,t}$, and similarly $d_t(s,a)$ and $d_{s,a,t}$, interchangeably.

By strong duality for linear programs~\citep{luenberger1984linear}, $d^\star$ is an optimal solution of~\eqref{cond_A_rewrite} if and only if there exist \boldindex{dual variables} $y$ and \boldindex{slack variables} $z$ such that the following \boldindex{Karush--Kuhn--Tucker conditions} (\boldindex{KKT conditions}) hold:
\begin{equation}\label{eq:kkt}
	\begin{split}
		&A_\nu d^\star=b,\quad A_\nu^\top y+z=c_\nu,\\
		& z^\top d^\star=0,\quad d^\star\geq 0,\quad z\geq 0.
	\end{split}
\end{equation}

Combining condition (A), written through~\eqref{eq:kkt}, with the consistency condition (B), we obtain the following characterization.
\begin{theorem}  \label{thm:nash-opt}
	Finding a Nash equilibrium of the MFG in Section~\ref{sec:MFG}  is equivalent to solving the following
	constrained optimization problem:
	\begin{equation}\label{mf-vmo-constr}
		\begin{array}{ll}
			\text{\rm minimize}_{y,z,\nu} \quad & 0                                                                                                \\
			\text{\rm subject to}               & A_\nu \nu=b,\quad A_\nu^\top y+z=c_\nu,                                                          \\
			                                    & z^\top \nu=0,\quad \nu_t\geq 0,\quad \mathbf{1}^\top\nu_t=1,\quad t\in\mathcal{T},\quad z\geq 0,
		\end{array}
	\end{equation}
	where $A_\nu$, $c_\nu$, and $b$ are defined in~\eqref{kkt_param1} and~\eqref{kkt_param2}.
\end{theorem}

The constrained feasibility problem \eqref{mf-vmo-constr} in Theorem~\ref{thm:nash-opt}
can be rewritten in a form that is more convenient for first-order optimization methods: the auxiliary variables $y$ and $z$ are bounded, and the equality and complementarity constraints are moved to the objective.
The remaining constraints are simple convex constraints.

\begin{theorem}  \label{thm:mfvmo-full}
	Finding a Nash equilibrium of the MFG in Section~\ref{sec:MFG} is equivalent to solving the following optimization problem, which we call \boldindex{Mean-Field Occupation Measure Optimization} (\boldindex{MF-OMO}):
	\begin{equation}\label{mfvmo-full}
		\begin{array}{ll}
			\text{minimize}_{y,z,\nu} & \|A_\nu \nu-b\|_2^2+\|A_\nu^\top y+z-c_\nu\|_2^2+z^\top \nu    \\
			\text{subject to}         & \nu_t\geq 0, \quad \mathbf{1}^\top \nu_t=1, \,t\in\mathcal{T}, \\
			                          & \mathbf{1}^\top z\leq SA(T+1)(T+2)r_{\max}, \quad z\geq 0,      \\ 
			                          & \|y\|_2\leq \frac{S(T+1)(T+2)}{2}r_{\max},
		\end{array}
	\end{equation}
	where $A_\nu$, $c_\nu$, and $b$ are defined in~\eqref{kkt_param1} and~\eqref{kkt_param2}, ${\bf 1}$ is the all-one vector with the appropriate dimension, and
	the reward bound
\[
  r_{\max}=\sup_{\substack{t\in\mathcal{T},\ s\in\mathcal{S},\ a\in\mathcal{A}\\
  \nu\in\Delta(\mathcal{S}\times\mathcal{A})}}|r_t(s,a,\nu)|
\]
is assumed to be finite.
\end{theorem}
	Note that in this theorem, the term $\|A_\nu \nu-b\|_2^2$ measures the flow/consistency residual, the term $\|A_\nu^\top y+z-c_\nu\|_2^2$ measures the Bellman-dual residual, and $z^\top \nu$ is the complementarity gap. The feasible set is convex and compact, but the objective is generally nonconvex because $A_\nu$ and $c_\nu$ depend on $\nu$. Existence of an MFNE is equivalent to the global minimum value being zero; positive local minima or stationary points may still occur.

The theorem says that if $(\pi,\nu)$ is a Nash equilibrium of the MFG, then there exist $y,z$ such that $(y,z,\nu)$ is feasible for~\eqref{mfvmo-full} and attains objective value $0$. Conversely, if a feasible point $(y,z,\nu)$ for~\eqref{mfvmo-full} attains objective value $0$, then for any $\pi\in \Poldist(\nu)$, the pair $(\pi,\nu)$ is a Nash equilibrium.
We denote the objective function in \eqref{mfvmo-full} by $f^{\text{\rm MF-OMO}}(y,z,\nu)$.

More generally, approximate Nash equilibria can be related to the objective value in~\eqref{mfvmo-full}, without assuming Lipschitz continuity of the dynamics and rewards in $\nu_t$.
	For a policy $\pi$, we define its \boldindex{exploitability} by
	\begin{equation}
		\label{eq:def-exploitability}
		\Expl(\pi)
		=
		\max_{\pi'\in\strategies} J(\pi';\nu^\pi)-J(\pi;\nu^\pi),
	\end{equation}
	\notationindex{Expl}{$\Expl$: exploitability}
	where $\nu^\pi$ is the mean field generated by $\pi$ and the population flow is kept fixed at $\nu^\pi$ while a representative player deviates to $\pi'$.

For $a\in\mathcal{A}$ and $t<T$, let $P_t^a(\nu_t)\in\mathbb{R}^{S\times S}$ denote the transition matrix with entries $[P_t^a(\nu_t)]_{s,s'}=P_t(s'|s,a,\nu_t)$.
\begin{theorem}\label{MFNE-vs-VMO}
	Let $\pi=(\pi_t)_{t\in\mathcal{T}}$ be an $\epsilon$-Nash equilibrium, \ie, $\Expl(\pi)\leq\epsilon$. Define $\nu=\nu^\pi=\EMF(\pi)$,
	\[
		y=[V_1^\star(\nu), V_2^\star(\nu), \dots,V_T^\star(\nu), -V_0^\star(\nu)]\in\mathbb{R}^{S(T+1)},
	\]
	and $z=[z_0,z_1,\dots,z_T]\in\mathbb{R}^{SA(T+1)}$ with $z_t=[z_t^1,\dots,z_t^A]$ ($t=0,\dots,T$), where
	\[
		z_t^a=V_t^\star(\nu)-r_t(\cdot,a,\nu_t)-P_t^a(\nu_t)V_{t+1}^\star(\nu)\in\mathbb{R}^S, \quad t=0,\dots,T-1,\,a\in\mathcal{A}
	\]
	and $z_T^a=V_T^\star(\nu)-r_T(\cdot,a,\nu_T)\in\mathbb{R}^S$ ($a\in\mathcal{A}$). Then $(y,z,\nu)$ is a feasible solution to \eqref{mfvmo-full} and $f^{\text{\rm MF-OMO}}(y,z,\nu)\leq \epsilon$.
\end{theorem}

\paragraph{Special cases.}
For MFGs with linear rewards and dynamics independent of the mean field, the optimization formulation becomes more explicit.

Specifically, consider MFGs with
	\[
		\begin{aligned}
			r_t(s,a,\nu_t)    & =\bar{r}_{s,a,t}+\nu_t^\top\bar{R}_{s,a,t}, &  & t\in\mathcal{T}, \\
			P_t(s'|s,a,\nu_t) & =\bar{P}_{s',s,a,t},                        &  & t=0,\dots,T-1,
		\end{aligned}
	\]
	for constants $\bar{r}_{s,a,t}\in\mathbb{R}$ and $\bar{R}_{s,a,t}\in\mathbb{R}^{SA}$ for $t\in\mathcal{T}$, and transition probabilities $\bar{P}_{s',s,a,t}$ for $t=0,\dots,T-1$ satisfying $\bar{P}_{s',s,a,t}\ge 0$ and $\sum_{s'\in\states}\bar{P}_{s',s,a,t}=1$.
Then  $A_{\nu}$ is independent of  $\nu$, denoted as $\bar{A}$, and $c_{\nu}$ is of the form $c_{\nu}=\bar{c}+\bar{C} \nu$, where $\bar{c}\in\mathbb{R}^{SA(T+1)}$ is the vectorization of $[-\bar{r}_{s,a,t}]_{s\in\mathcal{S},a\in\mathcal{A},t\in\mathcal{T}}\in \mathbb{R}^{S\times A\times (T+1)}$ (with column-major order), and $\bar{C}\in\mathbb{R}^{SA(T+1)\times SA(T+1)}$ is defined as
\[
	\bar{C}=\operatorname{diag}\left(
	\left[\begin{matrix}
			-\bar{R}_{1,1,0}^\top \\
			\vdots                \\
			-\bar{R}_{S,A,0}^\top
		\end{matrix}
		\right],
	\left[\begin{matrix}
			-\bar{R}_{1,1,1}^\top \\
			\vdots                \\
			-\bar{R}_{S,A,1}^\top
		\end{matrix}
		\right],\dots,
	\left[\begin{matrix}
			-\bar{R}_{1,1,T}^\top \\
			\vdots                \\
			-\bar{R}_{S,A,T}^\top
		\end{matrix}
		\right]
	\right).
\]
Using~\eqref{mf-vmo-constr}, finding Nash equilibria in this class is equivalent to solving the following linear complementarity problem:
\begin{equation*}
	\begin{array}{ll}
		\text{\rm minimize}_{y,z,\nu} \quad & 0                                                                                                   \\
		\text{\rm subject to}               & \left[\begin{array}{ll}
				                                            0            & \bar{A}  \\
				                                            \bar{A}^\top & -\bar{C}
			                                            \end{array}\right]\left[\begin{array}{l}
				                                                                    y \\
				                                                                    \nu
			                                                                    \end{array}\right]+\left[\begin{array}{l}
				                                                                                             0 \\
				                                                                                             z
			                                                                                             \end{array}\right]=\left[\begin{array}{l}
				                                                                                                                      b \\
				                                                                                                                      \bar{c}
			                                                                                                                      \end{array}\right], \\
		                                    & z^\top \nu=0,\quad \nu\geq 0,\quad z\geq 0.
	\end{array}
\end{equation*}
This finite-dimensional complementarity problem can, in principle, be solved by enumerating complementarity patterns and solving the resulting linear programs.

More generally, the connection between Nash equilibria and the feasible set of the optimization formulation makes it possible to use optimization algorithms to find multiple, and sometimes all, Nash equilibria of MFGs.

\section{Bibliographic notes}

	Contraction arguments for existence of mean field Nash equilibria, such as Theorem~\ref{thm:existenceMFNE-Banach}, have been used extensively in the continuous time and space literature at least since~\citep{huang2006large}. In many elementary presentations, these arguments require the time horizon or the relevant Lipschitz constants to be sufficiently small. This is not intrinsic to all contraction-based approaches: in some continuous-time models, additional structural assumptions can yield well-posedness for arbitrary horizons; see, for instance, the weak monotonicity condition with common noise studied by~\citep[Theorem~4.4]{ahuja2016wellposedness}. For existence results based on Kakutani fixed point theorem in finite or discounted settings, we refer to~\citep{gomes2010discrete,saldi2018markov,cui2021approximately,guo2024mf}.

	The Lasry-Lions monotonicity condition of Definition~\ref{def:LL-mono} was introduced by~\citep{lasry2007mean} and is one of the central structural conditions in the analysis of MFGs. It is used not only for uniqueness, but also in stability estimates, convergence of approximation schemes, and convergence of learning algorithms. Without a monotonicity condition, uniqueness can fail, as illustrated by Section~\ref{Example: multipleNE}; see also~\citep{bardi2019nonuniqueness} for continuous finite-horizon examples and~\citep{cecchin2019convergence} for a finite-state model where non-uniqueness affects the convergence problem. Other monotonicity notions have also been developed. In particular, displacement monotonicity plays a role distinct from, and sometimes in dichotomy with, the Lasry--Lions condition: it can yield uniqueness for arbitrary time horizons in first-order MFGs with nonseparable Hamiltonians, including models motivated by crowd motion; see~\citep{meszaros2024mean}.

The technical presentation of the primal-dual representation of the set of Nash equilibria follows closely~\citep{guo2024mf}, from which the example of multiple Nash equilibria in Section~\ref{Example: multipleNE} is borrowed. Extension to personalized MFGs accommodating heterogeneity
can also be formulated using the same primal-dual approach, as analyzed in~\citep{guo2024mf}.
The connection between some MFG systems and variational problems already appears in the original work~\citep{lasry2006jeux2,lasry2007mean}. For first-order MFG systems of partial differential equations (PDEs), weak solutions can be characterized through an optimal control problem constrained by Hamilton-Jacobi and continuity equations~\citep{cardaliaguet2015mean}. Related second-order formulations with local couplings are studied in~\citep{cardaliaguet2015second}. Numerical methods for these PDE-constrained or variational formulations include augmented Lagrangian, proximal, and primal-dual algorithms; see, for example,~\citep{benamou2015augmented,briceno2018proximal,briceno2019implementation}. For broader accounts of numerical methods for continuous-space MFGs and mean field type control, see~\citep{achdoulauriere2020mfgnumerical,lauriere2021numerical}. Neural-network and machine-learning approaches to high-dimensional continuous MFG or mean field control (MFC) systems include~\citep{ruthotto2020machine,carmonalauriere2021convergenceergodic,carmona2022convergence,lin2021alternating,cao2024connecting}. A primal-dual approach for more general continuous-time MFGs is developed in~\citep{Guo2025Primaldual} by removing the convexity of the Hamiltonian and the uniqueness of its minimizer. This approach does not require that the equilibrium flow have a full support.

\chapter{Algorithms}
\label{chap:algorithms}

\section{Fixed-point-based algorithms for MFGs}
Fixed point methods are natural for MFGs since the definition of MFNE can be phrased as a fixed point problem. In this section, we present four types of methods: pure fixed point iterations, fictitious play, policy iteration, and online mirror descent.

\subsection{Fixed point iterations}
\label{sec:algo-fixed-point}

The first approach relies on finding a fixed point of a mapping that characterizes the equilibrium of the MFG. This method typically leverages the Banach fixed-point theorem to ensure convergence under the same conditions as the ones ensuring existence (and uniqueness) in Theorem~\ref{thm:existenceMFNE-Banach}.

We recall Definition~\ref{def:mfne-fixedpoint} and the notations $\BR$ and $\EMF$. To implement a method based on this definition, a key step is to compute a best response, i.e., an optimal policy against a given mean field $\nu = (\nu_t)_{t \in \timesteps}$. This can be achieved by computing the Q-function of~\eqref{setup:MFG-MDP-nu}. It is defined as in~\eqref{eq:MDP-Qstar} but depends on $\nu$ which influences the transitions and the rewards. We thus denote it by $Q^{*,\nu}$. The Bellman equation~\eqref{eq:Qstar-Bellman} becomes:
\begin{equation}
	\label{eq:def-Qstar-nu-T}
	\begin{cases}
		\displaystyle
		Q^{*,\nu}_{T}(s,a) =  r_T(s,a,\nu_T), \quad s \in \states, a \in \actions,
		\\
		\displaystyle
		Q^{*,\nu}_{t}(s,a) = \EE_{s' \sim P_t(\cdot|s,a,\nu_t)}[ r_t(s,a,\nu_t) + \max_{a' \in \actions} Q^{*,\nu}_{t+1}(s',a')], \\
		\qquad\qquad \quad s \in \states, a \in \actions, t = T-1, \dots, 1, 0.
	\end{cases}
\end{equation}
From here, a best response $\pi^*$ against $\nu$ can be found by taking an element in the set $\mathrm{Greedy}(Q^{*,\nu})$ of greedy policies, i.e., a policy such that
$\pi^*_t(\cdot|s)$ is supported in $\argmax Q^{*,\nu}_t(s,\cdot)$ for every $t \in \timesteps$ and $s \in \states$.

The \boldindex{fixed point algorithm} consists in iteratively updating the mean field and the policy. It is summarized in Algorithm~\ref{algo:fixed-point}.

\begin{algorithm}[H]
	\caption{Pure fixed point iterations for MFG}
	\label{algo:fixed-point}
	\begin{algorithmic}[1]
		\Require initial mean field \( \nu^{(0)} = (\nu^{(0)}_t)_{t \in \timesteps}\), number of iterations $K$
		\Ensure approximate mean field Nash equilibrium
		\For{$k=0,\dots,K-1$}
		\State Update the policy $\pi^{(k+1)} \in \BR(\nu^{(k)})$
		\State Update the mean field \( \nu^{(k+1)} = \EMF(\pi^{(k+1)})\)
		\EndFor
		\State \Return \( \pi^{(K)} = (\pi^{(K)}_t)_{t \in \timesteps}\) and \( \nu^{(K)} = (\nu^{(K)}_t)_{t \in \timesteps}\)
	\end{algorithmic}
\end{algorithm}

At iteration $k$, there could be multiple elements in $\mathrm{Greedy}(Q^{*,\nu^{(k)}})$, \ie, multiple best responses to the given mean field $\nu^{(k)}$. As a consequence, an implementation needs to select one policy among them. Once such a selection rule is fixed, the iteration defines a map on mean fields. The most direct convergence mechanism is contraction of this selected map: in that case the iterates converge to its unique fixed point, and the corresponding policy is a mean field Nash equilibrium policy. This condition is restrictive, but it is useful because it identifies what can go wrong when the best-response map is too sensitive to the mean field.

However, in many scenarios, such fixed point iterations can fail to converge. We refer to Section~\ref{sec:crowd-motion-numerics} below for an illustration of the lack of convergence in the context of a grid world model for crowd motion.

Such convergence failure can be due to a lack of regularity of the best response map, i.e., a small change in $\nu$ can induce a drastic change in $\BR(\nu)$. To address this issue, one can replace the $\BR$ map by an approximate best-response. For example, at iteration $k$, after computing $Q^{(k)}=Q^{*,\nu^{(k)}}$, one can replace the greedy update by
\[
	\pi^{(k+1)}_t(a|s)  = \softmax_\tau\left(Q^{(k)}_t(s,\cdot)\right)(a)
\]
for every $(t,s,a) \in \timesteps \times \states \times \actions$, where $\softmax_\tau(q):=\softmax(q/\tau)$ and $\softmax: \RR^d \to \RR^d$ is defined by: for a vector $z = (z_1,\dots,z_d) \in \RR^d$,
\[
	\softmax (z) = \left(\frac{e^{z_i}}{ \sum_{j=1}^d e^{z_j}}\right)_{i=1,\dots,d},
\]
which can be interpreted as a probability distribution on $\{1,\dots,d\}$. This leads to the \boldindex{softmax fixed point algorithm}.

Another way to stabilize the iteration is to damp the updates. For instance, one can update at every iteration the mean field or the policy by taking a weighted average with past mean field or policy. For example, to average the mean fields, one can replace in Algorithm~\ref{algo:fixed-point} the policy update $\pi^{(k+1)} \in \BR(\nu^{(k)})$ by
\[
	\pi^{(k+1),\alpha} \in \BR(\overline\nu^{(k),\alpha}),
\]
where, for all $t \in \timesteps, s \in \states, a \in \actions$,
\begin{equation}
	\label{eq:damped-fixed-point-update-avgmf}
	\overline\nu^{(k+1),\alpha}_t(s,a)
	=
	(1-\alpha) \overline\nu^{(k),\alpha}_t(s,a)
	+ \alpha \nu^{(k+1),\alpha}_t(s,a),
\end{equation}
with $\nu^{(k+1),\alpha} = \EMF(\pi^{(k+1),\alpha})$.
Here $\alpha \in (0,1]$ is an update rate (or learning rate); $(1-\alpha)$ can be viewed as a damping coefficient. This is summarized in Algorithm~\ref{algo:damped-fixed-point}. We refer to this algorithm as \boldindex{damped fixed point algorithm}.

\begin{algorithm}[H]
	\caption{Damped fixed point iterations for MFG}
	\label{algo:damped-fixed-point}
	\begin{algorithmic}[1]
		\Require initial mean field \( \nu^{(0),\alpha} = \overline\nu^{(0),\alpha} = (\nu^{(0)}_t)_{t \in \timesteps} \), update rate $\alpha \in (0,1]$, number of iterations $K$
		\Ensure approximate mean field Nash equilibrium
		\For{$k=0,\dots,K-1$}
		\State Update the policy $\pi^{(k+1),\alpha} \in \BR(\overline\nu^{(k),\alpha})$
		\State Update the mean field \( \nu^{(k+1),\alpha} = \EMF(\pi^{(k+1),\alpha})\)
		\State Update the average mean field $\overline\nu^{(k+1),\alpha} = (1-\alpha) \overline\nu^{(k),\alpha} + \alpha \nu^{(k+1),\alpha}$
		\EndFor
		\State \Return \( \pi^{(K),\alpha} = (\pi^{(K),\alpha}_t)_{t \in \timesteps}\) and \( \nu^{(K),\alpha} = (\nu^{(K),\alpha}_t)_{t \in \timesteps}\)
	\end{algorithmic}
\end{algorithm}

If $\alpha=1$, the algorithm reduces to the ``pure'' fixed point iterations of Algorithm~\ref{algo:fixed-point}. On the other hand, taking $\alpha$ closer to $0$ slows down the updates and can improve stability in cases where the undamped fixed-point iteration oscillates or diverges.

	The stabilizing effect of damping can be made precise when the residual in the consistency condition is strongly monotone. The following proposition is stated for a single-valued policy response. This can arise from a fixed selection of the best response, or from replacing the greedy best response by a softmax response as above. The map $\Phi_\alpha$ below corresponds to the update of the averaged mean field $\overline\nu^{(k),\alpha}$ in Algorithm~\ref{algo:damped-fixed-point}, with $\Phi(\overline\nu^{(k),\alpha})=\nu^{(k+1),\alpha}$.

	\begin{proposition}[Damping with a strongly monotone residual]
		\label{prop:damped-monotone-residual}
		Let $\mathcal X=\Delta(\states\times\actions)^\timesteps$, identified with a finite-dimensional Euclidean space, and let
		\[
			\Phi(\nu)=\EMF(\mathcal R(\nu))
		\]
		be a single-valued mean-field update map, where $\mathcal R(\nu)\in\strategies$ is the policy response used against the mean field $\nu$. Define the consistency residual
		\[
			G(\nu)=\nu-\Phi(\nu).
		\]
		Assume that there exist constants $m>0$ and $M<\infty$ such that, for all $\nu,\nu'\in\mathcal X$,
		\[
			\langle G(\nu)-G(\nu'),\nu-\nu'\rangle
			\geq m\|\nu-\nu'\|_2^2,
		\]
		and
		\[
			\|G(\nu)-G(\nu')\|_2
			\leq M\|\nu-\nu'\|_2.
		\]
		Then, for every update rate $0<\alpha\leq 1$ such that
		\[
			\alpha<\frac{2m}{M^2},
		\]
		the damped fixed-point map
		\[
			\Phi_\alpha(\nu)=(1-\alpha)\nu+\alpha\Phi(\nu)
			=\nu-\alpha G(\nu)
		\]
		is a strict contraction. More precisely,
		\[
			\|\Phi_\alpha(\nu)-\Phi_\alpha(\nu')\|_2
			\leq
			\sqrt{1-2\alpha m+\alpha^2M^2}\,
			\|\nu-\nu'\|_2 .
		\]
		Consequently, the damped iteration converges geometrically to the unique fixed point of $\Phi$.
	\end{proposition}
	\begin{proof}
		For any $\nu,\nu'\in\mathcal X$,
		\[
			\Phi_\alpha(\nu)-\Phi_\alpha(\nu')
			=
			(\nu-\nu')-\alpha(G(\nu)-G(\nu')).
		\]
		Hence
		\[
			\begin{aligned}
				\|\Phi_\alpha(\nu)-\Phi_\alpha(\nu')\|_2^2
				 & =
				\|\nu-\nu'\|_2^2                               \\
				 & \quad
				-2\alpha\langle G(\nu)-G(\nu'),\nu-\nu'\rangle \\
				 & \quad
				+\alpha^2\|G(\nu)-G(\nu')\|_2^2                \\
				 & \leq
				(1-2\alpha m+\alpha^2M^2)\|\nu-\nu'\|_2^2.
			\end{aligned}
		\]
		The factor is strictly smaller than $1$ under the stated condition on $\alpha$, so $\Phi_\alpha$ is a strict contraction. Since $\mathcal X$ is convex and $\Phi(\mathcal X)\subseteq \mathcal X$, we also have $\Phi_\alpha(\mathcal X)\subseteq\mathcal X$. Moreover, because $\alpha>0$, the fixed points of $\Phi_\alpha$ and $\Phi$ coincide. Banach's fixed-point theorem then gives geometric convergence to the unique fixed point.
	\end{proof}
	The first assumption in Proposition~\ref{prop:damped-monotone-residual} has a simple interpretation. The map $G(\nu)=\nu-\Phi(\nu)$ measures the error in the consistency condition. The monotonicity assumption says that this residual points, on average, in a direction that reduces discrepancies between mean-field flows. The Lipschitz assumption controls how sensitive the residual is to changes in the input flow. Under these two assumptions, damping amounts to taking a sufficiently small step in the direction $-G(\nu)$, and the proposition shows that this step is contractive.

	The assumptions above are sufficient conditions, not general properties of damped fixed-point iterations. Damping can control overshooting when the response map already moves discrepancies in a stabilizing direction. But, for instance in one dimension, if the response map has derivative greater than $1$ at a fixed point, no positive damping rate makes the update locally contractive there. The following elementary example illustrates this distinction.

	\begin{example}[Congestion versus herding in a damped softmax update]
		Consider a one-stage model with one state $\states=\{s_0\}$ and two actions $\actions=\{0,1\}$. Since there is only one state, a state-action mean field can be identified with the mass $p\in[0,1]$ choosing action $1$; the mass choosing action $0$ is $1-p$. For a given $p$, suppose that the algorithm uses the softmax policy
		\[
			\pi^p(1|s_0)=B_\tau(p),
			\qquad
			\pi^p(0|s_0)=1-B_\tau(p),
		\]
		where, for a temperature $\tau>0$,
		\[
			B_\tau(p)=
			\frac{\exp(r(1,p)/\tau)}
			{\exp(r(1,p)/\tau)+\exp(r(0,p)/\tau)}.
		\]
		In this one-state model, the mean-field update generated by this policy is $\Phi(p)=B_\tau(p)$. The damped update is therefore
		\[
			D_\alpha(p)=(1-\alpha)p+\alpha B_\tau(p).
		\]
		First consider congestion rewards
		\[
			r(1,p)=a_1-\beta p,\qquad r(0,p)=a_0-\beta(1-p),
			\qquad \beta>0.
		\]
		Then
		\[
			B_\tau'(p)
			=
			-\frac{2\beta}{\tau}B_\tau(p)(1-B_\tau(p)),
		\]
		and hence $-\beta/(2\tau)\le B_\tau'(p)\le 0$. Therefore
		\[
			D_\alpha'(p)=1-\alpha+\alpha B_\tau'(p)
		\]
		lies between $1-\alpha(1+\beta/(2\tau))$ and $1-\alpha$. Thus, if
		\[
			0<\alpha\le 1
			\qquad\hbox{and}\qquad
			\alpha<\frac{2}{1+\beta/(2\tau)},
		\]
		then $|D_\alpha'(p)|<1$ uniformly in $p$, and the mean value theorem shows that $D_\alpha$ is a strict contraction.

		In contrast, for herding rewards
		\[
			r(1,p)=a_1+\beta p,\qquad r(0,p)=a_0+\beta(1-p),
		\]
		the sign is reversed. In the symmetric case $a_1=a_0$, the point $p=1/2$ is fixed and
		\[
			B_\tau'(1/2)=\frac{\beta}{2\tau}.
		\]
		If $\beta>2\tau$, then
		\[
			D_\alpha'(1/2)
			=
			1-\alpha+\alpha\frac{\beta}{2\tau}>1
		\]
		for every $\alpha>0$. Thus damping cannot make the update locally contractive in this case.
	\end{example}

We refer to Section~\ref{sec:crowd-motion-numerics} for a numerical illustration of the damped fixed point algorithm.

\subsection{Fictitious play}
\label{sec:algo-fictitious-play}

In the context of game theory, the \boldindex{fictitious play algorithm} (\boldindex{FP}) is a learning process where each agent assumes the other agents' strategies are fixed and iteratively best responds to the empirical distribution of past strategies. In MFGs, fictitious play incorporates averaging of the mean field distribution.

Given a mean field, the best response map $\BR$ and the extended mean field map $\EMF$ are computed as above. However, at each iteration the best response is computed against the average mean field over past iterations rather than the latest one. The method is summarized in Algorithm~\ref{algo:fictitious-play}. The average mean field at iteration $k$ is computed as:
\begin{equation*}
	\overline\nu^{(k+1)}_t(s,a) = \frac{k}{k+1}\overline\nu^{(k)}_t(s,a) + \frac{1}{k+1}\nu^{(k+1)}_t(s,a),
\end{equation*}
for all $t \in \timesteps, s \in \states, a \in \actions$. This can be viewed as a variant of the damped fixed point~\eqref{eq:damped-fixed-point-update-avgmf} when $\alpha$ is replaced by $1/(k+1)$ at iteration $k$.

\begin{algorithm}[H]
	\caption{Fictitious play iterations for MFG}
	\label{algo:fictitious-play}
	\begin{algorithmic}[1]
		\Require initial average mean field \( \overline\nu^{(0)} = (\overline\nu^{(0)}_t)_{t \in \timesteps} \), number of iterations $K$
		\Ensure approximate mean field Nash equilibrium
		\For{$k=0,\dots,K-1$}
		\State Update the policy $\pi^{(k+1)} \in \BR(\overline\nu^{(k)})$
		\State Update the mean field \( \nu^{(k+1)} = \EMF(\pi^{(k+1)})\)
		\State Update the average distribution \(\overline\nu^{(k+1)} = \frac{k}{k+1}\overline\nu^{(k)} + \frac{1}{k+1}\nu^{(k+1)}\)
		\EndFor
		\State When possible, recover an averaged policy $\overline\pi^{(K)}$ from $\overline\nu^{(K)}$ by disintegration
		\State \Return \( \overline\pi^{(K)} = (\overline\pi^{(K)}_t)_{t \in \timesteps}\) and \( \overline\nu^{(K)} = (\overline\nu^{(K)}_t)_{t \in \timesteps}\)
	\end{algorithmic}
\end{algorithm}

	The last step should be understood as a reconstruction step. Let $\overline\mu^{(K)}_t(s)=\sum_{a\in\actions}\overline\nu^{(K)}_t(s,a)$. For every state with $\overline\mu^{(K)}_t(s)\neq 0$, define
	\begin{equation}
	\label{eq:fic-play-update-policy}
		\overline\pi^{(K)}_t(a|s) = \frac{\overline\nu^{(K)}_t(s,a)}{\overline\mu^{(K)}_t(s)}, \qquad a \in \actions.
	\end{equation}
	For states with $\overline\mu^{(K)}_t(s)=0$, the action distribution can be chosen arbitrarily. In general, since $\overline\nu^{(K)}$ is an average of flows generated by different policies, the policy obtained by this disintegration need not regenerate exactly $\overline\nu^{(K)}$ when it is run through the dynamics. Thus the averaged mean field is the primary output of fictitious play; the associated policy is mainly useful once the iterates are close to consistency.

The natural convergence object for fictitious play is the averaged mean field $\overline\nu^{(K)}$, together with an associated averaged policy $\overline\pi^{(K)}$ when such a policy can be recovered as above, not the last best-response iterate $(\pi^{(K)},\nu^{(K)})$. Under monotonicity-type conditions, the averaging of best responses can be interpreted as progressively reducing exploitability. In potential MFGs, another useful viewpoint is that fictitious play behaves like a Frank-Wolfe-type method applied to an associated potential. These mechanisms explain why fictitious play is stable in some important classes of MFGs, but they do not give a general convergence guarantee for arbitrary models. We refer to the bibliographic notes at the end of the chapter for comments on the convergence of this algorithm.

As already discussed in the fixed point iterations presented above, the best response computation can be replaced by an approximately greedy policy such as a softmax policy, i.e.,
\[
	\pi^{(k+1)}_t(a|s)  = \softmax_\tau\left(Q_t^{*,\overline{\nu}^{(k)}}(s,\cdot)\right)(a),
\]
for every $(t,s,a) \in \timesteps \times \states \times \actions$, where $\tau$ is a given temperature hyperparameter. This variant has been referred to as the \boldindex{softmax fictitious play algorithm}.

\subsection{Policy iteration}
\label{sec:algo-policy-iteration}

The \boldindex{policy iteration algorithm} is a classical method to solve MDPs. We recall the main ideas before presenting its adaptation to solve MFGs.
It starts with an initial policy $\pi^{(0)}$ and then repeatedly performs the following two steps, for $k \ge 0$:
\begin{itemize}
	\item {\bf Policy evaluation:}\index{policy evaluation} the current policy $\pi^{(k)}$ is evaluated. This involves calculating $V^{\pi^{(k)}}$ (resp. $Q^{\pi^{(k)}}$) defined in~\eqref{eq:V-pi-def} (resp.~\eqref{eq:Q-pi-def}), i.e., the expected long-term reward for each state (resp. each state-action pair) when the system follows the current policy.
	\item {\bf Policy improvement:}\index{policy improvement} Based on the evaluation from the previous step, the policy $\pi^{(k)}$ is  replaced by a better policy, $\pi^{(k+1)}$, i.e., such that $V^{\pi^{(k+1)}}(s) \ge V^{\pi^{(k)}}(s)$ for every state $s \in \states$. Concretely, the policy improvement lemma says that it is sufficient to take $\pi^{(k+1)}$ as a greedy policy with respect to $Q^{\pi^{(k)}}$.
\end{itemize}
The process of policy evaluation and policy improvement is repeated iteratively. Each iteration involves first evaluating the current policy and then improving it. This loop continues until the value function no longer changes, indicating that an optimal policy has been found.

For MFGs, we cannot apply directly the policy iteration method to find a Nash equilibrium. Indeed, the mean field needs to be updated, and not only the policy. Various ways to intertwine these updates can be conceived. Here, we focus on the most common approach.

To facilitate the presentation, we introduce the following notation.
Given a mean field sequence $\nu$ and a policy $\pi$, the associated state-action value function $Q^{\nu,\pi}$ is defined as in~\eqref{eq:Q-pi-def} but with the transitions and the rewards depending on $\nu$. The Bellman equation~\eqref{eq:Q-pi-Bellman} becomes:
\begin{equation}
	\label{eq:Q-nu-pi-Bellman}
	\begin{cases}
		Q^{\nu,\pi}_{T}(s,a) =  r_T(s,a,\nu_T), \quad s \in \states, a \in \actions,
		\\
		Q^{\nu,\pi}_{t}(s,a) = \EE_{\substack{s' \sim P_t(\cdot|s,a,\nu_t), \\ a' \sim \pi_{t+1}(\cdot|s')}}[ r_t(s,a,\nu_t) + Q^{\nu,\pi}_{t+1}(s',a')], \,
		\\
		\qquad\qquad s \in \states, a \in \actions, t = T-1, \dots, 0.
	\end{cases}
\end{equation}
The policy improvement consists in taking one of the greedy policies of the Q-function. We recall the notation~\eqref{eq:greedy-policies} for the set of greedy policies.

The policy iteration algorithm is given in Algorithm~\ref{algo:policy-iteration}. In contrast with the fixed point algorithm and the fictitious play algorithm, it does not involve the computation of a best response.
\begin{algorithm}[H]
	\caption{Policy iteration for MFG}
	\label{algo:policy-iteration}
	\begin{algorithmic}[1]
		\Require initial policy \( \pi^{(0)} = (\pi^{(0)}_t)_{t \in \timesteps} \), number of iterations $K$
		\Ensure approximate mean field Nash equilibrium
		\State Update the mean field \( \nu^{(0)} = \EMF(\pi^{(0)})\)
		\For{$k=0,\dots,K-1$}
		\State Update the value function \( Q^{(k)} = Q^{\nu^{(k)}, \pi^{(k)}}\)
		\State Update the policy $\pi^{(k+1)} \in \mathrm{Greedy}(Q^{(k)})$
		\State Update the mean field \( \nu^{(k+1)} = \EMF(\pi^{(k+1)})\)
		\EndFor
		\State \Return \( \pi^{(K)} = (\pi^{(K)}_t)_{t \in \timesteps}\) and \( \nu^{(K)} = (\nu^{(K)}_t)_{t \in \timesteps}\)
	\end{algorithmic}
\end{algorithm}
With this indexing, each completed iterate satisfies the consistency relation $\nu^{(k)}=\EMF(\pi^{(k)})$. The policy update from $k$ to $k+1$ is computed using the value function evaluated at the previous consistent pair $(\pi^{(k)},\nu^{(k)})$.

Here again, the greedy policy can be replaced by an approximately greedy policy such as a softmax policy, i.e.,
\[
	\pi^{(k+1)}_t(a|s)  = \softmax_\tau\left(Q_t^{(k)}(s,\cdot)\right)(a)
\]
for every $(t,s,a) \in \timesteps \times \states \times \actions$, where $\tau$ is a given temperature hyperparameter. We refer to this softmax variant as the \boldindex{Boltzmann policy iteration algorithm} (\boldindex{BPI}).

Unlike standard policy iteration for a fixed MDP, this procedure changes the MDP itself at each step through the mean-field update. Hence the usual policy-improvement argument does not by itself imply convergence to a mean field Nash equilibrium. Convergence results for policy-iteration-type schemes require additional structure, such as monotonicity, smoothing, or contractivity assumptions. We refer to the bibliographic notes for more details.

\subsection{Online mirror descent}
\label{sec:algo-omd}

In the context of discrete MFGs, \boldindex{online mirror descent} (\boldindex{OMD}) can be viewed as a modification of policy iteration (see Section~\ref{sec:algo-policy-iteration}), in which the greedy policy with respect to the latest Q-function is replaced by a softmax policy with respect to the cumulative Q-function over past iterations.

To this end, we recall the notation $Q^{\nu,\pi}$ for policy evaluation introduced in~\eqref{eq:Q-nu-pi-Bellman} and we introduce the following notation for the regularized greedy improvement: Given a $Q$-function, a regularized greedy policy is computed by taking:
\[
	\pi_t(a|s) = \softmax_\tau\left(Q_t(s,\cdot)\right)(a), \qquad (t,s,a) \in \timesteps \times \states \times \actions,
\]
where $\tau>0$ is the temperature parameter. We use the shorthand $\pi = \softmax_\tau(Q)$.

At iteration $k$ of OMD, the current policy is denoted by $\pi^{(k)}$, and we use the same indexing convention as in policy iteration: $\nu^{(k)}=\EMF(\pi^{(k)})$. The first step consists in computing the Q-function $Q^{(k)} = Q^{\nu^{(k)}, \pi^{(k)}}$ associated with the current consistent pair, see~\eqref{eq:Q-nu-pi-Bellman}. Then, this Q-function is added to the cumulative Q-function with a coefficient $\alpha$:
\[
	\overline{Q}^{(k+1)}_t(s,a) = \overline{Q}^{(k)}_t(s,a) + \alpha Q^{(k)}_t(s,a).
\]
Finally, the new policy is defined as the regularized greedy policy with respect to this cumulative Q-function:
\[
	\pi^{(k+1)}_t(a|s) = \softmax_\tau\left(\overline{Q}^{(k+1)}_t(s,\cdot)\right)(a).
\]
The method is summarized in Algorithm~\ref{algo:omd}.

\begin{algorithm}[H]
	\caption{Online mirror descent for MFG}
	\label{algo:omd}
	\begin{algorithmic}[1]
		\Require initial policy \( \pi^{(0)} = (\pi^{(0)}_t)_{t \in \timesteps} \), initial cumulative value function $\overline Q^{(0)}=0$, coefficients $\alpha > 0$ and $\tau>0$, number of iterations $K$
		\Ensure approximate mean field Nash equilibrium
		\State Update the mean field \( \nu^{(0)} = \EMF(\pi^{(0)})\)
		\For{$k=0,\dots,K-1$}
		\State Update the value function \( Q^{(k)} = Q^{\nu^{(k)}, \pi^{(k)}}\)
		\State Update the cumulative value function \( \overline{Q}^{(k+1)} = \overline{Q}^{(k)} + \alpha Q^{(k)}\)
		\State Update the policy $\pi^{(k+1)} = \softmax_\tau \overline{Q}^{(k+1)}$
		\State Update the mean field \( \nu^{(k+1)} = \EMF(\pi^{(k+1)})\)
		\EndFor
		\State \Return \( \pi^{(K)} = (\pi^{(K)}_t)_{t \in \timesteps}\) and \( \nu^{(K)} = (\nu^{(K)}_t)_{t \in \timesteps}\)
	\end{algorithmic}
\end{algorithm}

The convergence mechanism is different from fictitious play: OMD does not solve a best-response problem at each iteration, but accumulates value-function information and uses the softmax update as a mirror-descent step. Under suitable monotonicity assumptions, Lyapunov arguments establish convergence for related continuous-time learning dynamics and help explain the behavior of OMD. Convergence of the discrete iteration above requires a separate analysis, including conditions on the step sizes.

\section{Optimization-based algorithms for MFGs}
\label{sec:algo-optimization}

The optimization formulation~\eqref{mfvmo-full} makes it possible to apply optimization algorithms to compute Nash equilibria. As an illustration, we present the projected gradient descent algorithm and state a local convergence guarantee for finite-state MFGs under the assumptions below.

	We write
	\[
		\theta=(y,z,\nu)\in
		\mathbb{R}^{S(T+1)}\times\mathbb{R}^{SA(T+1)}
		\times \Delta(\states\times\actions)^\timesteps.
	\]
	We denote by $\Theta$ the feasible set of \eqref{mfvmo-full}, namely the set of all triples $(y,z,\nu)$ satisfying the constraints in \eqref{mfvmo-full}. Finally, we denote the optimal solution set by
	\[
		\Theta^\star=\arg\min_{\theta\in\Theta} f^{\text{\rm MF-OMO}}(\theta).
	\]

\subsection{Projected gradient descent for MF-OMO}
The algorithm of \boldindex{projected gradient descent} (\boldindex{PGD}) updates $\theta$ with an iterative process described in Algorithm \ref{alg:pgd-mfomo}.
\begin{algorithm}[ht]
	\caption{Projected gradient descent for MF-OMO}
	\label{alg:pgd-mfomo}
	\begin{algorithmic}[1]
		\State {\bfseries Input:} feasible initialization $\theta_0=(y^0,z^0,\nu^0)\in\Theta$ and step-sizes $\eta_k>0$ ($k\geq 0$).
		\For{$k=0, 1, \dots$}
		\State Evaluate the gradient $\nabla_{\theta}f^{\text{\rm MF-OMO}}(\theta_k)$ and update $\theta_{k+1}$ as follows:
		\begin{equation}\label{pgd_mfvmo}
			\theta_{k+1}=\textbf{Proj}_{\Theta}(\theta_k-\eta_k \nabla_{\theta}f^{\text{\rm MF-OMO}}(\theta_k)).
		\end{equation}
		\EndFor
	\end{algorithmic}
\end{algorithm}

Here $\eta_k>0$ is the step-size of the $k$-th iteration, for which appropriate ranges are specified below. The projection operator $\textbf{Proj}_{\Theta}$ is the Euclidean projection onto the product constraint set: it projects $\nu$ onto $\Delta(\states\times\actions)^\timesteps$, $z$ onto the polyhedron $\{z\,|\,z\geq0, \|z\|_1\leq SA(T+1)(T+2)r_{\max}\}$, and $y$ onto the $\ell_2$-norm ball with radius $S(T+1)(T+2)r_{\max}/2$.
All these projections can be efficiently evaluated in linear time and (almost) closed-form. 

\subsection{Convergence to Nash equilibria}
To find a Nash equilibrium of the MFG defined by~\eqref{setup:MFG-MDP-nu} and \eqref{consistency_pop_flow} using the PGD algorithm and to guarantee its convergence, we will need a couple of assumptions.

\begin{Assumption}\label{C2smooth}
	For every $t=0,\dots,T-1$, $s,s'\in\mathcal{S}$, and $a\in\mathcal{A}$, the map $\nu_t\mapsto P_t(s'|s,a,\nu_t)$ is second-order continuously differentiable on some open set containing $\Delta(\mathcal{S}\times\mathcal{A})$. For every $t\in\timesteps$, $s\in\mathcal{S}$, and $a\in\mathcal{A}$, the map $\nu_t\mapsto r_t(s,a,\nu_t)$ has the same property.
\end{Assumption}
Here $\nu_t$ is viewed as a vector in $\mathbb{R}^{SA}$ indexed by state-action pairs. The phrase ``second-order continuously differentiable in $\nu_t$'' means that the map admits a $C^2$ extension to an open neighborhood of $\Delta(\mathcal S\times\mathcal A)$ in $\mathbb{R}^{SA}$, and the derivatives are the usual Euclidean derivatives of this extension. The topology is therefore the one induced by the Euclidean norm; since the space is finite-dimensional, any equivalent norm would lead to the same differentiability notion.
This assumption guarantees the existence of optimal solution(s) for \eqref{mfvmo-full}.

We will also need a definability assumption.

\begin{Assumption}
	\label{definability}
	For every $t=0,\dots,T-1$, $s,s'\in\mathcal{S}$, and $a\in\mathcal{A}$, the map $\nu_t\mapsto P_t(s'|s,a,\nu_t)$ is the restriction of a definable function on the log-exp structure\footnote{For simplicity, below we say that a function is definable if it is definable on the log-exp structure.} to $\Delta(\mathcal{S}\times\mathcal{A})$. For every $t\in\timesteps$, $s\in\mathcal{S}$, and $a\in\mathcal{A}$, the map $\nu_t\mapsto r_t(s,a,\nu_t)$ has the same property.
\end{Assumption}
This is one of the most commonly adopted assumptions in the convergence analysis for nonconvex optimization.
Definable functions cover a broad class of smooth and nonsmooth functions, including semialgebraic functions, analytic functions on definable compact sets, and functions built using the exponential and logarithm. In practice, many functions used in numerical implementations belong to this class.

	\begin{example}[Smooth definable finite-state models]
		The reward part of Assumptions~\ref{C2smooth} and~\ref{definability} is satisfied by affine and polynomial mean-field dependencies, since they extend to smooth semialgebraic functions on $\mathbb{R}^{SA}$. For instance, a reward of the form
		\begin{align*}
			r_t(s,a,m)
			 & =
			\rho^0_{t,s,a}
			+
			\sum_{(\bar s,\bar a)\in\states\times\actions}
			\rho^1_{t,s,a,\bar s,\bar a}m(\bar s,\bar a)
			\\
			 & \qquad +
			\sum_{\substack{(\bar s,\bar a)\in\states\times\actions \\
			(\tilde s,\tilde a)\in\states\times\actions}}
			\rho^2_{t,s,a,\bar s,\bar a,\tilde s,\tilde a}
			m(\bar s,\bar a)m(\tilde s,\tilde a)
		\end{align*}
		is smooth and semialgebraic as a function of $m\in\Delta(\states\times\actions)$.

		The transition kernel can also depend smoothly on the mean field. For example, let
		\[
			\ell_{t,s,a,s'}(m)
			=
			\beta^0_{t,s,a,s'}
			+
			\sum_{(\bar s,\bar a)\in\states\times\actions}
			\beta^1_{t,s,a,s',\bar s,\bar a}m(\bar s,\bar a)
		\]
		and define
		\[
			P_t(s'|s,a,m)
			=
			\frac{\exp(\ell_{t,s,a,s'}(m))}
			{\sum_{y\in\states}\exp(\ell_{t,s,a,y}(m))}.
		\]
		This softmax parametrization produces a probability distribution on $\states$, admits a $C^\infty$ extension to an open neighborhood of the simplex, and is definable on the log-exp structure.
	\end{example}

Under Assumptions~\ref{C2smooth} and~\ref{definability}, Theorem~\ref{thm:convergence-PGD} gives a local convergence guarantee for PGD with the specified initialization and step size. If the initial mean-field flow is sufficiently close to an equilibrium flow, the iterates approach the MF-OMO solution set and the exploitability of reconstructed policies tends to zero. If the equilibrium flow is isolated, the neighborhood can be chosen so that the mean-field iterates converge to that flow. Thus suitable initializations can, in principle, recover each isolated equilibrium mean-field flow.

Under Assumption~\ref{C2smooth}, the function $f^{\text{\rm MF-OMO}}$ is twice continuously differentiable on an open neighborhood of the compact set $\Theta$. Hence there exists $M<\infty$ such that
	\[
		\|\nabla_\theta f^{\text{\rm MF-OMO}}(\theta)-\nabla_\theta f^{\text{\rm MF-OMO}}(\theta')\|_2
		\leq M\|\theta-\theta'\|_2,
		\qquad \theta,\theta'\in\Theta.
	\]
	In the theorem below, $M$ denotes such a constant.

\begin{theorem}
	\label{thm:convergence-PGD}
	Under Assumptions \ref{C2smooth} and \ref{definability}, let $\theta_k=(y^k,z^k,\nu^k)$ {\rm(}$k\geq 0${\rm)} be the sequence generated by PGD \eqref{pgd_mfvmo} with $\eta_k=\eta\in(0,2/M)$. Then for any Nash equilibrium $(\pi^\star,\nu^\star)$, there exists $\epsilon_0>0$, such that for any $\nu^0\in\Delta(\mathcal{S}\times\mathcal{A})^\timesteps$ with $\|\nu^0-\nu^\star\|_1\leq \epsilon_0$, if we initialize PGD with $\theta_0=(y^0,z^0,\nu^0)$, where $y^0=[V_1^\star(\nu^0), V_2^\star(\nu^0), \dots,V_T^\star(\nu^0), -V_0^\star(\nu^0)]\in\mathbb{R}^{S(T+1)}$, $z^0=[z_{0,0},\dots,z_{0,T}]\in\mathbb{R}^{SA(T+1)}$, $z_{0,t}=[z_{0,t}^{1},\dots,z_{0,t}^{A}]$ {\rm(}$t\in\mathcal{T}${\rm)}, with
	\[
		\begin{gathered}
  z_{0,t}^{a}=V_t^\star(\nu^0)-r_t(\cdot,a,\nu_t^0)-P_t^a(\nu_t^0)V_{t+1}^\star(\nu^0)\in\mathbb{R}^S,\\
  t=0,\dots,T-1,\qquad a\in\mathcal{A}
\end{gathered}
	\]
	and
	$z_{0,T}^a=V_T^\star(\nu^0)-r_T(\cdot,a,\nu_T^0)\in\mathbb{R}^S$ {\rm(}$a\in\mathcal{A}${\rm)}, then
	we have the following:
	\begin{itemize}
		\item $\lim_{k\rightarrow\infty}{\bf dist}(\theta_k,\Theta^\star)=0$, and $\|\nu^k-\nu^\star\|_2\leq g(\epsilon_0)$ for any $k\geq 0$. Here 
		      $g:\mathbb{R}_+\rightarrow\mathbb{R}_+$ is a function such that $\lim_{\epsilon\rightarrow0}g(\epsilon)=0$.
		      Here ${\bf dist}$ is the point-set distance defined as ${\bf dist}(\theta_k,\Theta^\star)=\inf_{\theta\in\Theta^\star}\|\theta_k-\theta\|_2$.
		\item The limit point set of $(\theta_k)_{k\geq 0}$ is nonempty, and for any limit point $\bar{\theta}=(\bar{y},\bar{z},\bar{\nu})$ of the iteration sequence $(\theta_k)_{k\geq 0}$, any $\bar{\pi}\in\Poldist(\bar{\nu})$ is a Nash equilibrium of  the MFG defined by~\eqref{setup:MFG-MDP-nu} and~\eqref{consistency_pop_flow}.
		\item For any sequence of $\pi^k\in\Poldist(\nu^k)$,  $\lim_{k\rightarrow\infty}\Expl(\pi^k)=0$.
	\end{itemize}

	In addition, for any isolated Nash equilibrium $(\pi^\star,\nu^\star)$ (\ie, there exists $\epsilon>0$ such that any Nash equilibrium $(\pi,\nu)$ with $\|\nu-\nu^\star\|_2\leq \epsilon$ satisfies $\nu=\nu^\star$), if we initialize PGD with $\theta_0=(y^0,z^0,\nu^0)$ in the same way as above, then $\lim_{k\rightarrow\infty}\nu^k=\nu^\star$ and $\lim_{k\rightarrow\infty}\Expl(\pi^k)=0$ for any sequence $\pi^k\in\Poldist(\nu^k)$.
\end{theorem}

\section{Bibliographic notes}
\label{sec:algo-bibliographic-notes}

	The use of fixed point iterations to solve MFGs goes back to early works such as~\citep{gueant2011mean}. Pure fixed point iterations can fail to converge, and in the continuous time and space setting we refer, for instance, to Section~2 of~\citep{lauriere2021numerical} for illustrations in linear-quadratic MFGs. Under suitable local assumptions, however, pure fixed point iterations can converge; see, for example, the local convergence result in~\citep[Proposition 1]{li2021simple}. The damped and softmax variants discussed above are stabilization mechanisms for regimes where the undamped update is too sensitive.

Fictitious play has been introduced by~\cite{brown1951iterative} and further analyzed in~\citep{robinson1951iterative}. The adaptation to MFGs has been introduced by~\cite{cardaliaguet2017learning} in the continuous time and space setting. In the discrete time and space setting, it has been analyzed by~\cite{elie2020convergence} and its convergence in a continuous time version has been proved in~\citep{perrin2020fictitious} under the Lasry-Lions monotonicity condition. In the context of potential MFGs, fictitious play can be viewed as the classical Frank-Wolfe algorithm~\citep{frank1956algorithm}, as explained in~\citep{geist2022concave}; related conditional-gradient analyses for potential MFGs, including rates for exploitability, are developed in~\citep{lavigne2023generalized}.

Policy iteration for MDPs dates back to the work of~\cite{bellman1957DynamicProgramming} and~\cite{howard1960dynamic}. In the context of MFGs, it was introduced and analyzed in~\citep{cacace2021policy} in the continuous time and space setting. Convergence for more general cost structures was established in~\cite{lauriere2023policy}. 
Several variants have been considered. For example, smoothed policy-iteration algorithms for potential MFGs are studied in~\citep{tang2024learning}; a related analysis of Newton's method for MFGs with non-separable Hamiltonians is given in~\citep{camilli2024quadratic}. When the action space is finite, the greedy policy-improvement step can also be replaced by a softmax of the current policy's Q-function, as in the BPI variant presented in Section~\ref{sec:algo-policy-iteration}. The Boltzmann iteration studied in~\cite{cui2021approximately} instead uses the optimal Q-function of the frozen-flow MDP and a prior-weighted softmax update. Entropy-regularized MFG models and policy-gradient learning are studied, e.g., in~\citep{guo2022entropy}.

The mirror descent algorithm is a classical optimization method introduced by~\cite{nemirovskii1983problem}. Online mirror descent lets the objective or feedback change over iterations; see~\cite[Section 2.6]{shalev2012online} and~\cite[Section 5.3]{hazan2016introduction}. In MFGs, online mirror descent was introduced in the continuous time and space setting by~\cite{hadikhanloo2017learning}. For finite-state, discrete-time MFGs, the version presented above was introduced and studied in~\citep{perolat2022scaling}, with convergence under monotonicity assumptions. Implementations of fictitious play and online mirror descent algorithms with neural networks and reinforcement learning have been proposed in~\citep{lauriere22scalable-pmlr}. For non-stationary, discrete-time MFGs with continuous state and action spaces, a fictitious-play method combining deep reinforcement learning and population-density estimation is developed in~\citep{magnino2025continuous}. 
Related independent and decentralized learning formulations have also been studied directly in finite-population or decentralized environments, without relying on a population generative model~\citep{yardim2023policy,yongacoglu2022independent,subramanian2022decentralized}.

	Model-free and sample-based variants of MFG learning algorithms have been studied in several finite-state or discrete-time settings. Stationary, learning-while-playing, and single-sample-path approaches are analyzed in~\citep{subramanian2019reinforcementMFG,xie2021learning,zaman2023oracle}. Value-based and Q-learning methods, including regularized, discounted, and average-cost formulations, are studied in~\citep{guo2023general,angiuli2020unifiedRLMFGMFC,anahtarci2023q,anahtarci2021learningfittedaverage,anahtarci2023learning}. Recent work also studies model-based and statistical-efficiency questions for MFG and mean-field reinforcement learning (RL) under additional structural or function-approximation assumptions~\citep{huang2023statistical,huang2024model}. For ergodic finite-state MFGs, trust-region policy optimization and finite-sample convergence guarantees are studied in~\citep{ocello2025finite}. These results typically rely on structural assumptions such as contraction, regularization, monotonicity, realizability, or access to suitable samples. This reliance on structure is consistent with the broader algorithmic-game-theory picture: equilibrium computation in finite games is PPAD-complete in general~\citep{daskalakis2009complexity,chen2009settling}, with related hardness results even for restricted bimatrix games~\citep{mehta2014constant}. For MFGs, limitations of broad Lipschitz classes, both for approximation and computational tractability, are discussed in~\citep{yardim2024tractable}.

	Related reinforcement learning methods have been developed for mean field control (MFC) and mean-field MDP formulations, where the objective is social optimality rather than Nash equilibrium. These problems lead to dynamic programming on the population distribution and to different consistency conditions; see, for instance,~\citep{carmona2023model,gu2023dynamic,gu2021mean,wang2021globallqmfcg,pasztor2023efficientmodelbased,angiuli2022reinforcementmixedmfcg,angiuli2023convergence}. We do not develop this direction here because the focus of these notes is on non-cooperative MFGs.

	Another class of numerical methods directly minimizes exploitability, without explicitly using best-response or policy-evaluation updates as in the algorithms presented above. In the discrete-time, discrete-space, stationary setting, Bench-MFG~\citep{magnino2026benchmfg} introduces mean-field particle swarm optimization (MF-PSO), which is a black-box adaptation of particle swarm optimization to exploitability minimization, and compares it with fixed point, fictitious play, policy iteration, and OMD across several benchmark classes. This approach is useful as an empirical baseline because it does not rely on classical MFG assumptions such as contraction, monotonicity, or potential structure.

	Policy-space response-oracle (PSRO) methods form a related but conceptually different family. Mean-field PSRO~\citep{muller2022learningpsro}, inspired by PSRO in multi-agent reinforcement learning, maintains a finite set of policies, computes an equilibrium of the restricted MFG over this set, and adds best responses to the resulting population mixture. This policy-set viewpoint is useful when one wants to build a small library of strategically distinct policies, and it also leads naturally to mean-field correlated and coarse correlated equilibria~\citep{campi2022correlated,muller2022learningmfce}. The restricted-game subproblem is generally nonlinear in MFGs, so it requires additional optimization or regret-minimization tools; for this reason, we do not include it among the basic algorithms above.

	Another related direction is to learn from demonstrations rather than from direct interaction with a known MFG model. Inverse reinforcement learning and imitation learning for MFGs have been studied in~\citep{yang2018learningdeepmfg,chen2021agentinverseRL,chen2022individualinverseRL,ramponi2024imitation,anahtarci2025maximum,lambrecht2026populationaware}. These works address a different problem from equilibrium computation: the goal is to infer or imitate behavior generated by an equilibrium population, possibly under an unknown reward.

Assumption~\ref{C2smooth} has been adopted in  a number of  existing works, including \citep{guo2019learning}. The closed-form definition of the function $g$ appearing in the statement of Theorem~\ref{thm:convergence-PGD} is given in \cite{guo2024mf}.

General background for PGD algorithms can be found in~\cite{duchi2008efficient,liu2009efficient,songsiri2011projection,condat2016fast}
for projection onto the probability simplex and specialized polyhedra, and in~\cite{parikh2014proximal} for projection onto $\ell_2$-norm balls.
For PGD applied to MF-OMO, \cite{guo2024mf} prove local convergence results showing how different initializations can recover different isolated Nash equilibria.
For MFGs with linear rewards and dynamics independent of the mean field, enumerating complementarity patterns reduces equilibrium computation to a finite family of linear programs, as discussed in Section~\ref{sec:maths-MFG}. This class of MFGs includes the left-right (LR) and rock-paper-scissors (RPS) problems in \cite{cui2021approximately}.

\chapter{Applications}
\label{chap:applications}

This chapter presents four applications, with numerical results. 
Code and numerical data for the crowd-motion and traffic-routing examples in this chapter are available in an accompanying repository.\footnote{\url{https://github.com/mlauriere/optimization-learning-discrete-mfg-book-code}.} These first experiments are implemented using the open source library \texttt{OpenSpiel}\index{OpenSpiel}.\footnote{\url{https://openspiel.readthedocs.io/en/latest/intro.html}} The SIS and multiple-equilibrium experiments are implemented in MFGLib.\footnote{\url{https://github.com/radar-research-lab/MFGLib}.}

\section{Crowd motion}\index{crowd motion}
\label{sec:crowd-motion-numerics}

We compare learning algorithms in two variants of the grid-world crowd-motion model introduced in Section~\ref{sec:ex-mf-crowd-motion}: a concert room and an asymmetric maze. Both use logarithmic congestion, movement costs and positional rewards, with the same algorithms and hyperparameters.

\subsection{Concert room}

In the concert-room model, the positional reward attracts agents toward the lower part of the room. The agents do not exit the grid world.

\paragraph{Environment. } We take $d_1=d_2=d=15$ and use coordinates $(x,y)\in\{0,\dots,14\}^2$, with $x$ increasing to the right and $y$ downward. Walls are forbidden states. The initial distribution is uniform on the twenty starting cells shown in the first column of the distribution figures. The marked attraction locations are $g_1=(1,14)$ and $g_2=(13,14)$, with weights $b_1=b_2=10$. They lie on the bottom wall and influence preferences through distance; they are not reachable exits.

We use $T=20$ in the book's notation, with reward terms at times $0,\dots,20$. This corresponds to the implementation's horizon parameter $H=21$ and a zero-payoff terminal state at time $21$. Write $\mu_t(s)$ for the population mass at cell $s$. The reward has the equivalent immediate-movement-cost representation
\begin{equation}
\label{eq:crowd-numerical-reward}
\begin{aligned}
r_t(s,a,\mu_t)
&=2-\frac{2}{d}\sum_{j=1}^{2}b_j\|s-g_j\|_1\\
&\quad-c_{\mathrm{cong}}\log\bigl(\widetilde\mu_t(s)+\varepsilon_{\mathrm{reg}}\bigr)
-\mathbf{1}_{\{t<T\}}\frac{\|a\|_1}{d},
\end{aligned}
\end{equation}
where $\varepsilon_{\mathrm{reg}}\simeq 10^{-25}$, $\widetilde\mu_0(s)=1/d^2$, and $\widetilde\mu_t=\mu_t$ for $t\geq 1$. The initial uniform density in this formula is a convention for the initial reward, not the actual initial population distribution. Its contribution is independent of the policy. We use $c_{\mathrm{cong}}\in\{1,10\}$. The implementation records a movement cost at the following decision state; charging it immediately for $t<T$ preserves the total payoff. The final outgoing move has no subsequent nonzero reward.

The available displacements are stay and the four cardinal moves, restricted to those that do not enter a wall:
\[
\actions(s)=\bigl\{a\in\{(0,0),(-1,0),(1,0),(0,-1),(0,1)\}:s+a\text{ is not a wall}\bigr\}.
\]
After choosing $a_t\in\actions(s_t)$, an additional displacement $\epsilon_t$ is drawn at the intermediate position $s_t+a_t$. With probability $1-p$ it is zero; with probability $p$ it is uniform on $\actions(s_t+a_t)$, including stay. Thus, for $t<T$,
\[
s_{t+1}=s_t+a_t+\epsilon_t.
\]
We use $p\in\{0,0.5\}$. The boundary walls prevent wrap-around.

\paragraph{Algorithms. } We compare seven algorithm families:
 Fixed Point, Softmax Fixed Point, Damped Fixed Point, Fictitious Play, Softmax Fictitious Play, Boltzmann Policy Iteration, and Online Mirror Descent. Fixed Point iteratively computes a best response to the current mean field and then updates the mean field induced by this policy. Softmax Fixed Point replaces the greedy best response by a softmax of the optimal Q-function, while Damped Fixed Point averages the current mean field with the flow induced by a best response; see Section~\ref{sec:algo-fixed-point}. Fictitious Play uses a uniform average over past mean fields, and Softmax Fictitious Play replaces its best response by a softmax policy with respect to the optimal Q-function; see Section~\ref{sec:algo-fictitious-play}. Boltzmann Policy Iteration is the softmax version of policy iteration, and Online Mirror Descent evaluates the current policy, accumulates Q-functions, and then updates the policy by applying a softmax to the cumulative Q-function; see Sections~\ref{sec:algo-policy-iteration} and~\ref{sec:algo-omd}.
 
 All methods start from the policy that is uniform over legal actions at each state, recorded as update $0$. We perform $1000$ learning updates. Fictitious play includes the initial policy's occupation measure and assigns weight $1/(k+1)$ to the new response at update $k\geq 1$. The policy merges used for fictitious play and damping are weighted by state occupancies, rather than averaging action probabilities pointwise; see~\eqref{eq:fic-play-update-policy}. Since the transitions in this example are independent of the mean field, these merges induce the corresponding averages of population flows. Table~\ref{tab:crowd-learning-parameters} gives the parameters, also used in the maze example.

\begin{table}[t]
\centering
\caption{Parameters for the concert-room and maze comparisons.}
\label{tab:crowd-learning-parameters}
\small
\begin{tabularx}{\textwidth}{@{}>{\raggedright\arraybackslash}p{0.42\textwidth}Xr@{}}
\toprule
Method & Parameter & Value\\
\midrule
OMD & Step size $\eta$ & $0.1,\ 0.01$\\
Damped fixed point & New-flow weight $\alpha$ & $0.01$\\
Softmax fixed point and softmax fictitious play & Temperature $\tau$ & $0.1$\\
BPI & Inverse temperature $\beta$ & $0.1$\\
\bottomrule
\end{tabularx}
\end{table}

	Throughout the numerical experiments, exploitability refers to the quantity defined in~\eqref{eq:def-exploitability}. The concert-room and maze curves show raw, unnormalized exploitability. We plot every update through $20$ and every fifth update thereafter, using the same sampling for all methods. In Sections~\ref{sec:appli-sis} and~\ref{sec:finding-multiple-nes}, \boldindex{normalized exploitability} means the exploitability at a given update divided by its value at initialization.

Distribution panels show the physical-time evolution induced by the policy returned after the final learning update. For these crowd-motion examples, transition expectations are evaluated tabularly rather than by sampling population trajectories. The colorbar label, Population density, denotes probability mass per grid cell. Each grid uses a common color scale across its selected methods and all saved physical times, without renormalizing individual panels.

\paragraph{Numerical results. } The experiments are implemented using the open source library \texttt{OpenSpiel}\index{OpenSpiel}.

\begin{figure}[!htb]
\centering
\includegraphics[width=\textwidth]{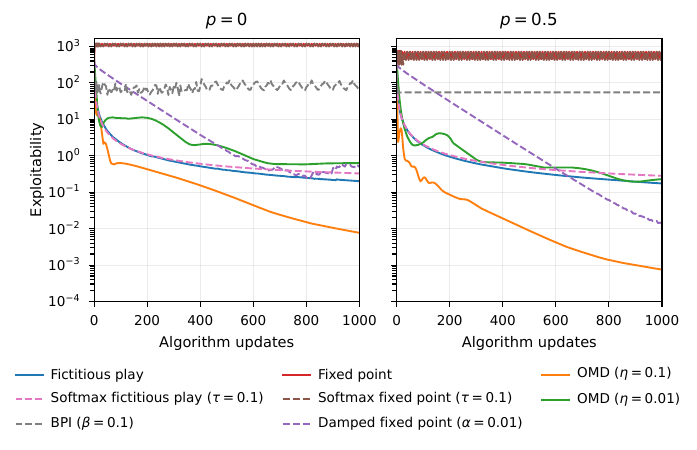}
\caption{Concert-room crowd motion with $c_{\mathrm{cong}}=1$: raw exploitability over $1000$ learning updates, with $p=0$ on the left and $p=0.5$ on the right. Update $0$ is the uniform legal-action policy.}
\label{fig:ex-crowd-evac-cong1-exploitability}
\end{figure}

\begin{figure}[htbp]
\centering
\includegraphics[width=0.9\textwidth]{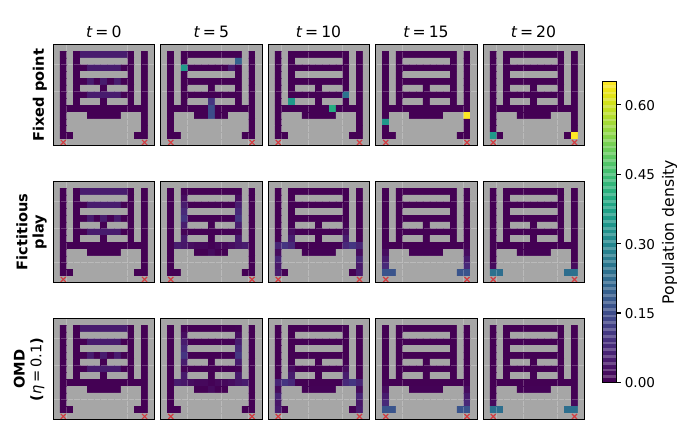}
\caption{Concert-room population flows with $c_{\mathrm{cong}}=1$ and $p=0$, under the policies returned after $1000$ updates. Rows: fixed point, fictitious play, and OMD with $\eta=0.1$. Columns: physical-time indices $t=0,5,10,15,20$. The common colorbar gives probability mass per cell; grey cells are walls and red crosses mark attraction locations.}
\label{fig:ex-crowd-evac-cong1-noise0p0-distribs}
\end{figure}

\begin{figure}[htbp]
\centering
\includegraphics[width=0.9\textwidth]{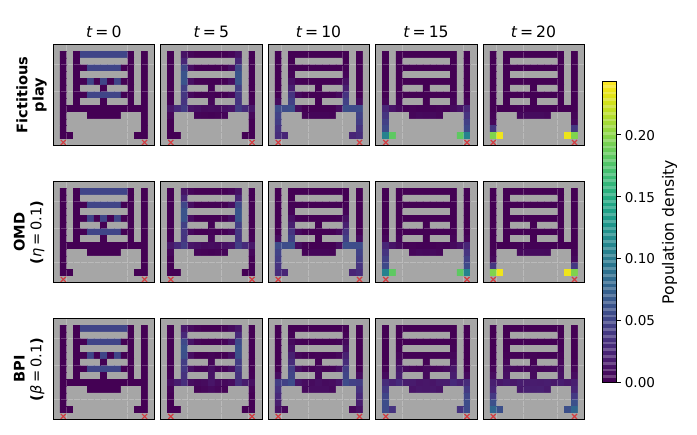}
\caption{Concert-room population flows with $c_{\mathrm{cong}}=1$ and $p=0.5$, under the policies returned after $1000$ updates. Rows: fictitious play, OMD with $\eta=0.1$, and BPI with $\beta=0.1$. Columns: physical-time indices $t=0,5,10,15,20$. The common colorbar gives probability mass per cell; red crosses mark attraction locations.}
\label{fig:ex-crowd-evac-cong1-noise0p5-distribs}
\end{figure}

\paragraph{Results with $c_{\mathrm{cong}}=1$.}
Figure~\ref{fig:ex-crowd-evac-cong1-exploitability} shows exploitability for the eight configurations, with $p=0$ on the left and $p=0.5$ on the right. Figures~\ref{fig:ex-crowd-evac-cong1-noise0p0-distribs} and~\ref{fig:ex-crowd-evac-cong1-noise0p5-distribs} show the induced flows for selected returned policies.
We make the following observations:
\begin{itemize}
\item After $1000$ updates, OMD with $\eta=0.1$ has the lowest exploitability among the tested configurations at both noise levels: approximately $0.0077$ for $p=0$ and $0.00075$ for $p=0.5$. Fictitious play finishes below its softmax variant at both noise levels.
\item Fixed point and softmax fixed point oscillate at high exploitability. With the selected weight $\alpha=0.01$, damping has a long initial transient and becomes oscillatory in the noiseless case, while it continues decreasing in the noisy case. BPI remains highly exploitable.
\item The population panels compare physical flows under the returned policies, including policies with large remaining exploitability. In the noisy case, fictitious play and OMD with $\eta=0.1$ move most of the population toward the lower corridors, while BPI produces a more diffuse flow.
\end{itemize}

\begin{figure}[htbp]
\centering
\includegraphics[width=\textwidth]{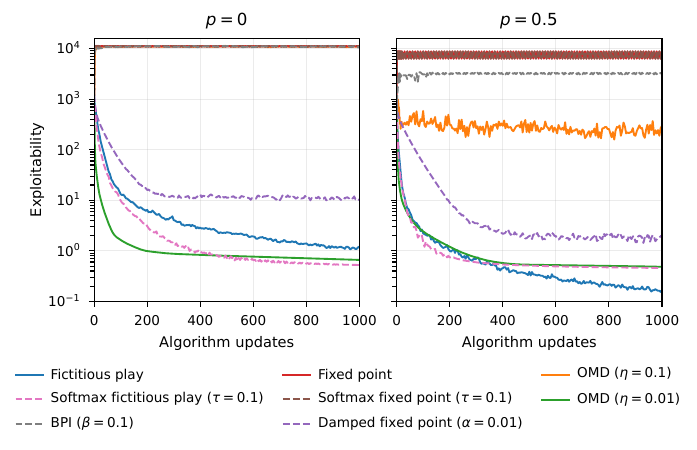}
\caption{Concert-room crowd motion with $c_{\mathrm{cong}}=10$: raw exploitability over $1000$ learning updates, with $p=0$ on the left and $p=0.5$ on the right. Update $0$ is the uniform legal-action policy.}
\label{fig:ex-crowd-evac-cong10-exploitability}
\end{figure}

\begin{figure}[htbp]
\centering
\includegraphics[width=0.9\textwidth]{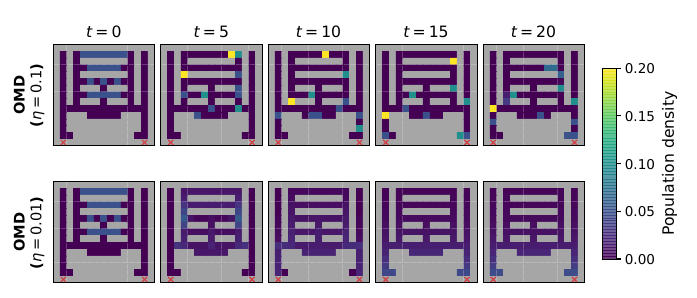}
\caption{Concert-room population flows with $c_{\mathrm{cong}}=10$ and $p=0$, under the policies returned after $1000$ updates. Rows: OMD with $\eta=0.1$ and $\eta=0.01$. Columns: physical-time indices $t=0,5,10,15,20$. The common colorbar gives probability mass per cell.}
\label{fig:ex-crowd-evac-cong10-noise0p0-distribs}
\end{figure}

\paragraph{Results with $c_{\mathrm{cong}} = 10$.}
We now increase the congestion coefficient and keep the same two noise levels.
Figure~\ref{fig:ex-crowd-evac-cong10-exploitability} shows exploitability for the eight configurations, with $p=0$ on the left and $p=0.5$ on the right. Figure~\ref{fig:ex-crowd-evac-cong10-noise0p0-distribs} shows the distribution evolution for a few selected algorithms.
\begin{itemize}
\item With $p=0$, softmax fictitious play has the lowest exploitability after $1000$ updates, approximately $0.53$, followed by OMD with $\eta=0.01$ at approximately $0.66$. OMD with $\eta=0.1$ remains highly exploitable. Figure~\ref{fig:ex-crowd-evac-cong10-noise0p0-distribs} compares the flows returned by the two OMD configurations.
\item With $p=0.5$, fictitious play, softmax fictitious play and OMD with $\eta=0.01$ finish at exploitabilities of approximately $0.17$, $0.46$ and $0.48$, respectively. OMD with $\eta=0.1$ has a much larger final value.
\item These cases illustrate sensitivity to the update rule and step size at a fixed learning budget. Oscillations and concentration in the high-exploitability flows are not caused by sampling population trajectories: the transition expectations are evaluated tabularly.
\end{itemize}

\subsection{Asymmetric maze}
\label{sec:asymmetric-maze}

\begin{figure}[!htb]
\centering
\includegraphics[width=\textwidth]{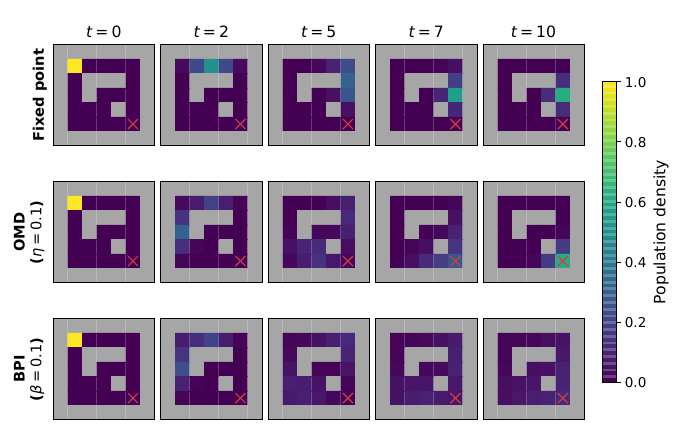}
\caption{Asymmetric maze with $p=0.5$: flows under the policies returned after $1000$ updates. Rows: fixed point, OMD with $\eta=0.1$, and BPI with $\beta=0.1$. Columns: physical-time indices $t=0,2,5,7,10$. The common colorbar gives probability mass per cell, including the initial point mass. Grey cells are walls and the red cross marks the attraction location.}
\label{fig:asymmetric-maze-distributions}
\end{figure}

As a second grid-world example, we change the geometry to the asymmetric $7\times 7$ maze shown in Figure~\ref{fig:asymmetric-maze-distributions}. The initial distribution is concentrated at $(1,1)$ and there is one reachable attraction location, $(5,5)$, with weight $10$. We use the same reward and transition conventions as in the concert-room example, with $d=7$, $c_{\mathrm{cong}}=1$, and a single term in the positional sum in~\eqref{eq:crowd-numerical-reward}. The horizon is $T=10$ in the book's notation, corresponding to $H=11$ in the implementation, and we perform $1000$ updates with the parameters in Table~\ref{tab:crowd-learning-parameters}.

\begin{figure}[htbp]
\centering
\includegraphics[width=\textwidth]{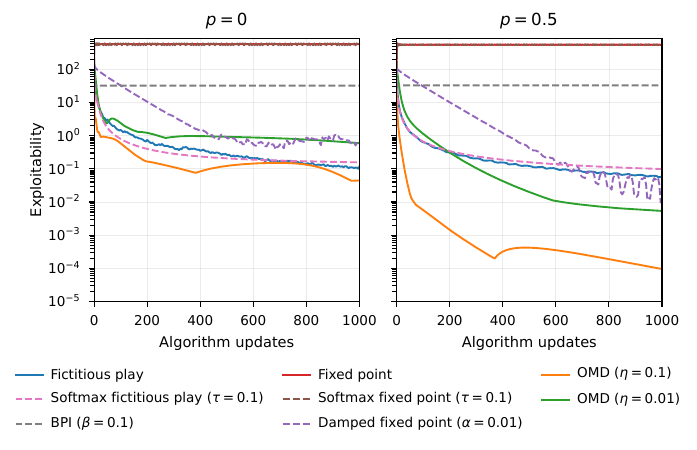}
\caption{Asymmetric maze: raw exploitability over $1000$ learning updates, with $p=0$ on the left and $p=0.5$ on the right.}
\label{fig:asymmetric-maze-exploitability}
\end{figure}

Figure~\ref{fig:asymmetric-maze-exploitability} compares the two noise levels. After $1000$ updates, OMD with $\eta=0.1$ has the lowest exploitability among the tested configurations in both cases, approximately $0.045$ for $p=0$ and $9.5\times 10^{-5}$ for $p=0.5$. The asymmetric corridors provide a different geometry for comparing the update rules, while the common initial point mass makes the subsequent spreading visible. Fixed point remains highly exploitable, and BPI also has substantial remaining exploitability. Figure~\ref{fig:asymmetric-maze-distributions} shows their returned-policy flows together with that of OMD.

\section{Traffic routing}\index{traffic routing}
\label{sec:appli-traffic}

We present an illustrative example of a traffic routing model.

\paragraph{Environment. }
The road network is represented by a directed graph, where $L$ denotes the set of edges or links. In the numerical example, the origin link is $(A,B)$ and the destination link is $\ell_{\mathrm{dest}}=(E,F)$; the letter $D$ in the network denotes an intermediate node. We use $H=17$ routing updates, indexed by $t=0,\dots,16$, and a terminal state at time $17$. Thus $T=H$ in the book's formulation, with a dummy terminal action. One update represents $\Delta=0.25$ units of physical time.

In this model, the agents are cars navigating the road network. The physical state at time $t$ is $s_t=(\ell_t,w_t)$, where $\ell_t\in L$ is the current link and $w_t$ is the remaining number of waiting steps. A wait at least as large as the number of remaining updates is represented by that number, since no further routing decision can then occur before the deadline. This gives a finite physical state space contained in $L\times\{0,\dots,H\}$.

An agent with $w_t=0$ who has not reached the destination chooses a successor link. Waiting agents and agents already at the destination have only a dummy action. The global action set can therefore be written as $\actions=L\cup\{a_\circ\}$, with the legal choices restricted according to the state. The destination is absorbing.

The initial state is $(\ell_0,w_0)=((A,B),0)$. If $w_t>0$, the agent stays on the current link and $w_{t+1}=w_t-1$. If $w_t=0$ and $\ell_t\neq\ell_{\mathrm{dest}}$, the agent chooses a successor link $a_t$ and sets $\ell_{t+1}=a_t$.

Congestion is evaluated after these link choices. Let $\ell^+(s,a)$ denote the resulting link: it is $a$ for an agent choosing a successor and the current link for a waiting or absorbed agent. From the current state-action mean field $\nu_t$, define the post-routing link mass
\[
\rho_t(\ell)=\sum_{s,a}\nu_t(s,a)\mathbf{1}_{\{\ell^+(s,a)=\ell\}}.
\]
For a population scale $N$, the mean-field traffic volume on $\ell$ is $N\rho_t(\ell)$. The corresponding finite-player game uses the realized number of vehicles on the link.

All capacities are fixed at $4$ as $N$ varies. The travel-time function is
\[
\tau_\ell(v)=\tau_\ell^0\left(1+\kappa_\ell\frac{v}{4}\right),
\]
where $\kappa_\ell=1$ on $(B,C)$ and $(D,E)$ and zero on the other links. The values $\tau_\ell^0=\tau_\ell(0)$ are the travel times when the traffic volume on the link is zero:
\[
\begin{array}{c|rrrrrrr}
\ell&(A,B)&(B,C)&(B,D)&(C,D)&(C,E)&(D,E)&(E,F)\\
\hline
\tau_\ell^0&0&1&2&0.25&2&1&0
\end{array}.
\]
On entering link $a_t$ at update $t$, the waiting counter becomes
\[
w_{t+1}
=\min\left\{H-t-1,\,
\max\left\{0,\left\lfloor\frac{\tau_{a_t}(N\rho_t(a_t))}{\Delta}\right\rfloor-1\right\}\right\}.
\]
There is no additional travel noise. The cap retains every vehicle in the population until the common horizon, including vehicles that cannot finish their journey.

Each player seeks an early arrival at $\ell_{\mathrm{dest}}$. Let $\tau_{\mathrm{arr}}$ be the index of the routing action that first enters this link; arrival during the move from $t$ to $t+1$ is recorded as $\tau_{\mathrm{arr}}=t$. The total payoff is
\[
G=
\begin{cases}
-\Delta\tau_{\mathrm{arr}},&\text{if the destination is reached before termination},\\
-H\Delta,&\text{otherwise}.
\end{cases}
\]
In particular, the nonarrival payoff is $-4.25$. An equivalent representation in our finite-horizon framework pays $-t\Delta$ on first entry into the destination at update $t$, pays $-H\Delta$ at the terminal time if the destination has not been reached, and gives zero reward otherwise. This representation uses the physical state above without storing the past arrival time. Arrived vehicles stay at the destination until the horizon. One could extend the model to handle multiple destinations by an additional destination component in the state.

\paragraph{Algorithm. }
For the population illustration, we perform ten OMD updates with step size $1$ and a uniform legal-action initial policy, using a population scale $N=4$. The implementation is based on \texttt{OpenSpiel} with local fixed-horizon and policy-evaluation adapters. Transition probabilities and expectations are evaluated tabularly, in finite precision, without sampling population trajectories.

\begin{figure}[htbp]
\centering
\includegraphics[width=0.8\textwidth]{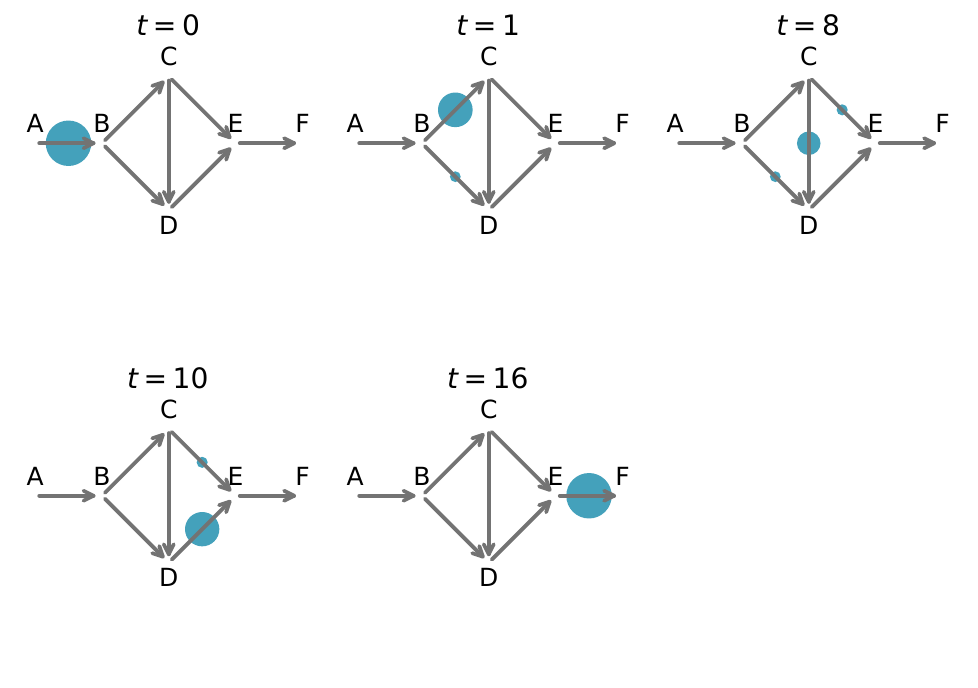}
\caption{Braess-network MFG with population scale $N=4$: flows under the policy returned after ten OMD updates. Panels show indices $t=0,1,8,10,16$ in row-major order, corresponding to physical times $0,0.25,2,2.5,4$. Circle radius is proportional to link probability. The origin and destination links are $(A,B)$ and $(E,F)$.}
\label{fig:braess_dynamics}
\end{figure}

\paragraph{Numerical results. }
We use the network associated with the \boldindex{Braess paradox}~\citep{braess2005paradox}. Figure~\ref{fig:braess_dynamics} shows the population flow induced by the policy returned after ten OMD updates. All mass reaches the destination by the last displayed slice.

We compare the cost of twenty learning updates in the finite-agent Braess game and its MFG approximation. The finite-player solver uses the external-sampling Monte Carlo variant of \boldindex{counterfactual regret minimization} (\boldindex{CFR}), with the library's \texttt{SIMPLE} averaging scheme; the MFG solver uses OMD with step size $1$. The timings compare fixed update budgets, without matching the final equilibrium accuracy.

\begin{figure}[htbp]
\centering
\includegraphics[width=0.6\textwidth]{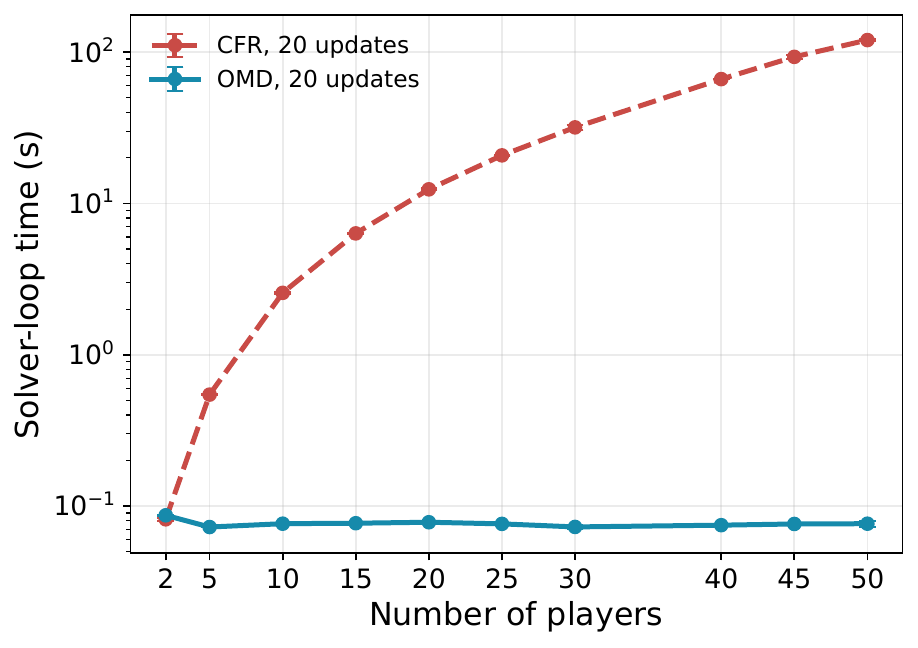}
\caption{Iteration-loop time for twenty updates of external-sampling CFR in the $N$-player Braess game and OMD in the corresponding MFG. Markers are means of three repetitions; error bars show sample standard deviations. Network capacities and horizon are fixed as described in the text. The comparison uses equal update counts.}
\label{fig:computation_time_algo_n_vs_mean_field}
\end{figure}

We use $N\in\{2,5,10,15,20,25,30,40,45,50\}$, keeping the network and capacities fixed. For each $N$, we run three repetitions of twenty updates per solver. Figure~\ref{fig:computation_time_algo_n_vs_mean_field} reports the mean iteration-loop time and sample standard deviation (construction, process startup, warm-up, exploitability evaluation and plotting are excluded).

The measured CFR loop cost increases strongly with $N$, while representative-agent OMD cost remains roughly constant on this fixed network. For example, the mean times at $N=50$ are approximately $120$ seconds and $0.076$ seconds, respectively, on the ARM64 macOS platform that we used for the simulations reported here. Notice that higher $N$ increases demand relative to the fixed capacities, and some journeys can remain unfinished at the deadline.

\section{SIS models}
\label{sec:appli-sis}

In this section, we return to the SIS model\index{SIS model} introduced earlier. In this model, a large population of players chooses at each time step whether to social distance or go out. Susceptible (S) players who go out may become infected (I), with probability proportional to the current mass of infected players, whereas susceptible players who social distance remain healthy. Infected players recover with a fixed probability at each time step. Each player seeks to minimize the costs associated with social distancing and infection.

\begin{figure}[h]
	\centering
	\includegraphics[width=0.9\textwidth]{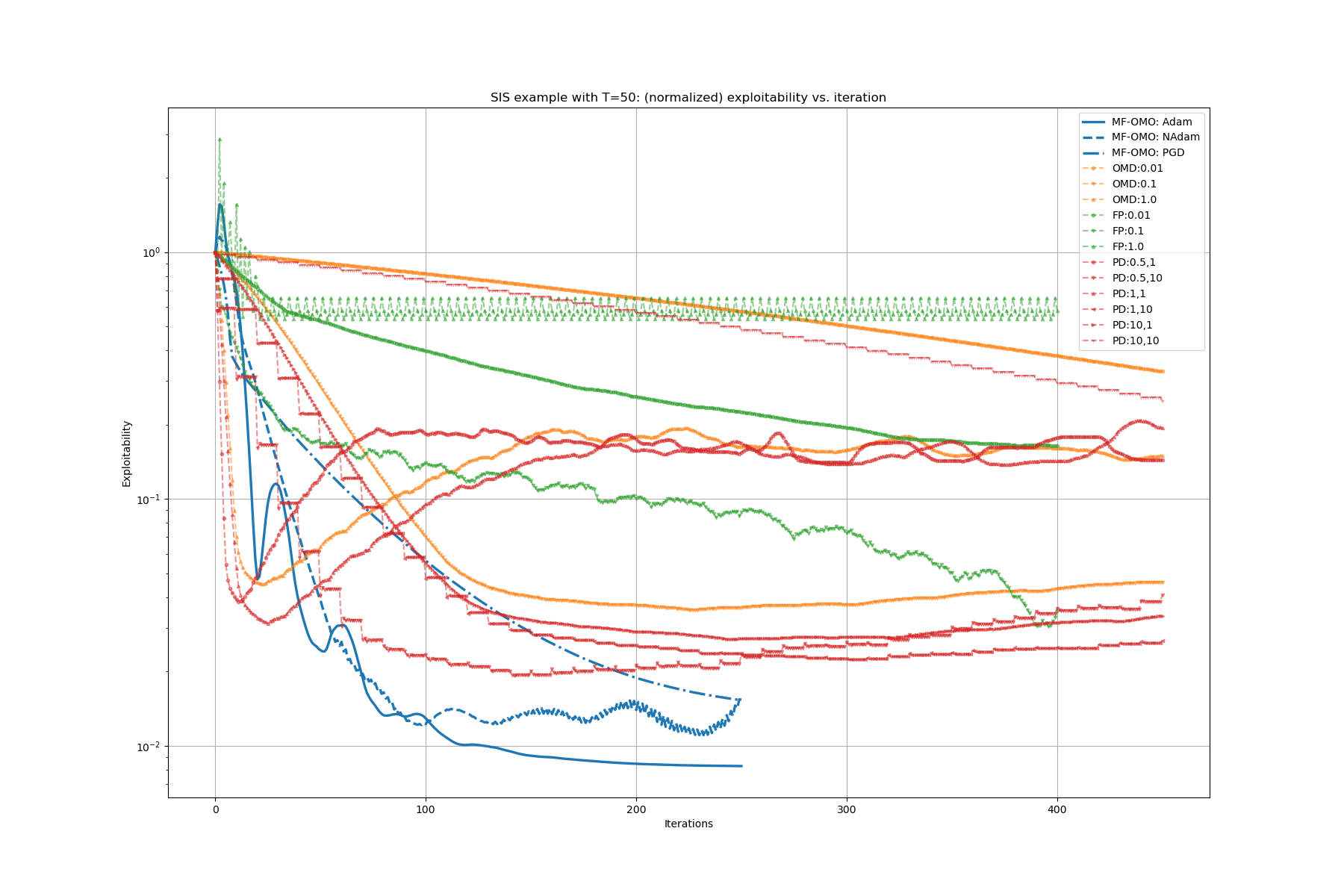}
	\caption{Comparison of different algorithms on SIS model with $T=50$: convergence against iterations}
	\label{fig:sis_iter_50}
\end{figure}

\begin{figure}[h]
	\centering
	\includegraphics[width=0.9\textwidth]{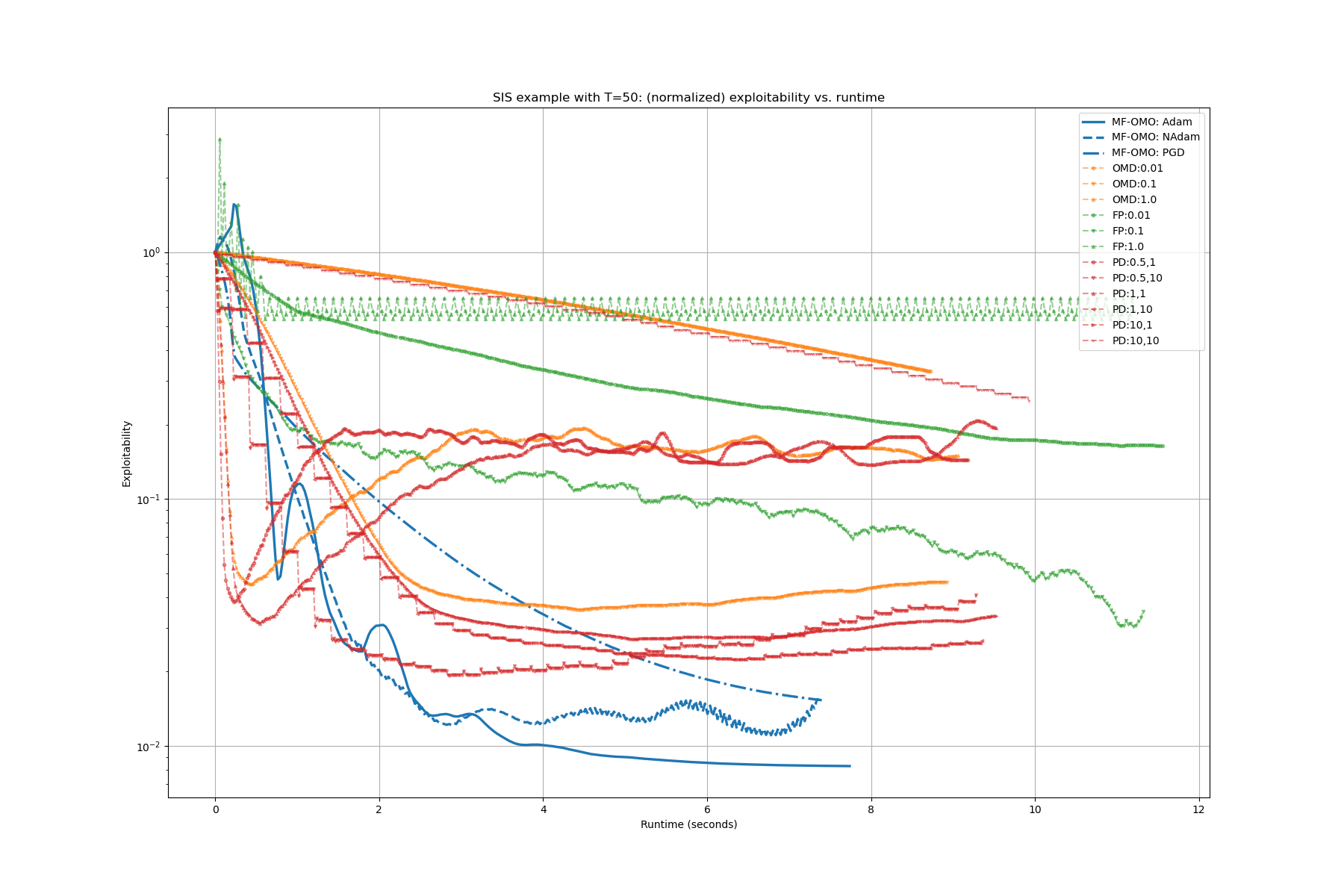}
	\caption{Comparison of different algorithms on SIS model with $T=50$: Convergence against runtime}
	\label{fig:sis_time_50}
\end{figure}

\begin{figure}[h]
	\centering
	\includegraphics[width=0.9\textwidth]{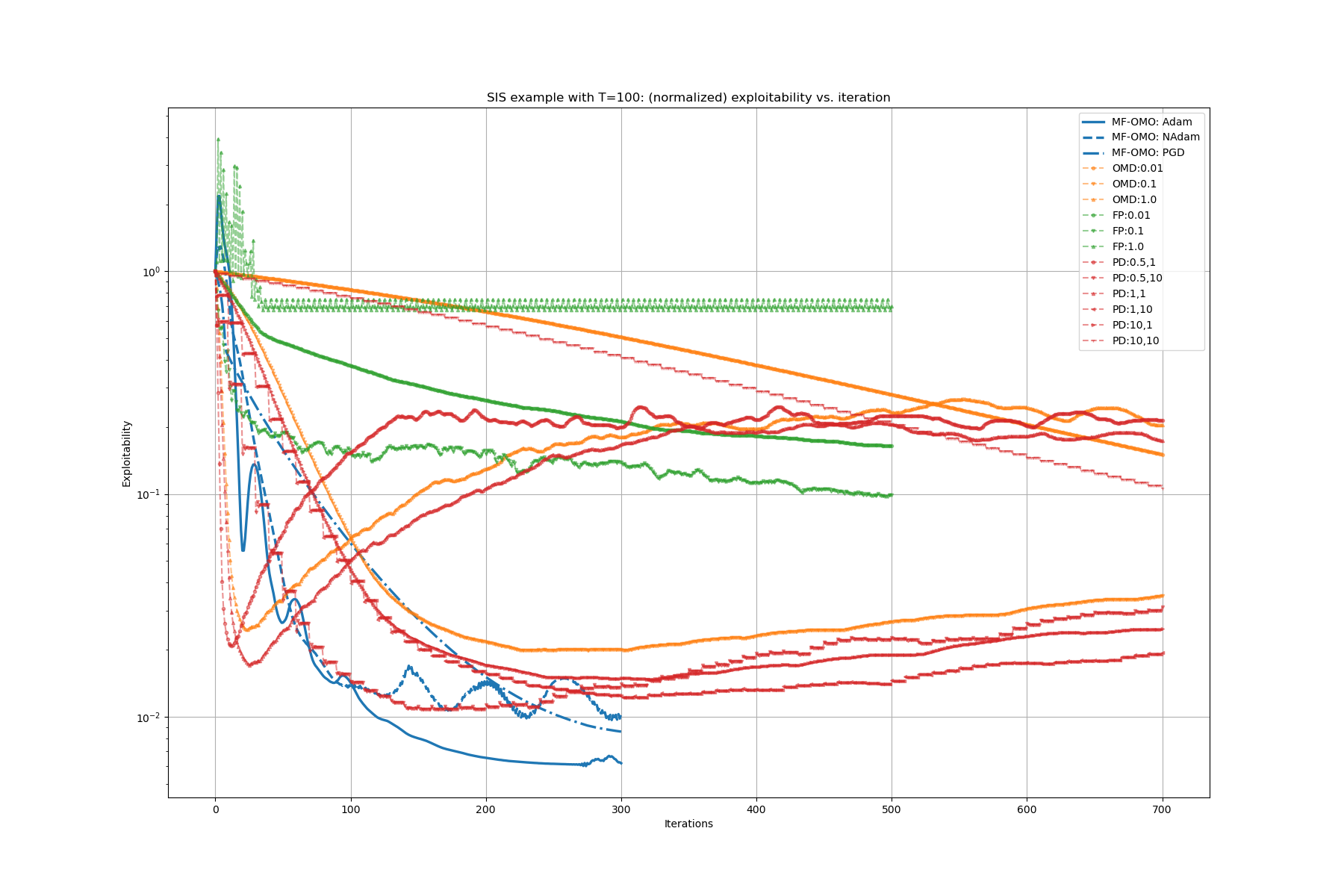}
	\caption{Comparison of different algorithms on SIS model with $T=100$: convergence against iterations}
	\label{fig:sis_iter_100}
\end{figure}

\begin{figure}[h]
	\centering
	\includegraphics[width=0.9\textwidth]{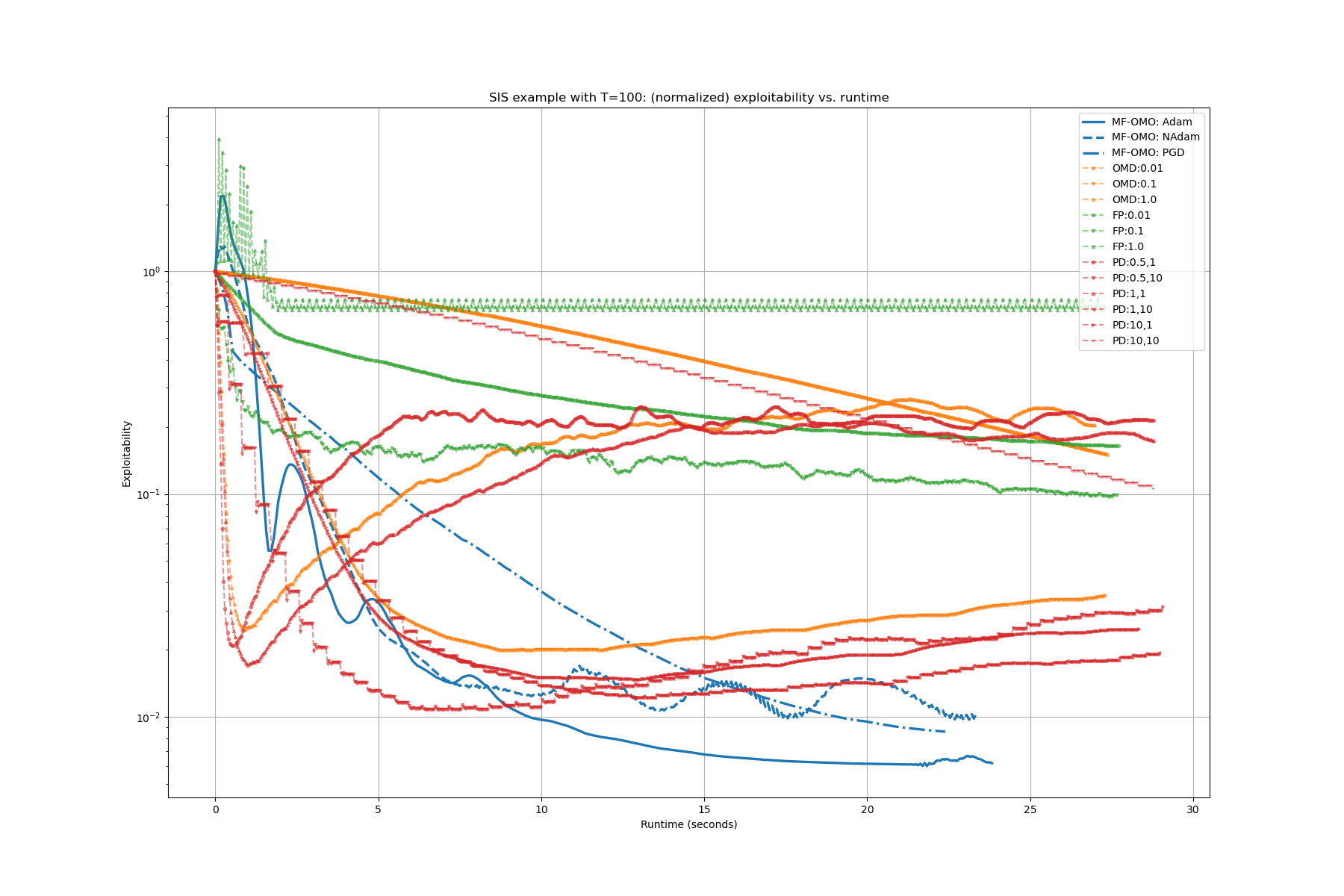}
	\caption{Comparison of different algorithms on SIS model with $T=100$: convergence against runtime}
	\label{fig:sis_time_100}
\end{figure}

We compare PGD, Adam, and Nesterov-accelerated Adam (NAdam), applied to the MF-OMO formulation, with online mirror descent (OMD), fictitious play (FP), and Prior Descent (PD), using several hyperparameter choices. Here, PD denotes the Prior Descent algorithm implemented in MFGLib, based on the prior-update approach of~\citep{cui2021approximately}. It is distinct from the policy-evaluation-based BPI variant in Section~\ref{sec:algo-policy-iteration}.
Figures~\ref{fig:sis_iter_50}, \ref{fig:sis_time_50},
\ref{fig:sis_iter_100}, and \ref{fig:sis_time_100}
compare these algorithms on this SIS model with $T=50$ and $T=100$, respectively.

The figure legends report the step sizes or learning rates used for OMD and FP and the hyperparameter settings used for PD. To make the runtime comparisons more balanced, the maximum number of iterations is chosen separately for each algorithm so that the total runtimes are comparable.
Starting from the same uniform initialization, MF-OMO with Adam gives the fastest decrease of normalized exploitability among the tested methods. The normalized exploitability of MF-OMO Adam quickly drops below $10^{-2}$ both as a function of the number of iterations and as a function of runtime. For $T=50$, only MF-OMO with Adam crosses this threshold within the displayed budgets. For $T=100$, the PGD and NAdam variants of MF-OMO also cross it. 

\section{Finding multiple NEs}
\label{sec:finding-multiple-nes}

In this section, we revisit the example with multiple Nash equilibria and illustrate how an optimization-based algorithm can converge to different equilibria from different initializations.

The discrete-time MFG is the one from Section~\ref{Example: multipleNE}, specialized to $n=5$. Take $\mathcal{S}=\mathcal{A}=\{1,\dots,n\}$, $T=10$, and an arbitrary initial distribution $\mu_0$. At time $0$, for any $\nu_0\in\Delta(\mathcal{S}\times\mathcal{A})$ and any $i,j=1,\dots,n$, the rewards are $r_0(i,j,\nu_0)=0$, and the transitions are deterministic:
\[
	P_0(i'|i,j,\nu_0)=\mathbf{1}_{\{i'=j\}}.
\]
For $t=1,\dots,T$, the rewards, and for $t=1,\dots,T-1$, the transition probabilities, are defined as in Section~\ref{Example: multipleNE}. In particular, if the population is concentrated on $(i,i)$, \ie, if $\nu_t=\delta_{(i,i)}$, then $r_t(i,j,\nu_t)=r^i\mathbf{1}_{\{i=j\}}$, and for $t<T$ the transition probabilities are deterministic, $P_t(i'|i,j,\nu_t)=\mathbf{1}_{\{i'=j\}}$.
We set $r^1=r^2=r^3=1.5$, and sample $r^4$, $r^5$, and $\{C_t\}_{t=1}^{T-1}$ independently and uniformly in $[0,1]$. With this choice, states $1$, $2$, and $3$ are maximizers of $i\mapsto r^i$, so the construction in Section~\ref{Example: multipleNE} gives at least three Nash equilibria, corresponding to the population remaining in state $1$, $2$, or $3$.

In the numerical experiment, let $\mathrm{RI}(i,\delta)$ denote a set of random initializations sampled from a $\delta$-neighborhood of the $i$-th equilibrium, for $i=1,2,3$. We test several such sets, with $20$ independent samples in each $\mathrm{RI}(i,\delta)$. We use the NAdam algorithm applied to the MF-OMO formulation and report the convergence behavior in Table~\ref{table:diff init}. For each neighborhood size $\delta$ and target equilibrium $i$, $p_0$ is the proportion of samples whose normalized exploitability remains above $10^{-3}$ after $400$ iterations, while $p_i$ is the proportion of samples that converge to an approximate NE closest, in Euclidean distance between mean-field flows, to the $i$-th equilibrium. The remaining column groups the proportions converging to one of the other two equilibria.

\begin{table}
	\centering
	\setlength{\tabcolsep}{5pt}
	\begin{tabular}{c|ccc|ccc|ccc}%
		\hline\hline
		\multirow{2}{*}{$\delta$} &
		\multicolumn{3}{c|}{NE 1} &
		\multicolumn{3}{c|}{NE 2} &
		\multicolumn{3}{c}{NE 3}                                                                                      \\
		                          & $p_1$ & $p_2+p_3$ & $p_0$ & $p_2$ & $p_1+p_3$ & $p_0$ & $p_3$ & $p_1+p_2$ & $p_0$ \\
		\hline\hline
		0.05                      & 70\%  & 30\%      & 0\%   & 75\%  & 20\%      & 5\%   & 65\%  & 35\%      & 0\%   \\
		\hline
		0.2                       & 50\%  & 50\%      & 0\%   & 65\%  & 30\%      & 5\%   & 60\%  & 35\%      & 5\%   \\
		\hline
		0.8                       & 45\%  & 50\%      & 5\%   & 40\%  & 55\%      & 5\%   & 30\%  & 70\%      & 0\%   \\
		\hline
		1.0                       & 40\%  & 60\%      & 0\%   & 45\%  & 55\%      & 0\%   & 25\%  & 75\%      & 0\%   \\
		\hline\hline
	\end{tabular}
	\caption{Convergence behavior of different initializations}
	\label{table:diff init}
\end{table}

Table~\ref{table:diff init} shows that most samples converge to an approximate NE with normalized exploitability below the target $10^{-3}$ tolerance for all initialization neighborhoods reported in the table. Different initializations can lead to different equilibria. At the smallest reported radius, $\delta=0.05$, a majority of runs converge to an approximate NE closest to the targeted equilibrium. This proportion decreases as the reported radius increases.

\section{Bibliographic notes}

The traffic-routing example in Section~\ref{sec:appli-traffic} is based on the routing MFG studied in~\citep{cabannes2022solving}. We use a discrete-time implementation derived from the \texttt{OpenSpiel} routing code~\citep{openspiel-routing-url}, with the fixed-horizon and payoff conventions specified in Section~\ref{sec:appli-traffic}. For a detailed discussion of realistic travel-time functions based on empirical studies, we refer to~\citep{cabannes2022solving}. More broadly, \texttt{OpenSpiel} is an open-source library of environments and algorithms for learning in multi-agent problems; it contains a variety of games, including examples of MFGs~\citep{lanctot2020openspiel}.

The numerical study for the SIS model is supported by MFGLib~\citep{guo2023mfglib}, with all parameters taken from~\citep{cui2021approximately}. The Adam algorithm follows~\citep{kingma2014adam}, and NAdam is adapted from~\citep{dozat2016incorporating}.
MFGLib is an open-source Python library designed to solve Nash equilibria (NEs) for general mean-field games (MFGs) through a user-friendly and customizable interface.
MFGLib provides a modular and extensible framework for the community to prototype and implement new algorithms, environments, and generalizations of MFGs. The multiple-equilibrium experiment in Section~\ref{sec:finding-multiple-nes} follows~\citep{guo2024mf}.

A related numerical study is~\citep{Perrin-2020-MFGflock}, which considers a flocking MFG and reports different learned (approximate) equilibria for different initializations. For more numerical examples and comparisons of learning algorithms, we refer to the survey~\citep{lauriere2022learningsurvey}.

\chapter{Conclusion and Outlook}
This monograph has developed a finite-dimensional perspective on discrete-time mean field games, in which the two defining requirements of equilibrium, optimality for a representative player and consistency of the population flow, can be studied through fixed-point arguments, learning dynamics, and optimization methods. Starting from the occupation-measure formulation of Markov decision processes, the primal-dual viewpoint leads naturally to MF-OMO and provides a common language for characterizing, computing, and comparing mean field Nash equilibria, including in models where equilibria are not unique.

	The chapters mirror this progression. Chapter~\ref{chap:mdp-markov-games} reviews Markov decision processes, their linear-programming formulation and duality, and finite-population Markov games. Chapter~\ref{chap:discrete-mfg} introduces the mean field framework, the solution concept of mean field Nash equilibrium, its links with finite-agent models, and the subclass of potential MFGs. Chapter~\ref{chap:analytical-MFG} develops existence, uniqueness under monotonicity, and primal--dual characterizations of the set of equilibria, including in models with multiple equilibria. Chapter~\ref{chap:algorithms} presents fixed-point-based and optimization-based algorithms, with a contraction result for damped fixed-point iterations and a local convergence theorem for PGD applied to MF-OMO. It also discusses convergence mechanisms for fictitious play, policy iteration, and online mirror descent.
	Chapter~\ref{chap:applications} illustrates the approach on crowd motion, traffic routing, SIS epidemics, and the computation of multiple equilibria. Readers wishing to go beyond the discrete setting can consult, e.g., \citet{carmona2018probabilistic,achdou2021mean} for the continuous-time theory and \citet{lauriere2022learningsurvey} for a survey of learning methods; the bibliographic notes at the end of each chapter provide further entry points.

	Table~\ref{tab:mfg-optimization-taxonomy} summarizes the main technical structures for discrete-time MFGs discussed in the manuscript. The rows are not mutually exclusive: an MFG may be both potential and monotone, while contractivity is a property of a selected response or update map. The last row records the general primal--dual characterization available in the finite-state setting without assuming any of these special structures. Across methods, exploitability, defined in~\eqref{eq:def-exploitability}, provides a measure of the incentive to deviate from equilibrium; the numerical experiments of Chapter~\ref{chap:applications} use it to compare algorithms across models, including cases where no convergence guarantee is available.

	\begin{table}[t]
		\centering
		\caption{Structural and optimization viewpoints for discrete-time MFGs.}
		\label{tab:mfg-optimization-taxonomy}
		\scriptsize
		\setlength{\tabcolsep}{2.4pt}
		\renewcommand{\arraystretch}{1.06}
\begin{tabularx}{\textwidth}{
  >{\raggedright\arraybackslash}p{0.14\textwidth}
  >{\raggedright\arraybackslash}p{0.20\textwidth}
  >{\raggedright\arraybackslash}X
  >{\raggedright\arraybackslash}X
  >{\raggedright\arraybackslash}p{0.155\textwidth}}
			\toprule
			\textbf{Structure}
			 & \textbf{Mathematical object}
			 & \textbf{Representative methods}
			 & \textbf{Guarantee and limitation}
			 & \textbf{Location in the monograph}                                                                                                                                                                                                                                          \\
			\midrule

			\textbf{Contractive response}
			 & Fixed point
			for single-valued best-response or regularized response.
			 & Pure fixed-point iteration, softmax response, and damping.
			 & Unique fixed point and global convergence when the update map is a contraction; the condition typically requires a short horizon or weak interactions, and a single-valued response (selection or regularization).
			 & Sections~\ref{sec:MFG} and~\ref{sec:existence-mfne} (Theorem~\ref{thm:existenceMFNE-Banach}); Section~\ref{sec:algo-fixed-point} (Proposition~\ref{prop:damped-monotone-residual}).
			\\
			\addlinespace

			\textbf{Potential MFG}
			 & Scalar potential over feasible occupation-measure or population flows; its first variation reproduces the coupling incentives.
			 & Fictitious play as a Frank--Wolfe-type method, conditional-gradient methods, and smoothed policy iteration.
			 & Reduction of MFNE computation to scalar optimization; convergence and exploitability guarantees under regularity assumptions; no universal guarantee for each update rule.
			 & Section~\ref{sec:potential-mfg}; Sections~\ref{sec:algo-fictitious-play} and~\ref{sec:algo-policy-iteration}; bibliographic notes in Section~\ref{sec:algo-bibliographic-notes}.                                                                                            \\
			\addlinespace

			\textbf{Monotone MFG}
			 & Lasry--Lions monotonicity: mean-field-independent transitions and monotone reward coupling.
			 & Fictitious play with averaged best responses; online mirror descent.
			 & Uniqueness of the equilibrium mean-field flow; convergence of fictitious play and online mirror descent under additional regularity assumptions.
			 & Section~\ref{sec:uniqueness-mfne} (Theorem~\ref{thm:unique-mfne-monotonicity}); Sections~\ref{sec:algo-fictitious-play} and~\ref{sec:algo-omd}; bibliographic notes in Section~\ref{sec:algo-bibliographic-notes}.                                                          \\
			\addlinespace

			\textbf{General finite-state MFG}
			 & Occupation-measure primal--dual conditions with KKT complementarity; no potentiality, monotonicity, or contractivity needed.
			 & Complementarity or feasibility solvers; MF-OMO with projected gradient descent; Adam/NAdam; multiple initializations for distinct equilibria.
			 & Objective value zero is equivalent to an MFNE; generally nonconvex; local PGD theory gives local convergence, while Theorem~\ref{MFNE-vs-VMO} gives a one-sided implication from small exploitability to small MF-OMO objective.
			 & Section~\ref{sec:LP-MDP}; Sections~\ref{sec:mfne-characterization} and~\ref{sec:maths-MFG} (Theorems~\ref{thm:nash-opt}, \ref{thm:mfvmo-full} and~\ref{MFNE-vs-VMO}); Section~\ref{sec:algo-optimization}; Sections~\ref{sec:appli-sis} and~\ref{sec:finding-multiple-nes}. \\
			\bottomrule
		\end{tabularx}
	\end{table}

	The finite-state and finite-horizon framework makes these connections explicit, while pointing toward several open directions: understanding the global optimization landscape of MF-OMO beyond contractive, monotone, and potential models; developing mechanisms for equilibrium selection when equilibria are multiple; sharpening the finite-population approximation guarantees behind Theorem~\ref{thm:approx-Nash-mf} into finite-sample and finite-iteration bounds; and relaxing perfect homogeneity and anonymity. Common noise, partial observations, large or continuous state--action spaces, and the design of model-free, sample-efficient algorithms raise further challenges for which the finite-state toolbox developed here may serve as a template.

	A natural continuation is the continuous-time counterpart of the occupation-measure approach. Recent work by~\citet{Guo2025Primaldual} formulates the continuous-time representative-agent control problem as an infinite-dimensional linear program over occupation measures and identifies a dual problem involving subsolutions of the associated Hamilton--Jacobi--Bellman equation. Strong duality then yields a characterization of the full set of mean field Nash equilibria without requiring convexity of the Hamiltonian, uniqueness of its optimizer, or full support of the equilibrium flow. This development opens the way to computational continuous-time counterparts of MF-OMO, to discretization and approximation results that rigorously connect the discrete- and continuous-time formulations, and to scalable primal-dual algorithms that combine ideas from stochastic control, partial differential equations, optimization, and learning.

	A second set of open directions concerns the relaxation of perfect
	homogeneity and anonymity. Multi-population and personalized MFGs allow
	agents to have type-dependent dynamics, rewards, and other primitives,
	while graphon games and graphon MFGs replace the single aggregate
	population flow by a type-dependent network aggregate; see, for instance,
	\citep{parise2023graphon,caines2021graphon,carmona2022stochasticgraphonstatic}.
	In a graphon model, each agent is identified with a type, and the population that an agent interacts with is a weighted aggregate of the population flows of the different types, with weights describing the network, rather than one common distribution shared by all agents. This suggests extending the occupation-measure and primal--dual formulations from a single flow to a collection of type-indexed flows, and developing finite-type and graphon approximations with quantitative equilibrium guarantees for large network games. Learning algorithms for such graphon mean-field models are an active direction; see, for example,
	\citep{cuikoeppl2022}. For a survey of extensions of MFGs beyond perfect homogeneity and anonymity, see~\citep{lauriere2025overview}.

	A complementary direction, developed in recent work on finite-agent
	games and their mean-field limits, is to relax exact potentiality. In an
	$\alpha$-potential game, variations in individual objectives under
	unilateral deviations are represented by variations of a scalar
	potential up to an error $\alpha$. A global optimizer of this approximate potential therefore yields an approximate Nash equilibrium, with $\alpha$ and any global optimality gap adding to the approximation error; the magnitude of $\alpha$ quantifies the departure from exact potential structure caused by finite population size, heterogeneity, or asymmetric interactions~\citep{guolizhang2025,guoliMSW2025}.
	Recent limit theory showing that potential MFGs arise from
	$N$-player $\alpha_N$-potential games when $\alpha_N\to0$ provides a
	natural bridge between finite-agent dynamic games, potential MFGs,
	and mean-field control
	\citep{guowangzhang2026}. Developing learning and
	primal-dual methods along this hierarchy would extend the optimization perspective
	developed here.

	Across these directions, a common goal is to obtain existence, uniqueness, convergence, and approximation guarantees while relaxing the structural assumptions, the setting, or both. A unifying theme of the monograph is that in finite state and action spaces, equilibria are characterized by finite-dimensional optimization problems and their optimality conditions, in contrast with the infinite-dimensional equations arising in continuous time, so that analysis, computation, and learning rest on a common set of tools. Notions such as fixed-point formulations, monotonicity, potential structure, exploitability, and complementarity recur throughout the mean field literature, including settings beyond the scope of these notes. We hope this discrete toolbox serves the reader as an entry point into the broader literature.

\begin{acknowledgements}

	The authors would like to thank the anonymous reviewers for their helpful comments.
	X.G. would like to thank all her collaborators, especially Anran Hu and Junzi Zhang for  their  contributions to the work presented in this monograph.
	M.L. would like to thank his collaborators for all the works cited in this monograph and for all the constructive comments on preliminary versions of the monograph.
\end{acknowledgements}

\printindex
\printindex[notation]

\appendix

\chapter{Auxiliary results and proofs}

\section[Proof of Theorem~\ref*{thm:existence-NE}]{Proof of Theorem~\ref{thm:existence-NE}}
\label{app:proof-existence-NE}

\begin{proof}
	Since the state and action spaces are finite and the horizon is finite, the result follows from backward induction. For a state-action profile $(\bfs,\bfa)$, let $\nu^N(\bfs,\bfa)$ denote the empirical distribution of the pairs $(s^j,a^j)$, $j=1,\dots,N$. At time $T$, for each state profile $\bfs\in\states^N$, consider the finite normal-form game in which player $i$ chooses $a^i\in\actions$ and receives payoff
	\[
		r^i_T(s^i,a^i,\nu^N(\bfs,\bfa)).
	\]
	By Nash's theorem, this finite game admits a mixed Nash equilibrium. Choose one such equilibrium for each $\bfs$ and denote the corresponding value for player $i$ by $V^i_T(\bfs)$.

	Suppose now that continuation values $V^i_{t+1}$ have been constructed. For each state profile $\bfs\in\states^N$, define a finite normal-form game with action profile $\bfa=(a^1,\dots,a^N)$ and payoff
	\[
		r^i_t(s^i,a^i,\nu^N(\bfs,\bfa))
		+
		\sum_{\bfs'\in\states^N}
		\left(\prod_{j=1}^N P_t(s'^j|s^j,a^j,\nu^N(\bfs,\bfa))\right)
		V^i_{t+1}(\bfs').
	\]
	Again, Nash's theorem gives a mixed Nash equilibrium for this one-stage game. Choosing one such equilibrium for every $\bfs$ defines the Markov policy at time $t$ and the value $V^i_t(\bfs)$.

	Iterating this construction backward from $T$ to $0$ defines a Markov policy profile. By the one-stage Nash equilibrium property at every time and state profile, no player can improve her expected total reward by a unilateral Markov deviation. Hence the constructed profile is a Nash equilibrium.
\end{proof}

\section{Proof of approximate Nash property}

Here we present a proof of the fact that mean field equilibrium policies provide approximate equilibria in finite-population games, i.e., Theorem~\ref{thm:approx-Nash-mf}, which assumes that $P$ is independent of $\nu$. In this situation, the agents' dynamics are not coupled. We refer to the literature for more details.
\begin{proof}[Proof of Theorem~\ref{thm:approx-Nash-mf}]
	We start with a sketch of proof. Let $\hat\pi$ be a mean field Nash equilibrium policy. Denote by $\hat\bspi = (\hat\pi,\dots,\hat\pi)$ the $N$-player policy profile consisting of all $\hat\pi$. Fix $\epsilon>0$ and let $\eta=\epsilon/2$. Let $\pi \in \Pi$ be a decentralized policy.

	To distinguish between the finite-population game objective and the mean field objective, we denote them respectively by $J^N$ and $J$.  We will justify each of the following inequalities:
	\begin{align*}
		J^N(\pi,\hat\bspi^{-i})
		 & \le J(\pi;\nu^{\hat\pi}) + \eta
		\\
		 & \le J(\hat\pi;\nu^{\hat\pi}) + \eta
		\\
		 & \le J^N(\hat\pi,\hat\bspi^{-i}) + 2\eta
		= J^N(\hat\pi,\hat\bspi^{-i}) + \epsilon.
	\end{align*}
	The second inequality follows from the fact that $\hat\pi$ is a best response against $\nu^{\hat\pi}$, since $\hat\pi$ is a mean field Nash equilibrium. The first and third inequalities will be obtained by showing the following more general property: for every fixed $\pi' \in \strategies$ and every $\eta>0$, there exists $N_0$ large enough such that for all $N \ge N_0$,
	\[
		\sup_{\pi\in\strategies}|J^N(\pi, \underline\pi') - J(\pi;\nu^{\pi'})| \le \eta,
	\]
	where $\underline\pi' = (\pi',\dots,\pi')$ denotes the vector of $N-1$ copies of $\pi'$.

		We now prove this estimate.
		Since the players are interchangeable, it is enough to consider a deviation by the first player. Let the first player use $\pi$ while the other players use $\underline\pi'$. We denote by
		\[
			\nu^{N,\pi,\underline\pi'}_t
			=
			\frac{1}{N}\sum_{i=1}^N \indic_{(s^i_t,a^i_t)}
		\]
		the induced empirical population distribution at time $t$. Since $P$ is independent of the distribution, the trajectory of the deviating player is the same as in the mean field problem with flow $\nu^{\pi'}$, and the variables $(s^i_t,a^i_t)$, $i=2,\dots,N$, are i.i.d. with law $\nu^{\pi'}_t$. If
		\[
			\overline\nu^{N-1,\pi'}_t
			=
			\frac{1}{N-1}\sum_{i=2}^N \indic_{(s^i_t,a^i_t)},
		\]
		then
		\[
			\left\|\nu^{N,\pi,\underline\pi'}_t-\nu^{\pi'}_t\right\|_\infty
			\le
			\frac{1}{N}
			+
			\left\|\overline\nu^{N-1,\pi'}_t-\nu^{\pi'}_t\right\|_\infty.
		\]
		The right-hand side does not depend on $\pi$ and converges to zero in probability by the law of large numbers; since it is bounded, it also converges to zero in expectation. Since $\states$, $\actions$, and $\timesteps$ are finite and $r_t(s,a,\cdot)$ is uniformly continuous and bounded on the compact simplex, the expected reward difference is controlled by the modulus of continuity of the reward. Hence
		\[
			\sup_{\pi\in\strategies}|J^N(\pi, \underline\pi') - J(\pi;\nu^{\pi'})| \to 0.
		\]
		This proves the claimed estimate. Applying it with $\pi'=\hat\pi$, together with the optimality of $\hat\pi$ against $\nu^{\hat\pi}$, yields the result for the first player, and symmetry gives the same conclusion for every player.
\end{proof}

\section{Existence of MFNE}

\subsection{Proof of Theorem~\ref{thm:existenceMFNE-Banach}}
\label{Appendix: MFNEexistence-Banach}
We start with a useful preliminary result.
\begin{lemma}
	[Lipschitz continuity of $\EMF$]
	\label{lem:EMF-Lip}
	Assume that $P$ is Lipschitz continuous in the mean field argument with constant $L_P$, for the sup norm. Then $\EMF$ is Lipschitz continuous with Lipschitz constant
	\[
		L_{\EMF}:=\sum_{k=0}^{T}(L_P+|\states||\actions|)^k.
	\]
\end{lemma}

\begin{proof}
		Let $\pi,\pi'$ be two policies. Let $\mu^\pi=\MF(\pi)$ and $\mu^{\pi'}=\MF(\pi')$, and let $\nu^\pi=\EMF(\pi)$ and $\nu^{\pi'}=\EMF(\pi')$. Set
		\[
			B:=L_P+|\states||\actions|.
		\]
		We prove by induction that, for every $t\in\timesteps$,
		\[
			\|\nu_t^\pi-\nu_t^{\pi'}\|_\infty
			\le C_t\|\pi-\pi'\|_\infty,
			\qquad
			C_t:=\sum_{k=0}^{t}B^k.
		\]
		At time $0$, since $\mu_0^\pi=\mu_0^{\pi'}=\mu_0$,
		\[
			|\nu_0^\pi(s,a)-\nu_0^{\pi'}(s,a)|
			=
			\mu_0(s)|\pi_0(a|s)-\pi'_0(a|s)|
			\le \|\pi-\pi'\|_\infty.
		\]
		Hence the claim holds at time $0$, with $C_0=1$.

		Assume that the claim holds at time $t\in\{0,\dots,T-1\}$. For every $s'\in\states$,
		\begin{align*}
			 & |\mu_{t+1}^{\pi}(s')-\mu_{t+1}^{\pi'}(s')|
			\\
			 & =
			\left|
			\sum_{s\in\states}\sum_{a\in\actions}
			\left[
				P_t(s'|s,a,\nu_t^\pi)\nu_t^\pi(s,a)
				-
				P_t(s'|s,a,\nu_t^{\pi'})\nu_t^{\pi'}(s,a)
				\right]
			\right|
			\\
			 & \le
			\sum_{s\in\states}\sum_{a\in\actions}
			|P_t(s'|s,a,\nu_t^\pi)-P_t(s'|s,a,\nu_t^{\pi'})|\nu_t^\pi(s,a)
			\\
			 & \quad+
			\sum_{s\in\states}\sum_{a\in\actions}
			P_t(s'|s,a,\nu_t^{\pi'})
			|\nu_t^\pi(s,a)-\nu_t^{\pi'}(s,a)|
			\\
			 & \le
			\left(L_P+|\states||\actions|\right)
			\|\nu_t^\pi-\nu_t^{\pi'}\|_\infty
			\\
			 & \le
			B C_t \|\pi-\pi'\|_\infty.
		\end{align*}
		Therefore, for every $(s',a')\in\states\times\actions$,
\begin{align*}
  & |\nu_{t+1}^{\pi}(s',a')-\nu_{t+1}^{\pi'}(s',a')|\\
  &= |\mu_{t+1}^{\pi}(s')\pi_{t+1}(a'|s')
     -\mu_{t+1}^{\pi'}(s')\pi'_{t+1}(a'|s')|\\
  &\le |\mu_{t+1}^{\pi}(s')-\mu_{t+1}^{\pi'}(s')|
     +|\pi_{t+1}(a'|s')-\pi'_{t+1}(a'|s')|\\
  &\le (BC_t+1)\|\pi-\pi'\|_\infty.
\end{align*}
		Since $C_{t+1}=BC_t+1$, the induction is complete. Taking the maximum over $t\in\timesteps$ gives
\[
  \begin{aligned}
  \|\EMF(\pi)-\EMF(\pi')\|_\infty
  &\le C_T\|\pi-\pi'\|_\infty\\
  &=\left(\sum_{k=0}^{T}(L_P+|\states||\actions|)^k\right)\|\pi-\pi'\|_\infty.
  \end{aligned}
\]
\end{proof}

We now turn to the proof of existence of an MFNE using the Banach fixed-point theorem.
\begin{proof}[Proof of Theorem~\ref{thm:existenceMFNE-Banach}]
	We apply the Banach fixed-point theorem (see, for instance, \citep[Theorem 5.7]{brezis2011functional}) to $\Phi =  \EMF \circ \BR$.

	By assumption, for every mean field, there is a unique best response, so $\Phi$ is a function from $\Delta(\states\times\actions)^\timesteps$ to itself.
	We note that, since $\states$ and $\actions$ are finite, $\Delta(\states\times\actions)^\timesteps$ is a complete metric space when endowed (for instance) with the distance induced by the sup norm, denoted by $\|\nu\|_\infty = \max_{t,s,a} |\nu_t(s,a)|$.

	By assumption, $\BR$ is Lipschitz continuous with Lipschitz constant $L_{\BR}$.
	By Lemma~\ref{lem:EMF-Lip}, $\EMF$ is Lipschitz continuous with Lipschitz constant $L_{\EMF}$.

	Then $\Phi$ is Lipschitz continuous with Lipschitz constant $L_{\EMF} L_{\BR}$, which is assumed to be smaller than $1$. Thus $\Phi$ is a strict contraction and, by the Banach fixed-point theorem, it has a unique fixed point $\hat\nu$. Set $\hat\pi=\BR(\hat\nu)$. Then $\hat\nu=\EMF(\hat\pi)$, so $(\hat\pi,\hat\nu)$ is an MFNE. Conversely, any MFNE $(\pi,\nu)$ satisfies $\pi=\BR(\nu)$ and $\nu=\Phi(\nu)$, hence $\nu=\hat\nu$ and, by uniqueness of the best response, $\pi=\hat\pi$.
\end{proof}

\subsection{Proof of Theorem~\ref{thm:existenceMFNE-Kakutani}}
\label{Appendix: MFNEexistence-Kakutani}

\begin{proof}[Proof of Theorem~\ref{thm:existenceMFNE-Kakutani}]
		Let $\Gamma^{\rm g}: \strategies \to 2^\strategies$ be the correspondence defined by
		\[
			\Gamma^{\rm g}(\pi)=\mathrm{Greedy}(Q^{*,\nu^\pi}),
		\]
		where $Q^{*,\nu^\pi}$ is the optimal Q-function of the MDP faced by a representative agent when the mean field is fixed to $\nu^\pi=\EMF(\pi)$.

		We will apply Kakutani's fixed point theorem to $\Gamma^{\rm g}$. This is sufficient for existence of an MFNE: if $\hat\pi\in\Gamma^{\rm g}(\hat\pi)$ and $\hat\nu=\nu^{\hat\pi}$, then $\hat\pi$ is optimal for the representative-agent MDP with mean field $\hat\nu$, and $\hat\nu=\EMF(\hat\pi)$ by definition. Hence $(\hat\pi,\hat\nu)$ is an MFNE.

		First, $\strategies$ is non-empty, compact, and convex, because it is a finite product of probability simplexes.

		Second, for every $\pi\in\strategies$, $\Gamma^{\rm g}(\pi)$ is non-empty. Indeed, for every $(t,s)$, the set $\argmax Q^{*,\nu^\pi}_t(s,\cdot)$ is non-empty because $\actions$ is finite.

		Third, for every $\pi\in\strategies$, $\Gamma^{\rm g}(\pi)$ is convex. If $\pi^1,\pi^2\in\Gamma^{\rm g}(\pi)$ and $\lambda\in[0,1]$, define $\pi^\lambda=\lambda\pi^1+(1-\lambda)\pi^2$. If $\pi^\lambda_t(a|s)>0$, then either $\pi^1_t(a|s)>0$ or $\pi^2_t(a|s)>0$, and therefore $a\in\argmax Q^{*,\nu^\pi}_t(s,\cdot)$. Hence $\pi^\lambda\in\Gamma^{\rm g}(\pi)$.

		Last, $\Gamma^{\rm g}$ has a closed graph. Let $(\pi_n,\tilde\pi_n)\to(\pi,\tilde\pi)$ with $\tilde\pi_n\in\Gamma^{\rm g}(\pi_n)$. Suppose that $\tilde\pi_t(a|s)>0$. Then, for all $n$ large enough, $\tilde\pi_{n,t}(a|s)>0$, so
		\[
			a\in\argmax Q^{*,\nu^{\pi_n}}_t(s,\cdot).
		\]
		The map $\pi\mapsto\nu^\pi$ is continuous in the finite-state finite-horizon setting. Since $P_t(s'|s,a,\nu)$ and $r_t(s,a,\nu)$ are continuous in $\nu$ and $\actions$ is finite, the backward dynamic-programming equations imply, by induction on $t$, that $\nu\mapsto Q^{*,\nu}_t(s,a)$ is continuous for every $(t,s,a)$. Passing to the limit, $a\in\argmax Q^{*,\nu^\pi}_t(s,\cdot)$. Thus $\tilde\pi\in\Gamma^{\rm g}(\pi)$.

		The assumptions of Kakutani's theorem are therefore satisfied, so $\Gamma^{\rm g}$ admits a fixed point. As observed above, this fixed point yields an MFNE.
\end{proof}

\section{Uniqueness of MFNE}
\label{app:auxi-uniqueness-mfne}

\begin{proof}[Proof of Theorem~\ref{thm:unique-mfne-monotonicity}]
	For the sake of contradiction, assume that there are two MFNEs, say $(\hat\pi,\hat\nu)$ and $(\hat\pi',\hat\nu')$, with different mean-field flows: $\hat\nu\neq\hat\nu'$. By optimality of $\hat\pi$ against mean field $\hat\nu$, we have:
	\begin{align*}
		0 & \le J(\hat\pi; \hat\nu) - J(\hat\pi'; \hat\nu)
		\\
		  & = \EE\left[ \sum_{t=0}^T r_t(s^{\hat\pi,\hat\nu}_t, a^{\hat\pi,\hat\nu}_t, \hat\nu_t) - r_t(s^{\hat\pi',\hat\nu}_t, a^{\hat\pi',\hat\nu}_t, \hat\nu_t)\right]
		\\
		  & =  \sum_{t=0}^T \sum_{(s,a)}\left[r_t(s, a, \hat\nu_t) \nu^{\hat\pi,\hat\nu}_t(s,a) - r_t(s, a, \hat\nu_t) \nu^{\hat\pi',\hat\nu}_t(s,a) \right].
	\end{align*}
	Similarly,
	\begin{align*}
		0 & \le J(\hat\pi'; \hat\nu') - J(\hat\pi; \hat\nu')
		\\
		  & =  \sum_{t=0}^T \sum_{(s,a)}\left[r_t(s, a, \hat\nu_t') \nu^{\hat\pi',\hat\nu'}_t(s,a) - r_t(s, a, \hat\nu_t') \nu^{\hat\pi,\hat\nu'}_t(s,a) \right].
	\end{align*}
	By the consistency condition in the definition of MFNE (see Definition~\ref{def:mfNE}), $\nu^{\hat\pi,\hat\nu}_t = \hat\nu_t$ and $\nu^{\hat\pi',\hat\nu'}_t = \hat\nu'_t$. Furthermore, since $P$ does not depend on the mean field, $\nu^{\hat\pi,\hat\nu'}_t = \nu^{\hat\pi,\hat\nu}_t = \hat\nu_t$ and $\nu^{\hat\pi',\hat\nu}_t = \nu^{\hat\pi',\hat\nu'}_t = \hat\nu'_t$.

	Summing the terms and using the above equalities, we obtain:
	\begin{align*}
		0 & \le J(\hat\pi; \hat\nu) - J(\hat\pi'; \hat\nu) + J(\hat\pi'; \hat\nu') - J(\hat\pi; \hat\nu')
		\\
		  & =  \sum_{t=0}^T \sum_{(s,a)}\Big[r_t(s, a, \hat\nu_t) \hat\nu_t(s,a) - r_t(s, a, \hat\nu_t) \hat\nu'_t(s,a)
		\\
		  & \qquad\qquad + r_t(s, a, \hat\nu_t') \hat\nu'_t(s,a) - r_t(s, a, \hat\nu_t') \hat\nu_t(s,a) \Big]
		\\
		  & =  \sum_{t=0}^T \sum_{(s,a)}\Big[\Big(r_t(s, a, \hat\nu_t)  - r_t(s, a, \hat\nu'_t) \Big) \Big(\hat\nu_t(s,a) - \hat\nu'_t(s,a) \Big) \Big],
	\end{align*}
	which is strictly negative by monotonicity because $\hat\nu\neq\hat\nu'$. This contradiction shows that any two MFNEs must have the same mean-field flow. If, in addition, the best response to this flow is unique, then their policies also coincide.
\end{proof}

\backmatter  %

\printbibliography

\end{document}